\documentclass[11pt,letterpaper]{article}
\usepackage{amsmath,amssymb,amsthm,mathrsfs}
\usepackage{comment}
\usepackage[shortlabels]{enumitem}
\usepackage{amsrefs}
\usepackage{esvect}

\def\r{\mathbb{R}}
\def\s{\mathbb{S}}
\def\<{\langle}
\def\>{\rangle}
\def\norm#1#2{\|#2\|_{#1}}
\def\loc{\text{loc}}

\def\Vol{\text{Vol}}

\DeclareMathOperator{\spa}{span}
\DeclareMathOperator{\supp}{supp}
\DeclareMathOperator{\Proj}{Proj}
\DeclareMathOperator{\Id}{Id}
\DeclareMathOperator{\trace}{tr}
\DeclareMathOperator{\dist}{dist}
\DeclareMathOperator{\codim}{codim}

\newcommand{\RN}[1]{%
\textup{\uppercase\expandafter{\romannumeral#1}}%
}

\newtheorem{thm}{Theorem}[section]
\newtheorem{prop}[thm]{Proposition}
\newtheorem{coro}[thm]{Corollary}
\newtheorem{lem}[thm]{Lemma}
\newtheorem{rem}[thm]{Remark}

\newcommand{\Addresses}{{
  \bigskip
  \footnotesize

  L.~Gong, \textsc{Department of Mathematics, The Chinese University of Hong Kong, Hong Kong, and Department of Mathematics, Rutgers University, 110 Frelinghuysen Road, Piscataway, NJ, USA.}
  \par\nopagebreak
  \textit{E-mail address}: \texttt{lg689@rutgers.edu}

  \medskip

  Y.Y.~Li, \textsc{Department of Mathematics, Rutgers University, 110 Frelinghuysen Road, Piscataway, NJ, USA.}\par\nopagebreak
  \textit{E-mail address}: \texttt{yyli@math.rutgers.edu}

}}

\providecommand{\keywords}[1]
{
  \small	
  \textit{Key Words---} #1
}

\title{Conformal metrics of constant scalar curvature with unbounded volumes}

\author{Liuwei Gong and Yanyan Li\\
}
\date{}
\begin{document}

\maketitle

\begin{abstract}
For $n\geq 25$, we construct a smooth metric $\tilde{g}$ on the standard $n$-dimensional sphere $\s^n$ such that there exists a sequence of smooth
metrics $\{\tilde{g}_k\}_{k\in\mathbb{N}}$ 
conformal to $\tilde g$
where each $\tilde g_k$ has
 scalar curvature $R_{\tilde{g}_k}\equiv 1$ and
their volumes $\Vol(\s^n,\tilde{g}_k)$ tend to infinity
as $k$ approaches infinity.
\end{abstract}


\begin{center}
\end{center}

\tableofcontents

\setlength{\parskip}{1em}

\section{Introduction}

Let $(M,g)$ be a closed, connected, smooth 
Riemannian manifold of dimension $n$. 
For $n=2$, the 
uniformization theorem of Poincar\'e and Koebe
implies that there exists a constant Gaussian curvature 
 metric in the conformal class $[g]: = \{e^{2f}g \ |\  f\in C^{\infty}(M)\}$.

The Yamabe problem is a natural generalization of the uniformization theorem to higher dimensions. For $n\geq 3$, it is to find metrics with constant scalar curvature in the given conformal class $[g]$. 
For a conformal metric  $g_u:=u^{\frac{4}{n-2}}g$, where $0<u\in C^{\infty}(M)$, the scalar curvature  of $g_u$ satisfies
$$
R_{g_u}=-c(n)^{-1}u^{-\frac{n+2}{n-2}}L_gu,\quad L_gu:=\Delta_{g}u-c(n)R_{g}u,
$$
where $\Delta_{g}$ is the Laplace-Beltrami operator associated with $g$, $R_g$ is the scalar curvature of $g$, and $c(n)=\frac{n-2}{4(n-1)}$.

If a constant scalar curvature metric exists 
in the conformal class $[g]$, it is evident that  
the sign of the  scalar curvature  must be the same as that of $\lambda_1(-L_g)$, where $\lambda_1(-L_g)$ denotes 
the first eigenvalue of $-L_g$.

The Yamabe problem is therefore equivalent to finding a positive smooth solution of the following semilinear elliptic partial differential equation:
\begin{equation}\label{yamabe eqn}
-L_{g}u=\lambda_1(-L_g)u^{\frac{n+2}{n-2}}.
\end{equation} 

The Yamabe problem was solved through the work of Yamabe \cite{yamabe1960deformation}, Trudinger \cite{trudinger1968remarks}, Aubin \cite{aubin1976equations} and Schoen \cite{schoen1984conformal}.
They proved the existence of a positive solution of 
(\ref{yamabe eqn}) by finding a minimizer of the following Yamabe functional:
\begin{equation*}
E_g[u]=\left(\int_M u^{\frac{2n}{n-2}}dV_g\right)^{-\frac{n-2}{n}}\int_M uL_gudV_g, \qquad u\in H^1(M)\setminus\{0\}.
\end{equation*}

Consider the solution space
$$
\mathcal{M}_{g}:=\{ 0<u\in C^{\infty}(M)\mid u \text{ solves equation (\ref{yamabe eqn})}\}.
$$
It is straightforward 
to prove that the solution of (\ref{yamabe eqn}) is unique when 
$\lambda_1(-L_g)<0$, 
and is unique up to multiplication by a positive constant when
$\lambda_1(-L_g)=0$.
However, when $\lambda_1(-L_g)>0$, $\mathcal{M}_{g}$ may contain non-minimizing solutions. Schoen proved in  \cite{Schoen1989variational}
that  the product manifold $ S^{n-1}(1)\times S^1(L)$ has non-minimizing rotationally symmetric solutions for large   $L$.  Hebey and Vaugon proved in  \cite{hebey1993probleme} 
that when $(M,g)$ has a non-trivial isometry subgroup $G$, there exist  
non-minimizing 
 $G$-invariant solutions.
Pollack proved in \cite{pollack1993nonuniqueness}  that  the set of conformal classes $[g]$ containing non-minimizing solutions is dense in $C^0$ topology.

Consider  the conformal diffeomorphism group
$$
C(M,g):=\{\Phi:M\rightarrow M \text{ diffeomorphism }\mid \Phi^*g\in [g]\}.
$$
For  $\Phi \in C(M,g)$ and $u\in \mathcal{M}_g$, 
$(u\circ \Phi)|\Phi'|^{\frac{n-2}{2}} \in \mathcal{M}_g$
 since $\Phi^*g=|\Phi'|^{2}g$ where $|\Phi'|$ denotes the linear stretch factor of $\Phi$ with respect to $g$. 

For the standard sphere $\s^n$, since $\{|\Phi'|\mid \Phi\in C(\s^n,g_{\text{stand}})\}$ is noncompact in $C^0$ topology, the solution space  $\mathcal{M}_{g_{\text{stand}}}$ is also noncompact in $C^0$ topology. On the other hand, by \cite{ferrand1971transformations} and \cite{lelong1976geometrical}, if $(M,g)$ is not conformally diffeomorphic to the standard sphere, the conformal diffeomorphism group $C(M,g)$ is inessential, namely,  it can be reduced to a group of isometries by a conformal change of the metric $g$; see also \cite{obata1970conformal}, \cite{obata1971conjectures}, \cite{lafontaine1988theorem}, and \cite{schoen1995conformal}. Consequently, $\{|\Phi'|\mid \Phi\in C(M,g)\}$ is compact in $C^0$ topology if and only if $(M,g)$ is not conformally diffeomorphic to the standard sphere.

Schoen initiated the investigation of the compactness of the solution space $\mathcal{M}_g$ for a closed Riemmanian manifold $M$ with
$\lambda_1(-L_g)>0$ and not conformally diffeomorphic to 
the standard sphere. He proved in \cite{Schoen1991}  the compactness when the manifolds are locally conformally flat.  Moreover he   gave in  \cites{Schoenlecturenote,Schoen1989variational,schoen1990report} insightful investigations for the problem on general manifolds.

For general manifolds of dimension $n=3$, the compactness
 was proved by Li and Zhu \cite{li1999yamabe}.  For  $n=4$, 
 the compactness was proved as a combination
 of the  results by 
 Druet \cite{druet2003one} and Li and Zhang \cite{li2004harnack}, 
 with the $H^1$ bound given in \cite{li2004harnack}
 and the $L^\infty$ bound under the assumption of an $H^1$ bound given in 
 \cite{druet2003one}.
 Subsequently, the following compactness results were 
 proved by independent works of three research groups.
Druet \cite{druet2004compactness} proved the compactness for  $n\leq 5$; 
 Li and Zhang \cite{li2005compactness}
 proved the compactness if one of the following holds:   $n\leq 9$, or $n\ge 6$ with $|W|>0$ on $M$, or $n\geq 8$ with $|W|+|\nabla_g W|>0$ on $M$, where $W$ is the Weyl tensor of $g$; Marques \cite{marques2005priori} proved the compactness for 
 $n\leq 7$ as well as for    $n\geq 6$ with $|W|>0$ on $M$. 
 
 Li and Zhang \cite{li2007compactness}  continued to prove the compactness for $n\le 11$, as well as for $n=10, 11$ with
 $|W|+|\nabla_g W|+|\nabla^2_g W|>0$ on $M$.
 Aubin gave discussions in 
 \cite{aubin2006quelques} and \cite{aubin2007demonstration} on 
  the problem and related ones.

It was a surprise to researchers in the field when  Brendle  provided examples in \cite{brendle2008blow} demonstrating that $\mathcal{M}g$ is not necessarily bounded in $L^\infty(M)$ for dimensions $n\geq 52$. Brendle and Marques extended this result to $n\geq 25$ in another paper \cite{brendle2009blow}. Conversely, Khuri, Marques, and Schoen established compactness for $n\leq 24$ in \cite{khuri2009compactness}. Furthermore, they demonstrated compactness for $6\leq n\leq 24$ with $\sum_{i=0}^{\left[\frac{n-6}{2}\right]}|\nabla_g^i W|>0$ on $M$.
It's important to note that a crucial ingredient in the proofs of the aforementioned compactness results in dimensions $n\leq 24$ (without any assumption on the Weyl tensor $W$) is the deep positive mass theorem. For the proof of the positive mass theorem, refer to \cites{schoen1979complete,schoen1979proof,schoen1981energy} and \cite{Schoen1989variational} for $n\leq 7$; see \cite{witten1981new} for spin manifolds; and consult \cite{Schoen2022PMT}, as well as \cites{lohkamp2015hyperbolic, lohkamp2015skin, lohkamp2015skin2,lohkamp2016higher}, for $8\leq n\leq 24$.
On the other hand,
 the mentioned
 compactness results under $|W|>0$ or $|W|+|\nabla_g W|>0$ or $|W|+|\nabla_g W|+|\nabla^2_g W|>0$ or $\sum_{i=0}^{\left[\frac{n-6}{2}\right]}|\nabla_g^i W|>0$
do not rely on the positive mass theorem.

It's worth mentioning that another crucial ingredient in the proofs of all the aforementioned compactness results is the classic Liouville-type theorem of Caffarelli, Gidas, and Spruck in \cite{caffarelli1989asymptotic}.

In the Weyl non-vanishing direction, Marques provided examples in \cite{marques2009blow} demonstrating that 
$\mathcal{M}_g$ is not necessarily bounded in $L^\infty(M)$ 
under any one of the following three conditions: 
(1) $25 \leq n \leq 29$ and $\sum_{i=0}^l|\nabla_{g}^i W|>0$ for some $l \geq 7$; (2) $30 \leq  n \leq 51$ and $\sum_{i=0}^l|\nabla_{g}^i W|>0$ for some $l \geq 5$; (3) $n \geq  52$ and $\sum_{i=0}^l|\nabla_{g}^i W|>0$ for some $l \geq 3$.

We note that 
Berti and Malchiodi constructed in \cite{berti2001non} a
metric $\tilde{g}\in C^l(\s^n)$, $n\geq 4l+3$, so that $\mathcal{M}_{\tilde{g}}$ is not bounded in $L^\infty(\s^n)$.

The main result in this paper is as follows.

\begin{thm}\label{blowup2}
For $n\geq 25$, there  exist a smooth Riemannian metric $\tilde{g}$ on $\s^n$ and a sequence of positive 
functions $\{\tilde{v}_k\}$ solving
\begin{equation}\label{yamabe eqn on sphere}
-L_{\tilde{g}}\tilde{v}_k=n(n-2)\tilde{v}_k^{\frac{n+2}{n-2}}, \text{ on } \s^n,
\end{equation}
such that $
\Vol(\s^n,\tilde{g}_k)
\equiv 
\int_{\s^n}\tilde{v}_k^{\frac{2n}{n-2}}dV_{\tilde{g}}
\rightarrow\infty\text{ as } k\rightarrow \infty.$ 
\end{thm}

\begin{rem} The above theorem reinforces the findings of
Brendle and Marques in \cite{brendle2009blow}, where the examples  demonstrated $\max_{\mathbb{S}^n} \tilde{v}_k \rightarrow \infty$ as $k \rightarrow \infty$.
\end{rem}

We construct solutions 
$\{\tilde v_k\}$ near summations of an increasing number of Aubin-Talenti bubbles. The fact that the number of bubbles goes to infinity as $k$ goes to infinity leads to our unbounded-volume result. In comparison, the solutions in \cite{brendle2008blow} and  \cite{brendle2009blow} are near exactly one bubble with increasing heights.

We introduce weighted H{\"o}lder spaces which play a significant role in our proof of Theorem \ref{blowup2}. In \cite{brendle2008blow} and
\cite{brendle2009blow}, a pair of norms, $L^{\frac{2n}{n-2}}$-norm and $L^{\frac{2n}{n+2}}$-norm, were used.  Such  integral norms were used in carrying out
Lyapunov–Schmidt reductions in the study of other problems, see e.g.  \cite{bahri1991scalar},  \cite{bahri1995variational}, \cite{li1995prescribing} and
\cite{li1998dirichlet}.
However, the integral norms are more difficult
 to use to establish certain
 $k$-independent estimates. The weighted H{\"o}lder norms we use are sup norms in nature rather than integral norms. The idea of using weighted sup norms to overcome the difficulties caused by bubble interactions
 has appeared in different setups
 such as \cite{zelati1991homoclinic}, \cite{zelati1992homoclinic}, \cite{del2003two}, \cite{PISTOIA2007325}, and,  closer to our setup, \cite{wei2010infinitely} and \cite{li2018multi}. Bubble-tower radial solutions were constructed in 
  \cite{lin1999} and  \cite{DELPINO2003280} for 
 critical and slightly
super-critical elliptic equations. 

A non-negative critical point of the Yamabe functional $E_g$ gives a smooth 
positive solution of (\ref{yamabe eqn}).
Aubin \cite{aubin1976equations} proved the existence of a minimizer for the case $(M,g)$ non-locally conformally flat with $n\ge 6$. For the remaining two cases $n\leq 5$ and $(M,g)$ locally conformally flat with $n\ge 6$, Schoen \cite{schoen1984conformal} proved the existence of a minimizer.
The proof made use of the positive mass theorem of Schoen and Yau (see \cite{schoen1979proof} and  \cite{witten1981new}). It is worth noting that in these cases, Bahri and Brezis \cite{bahri1988equations} and Bahri \cite{Bahri1993}  later proved the existence of a critical point of $E_g$, not necessarily a minimizer, using the method of ``critical points at infinity''
without using the positive mass theorem.  

In CR geometry, there is an analogous CR Yamabe problem. Let the dimension of a CR manifold $M$ be $2n+1$. For the case $M$ non-locally CR equivalent to $\s^{2n+1}$ with $n\geq  2$, Jerison and Lee proved  in \cite{jerison1987yamabe} and  \cite{jerison1989intrinsic}  the existence of a minimizer of the corresponding CR Yamabe functional.

For the remaining two cases $n=1$ and $M$ locally CR equivalent to $\s^{2n+1}$ with $n\geq  2$, Gamara \cite{gamara2001cr} and Gamara and Yacoub \cite{gamara2001crconformallyflat} proved the existence of a critical point of the corresponding CR Yamabe functional using the method of ``critical points at infinity''. Cheng, Malchiodi and Yang showed in  \cite{cheng2023sobolev} that the corresponding CR Yamabe functional may not have a minimizer when $n=1$. They also showed that the corresponding CR positive mass theorem may fail when $n=1$.
The Liouville-type theorem  in the Heisenberg group,
corresponding to the above mentioned  Liouville-type theorem in 
\cite{caffarelli1989asymptotic},  was recently proved 
by Catino, Li, Monticelli and Roncoroni in \cite{CLMR} when $n=1$.

\noindent\textbf{Two Open Problems.} 

{\bf Problem 1}: By Theorem \ref{blowup2}, there exists $(M,\tilde{g})$ when $n\geq 25$ such that the solution space $\mathcal{M}_{\tilde{g}}$ is unbounded in $L^{\frac{2n}{n-2}}(M)$. On the other hand, it is easy to see that $\mathcal{M}_g$ is bounded in $L^{\frac{n+2}{n-2}}(M)$ for any closed manifolds $(M,g)$. So, there exists $p^*(n)\in \left[\frac{n+2}{n-2}, \frac{2n}{n-2}\right]$ such that for any $p>p^*(n)$, there exists a manifold $(M,\tilde{g})$ such that $\mathcal{M}_{\tilde{g}}$ is unbounded in $L^p(M)$; and for any $1\leq p<p^*(n)$, the solution space $\mathcal{M}_g$ is bounded in $L^p(M)$ for any closed manifolds $(M,g)$. The open problem is to identify the value of $p^*(n)$.

{\bf Problem 2}: As mentioned earlier, it has been proved,
without using the positive mass theorem, that $\mathcal{M}_{g}$ is compact in dimensions $n\geq 6$ under the assumption that 
$|W|+|\nabla_g W|>0$ on $M$, and $\mathcal{M}_{g}$ is compact in dimensions $6\leq  n\leq 51$ 
under the assumption that 
$|W|+|\nabla_g W|+|\nabla_g^2 W|>0$ on $M$. 
On the other hand, $\mathcal{M}_{g}$ is not necessarily compact under
the assumption that  $ |\nabla_{g}^3 W|>0$ on $M$ when $n\geq 52$.  
 There are still remaining open issues in this direction. For instance, whether or not
 $\mathcal{M}_{g}$ is compact 
  under the assumption $|W|+|\nabla_g W|+|\nabla_g^2 W|>0$ on $M$ 
in dimensions $n\ge 52$?

\noindent\textbf{Organization of the paper.} 

In Section 2, we provide the preliminaries,  including the approximate solutions, weighted function spaces, the nonlinear operator $F$, and its linearized operator $F'$. In Section 3, we present the basic properties of the weighted function spaces. 
In Section 4, we discuss the properties of the linearized operator $F'$, including its boundedness, Schauder estimates, and Fredholmness. Most importantly, we establish the invertibility of a projection of $F'$ restricted to a subspace. In Section 5.1, we show that this problem can be reduced to finding a critical point of a finite-dimensional function. In Section 5.2, we determine the leading term of the finite dimensional function. In Section 5.3, we conclude the paper by finding a local minimum point of this leading term, which in turn provides a local minimum point for the finite-dimensional function.

\noindent\textbf{Notations.} 

Here we list the notations which are used in the introduction and the rest of the paper.

1. Double indices follow the summation convention. We use Greek letters $\alpha, \beta, \mu, \nu$ taking on values from $\{1,\cdots,n\}$ or $\{0,\cdots,n\}$.

2. We use Latin letters $k,l\in \mathbb{N}$ to represent sequential indices, where $\mathbb{N}$ denotes the set of non-negative integers. We use Latin letters $i, j, p, q$ taking on values from $\{1,\cdots,k\}$ or $\{1,\cdots,l\}$.

3. For $\bar{x}\in \r^n$ and $\bar{r}\in \r^+$, we use $B(\bar{x},\bar{r}):=\{x\in\r^n\mid |x-\bar{x}|<\bar{r}\}$ to denote open balls. 

4. For a differentiable function $u=u(x)$ defined on an open subset $U\subseteq \r^n$, we use $D_{\mu}u$ to denote Euclidean partial derivatives $\partial_{x_{\mu}}u(x)$.

5. Use $\mathcal{D}^{1,2} (\r^n)$ to denote the completion of $C^{\infty}_c(\r^n)$ under the norm $\norm{\mathcal{D}^{1,2} }{v}:=\left(\int_{\r^n}|Dv|^2\right)^{\frac{1}{2}}$. It is a Hilbert space with the inner product $\<v_1,v_2\>_{\mathcal{D}^{1,2} }:=\int_{\r^n}Dv_1\cdot Dv_2$.

6. For $\alpha\in (0,1)$, an open subset $U\subseteq \r^n$, and a real valued function $f$ defined on $U$, we use $[f]_{\alpha,U}:=\sup_{y_1,y_2\in U}\frac{|f(y_1)-f(y_2)|}{|y_1-y_2|^{\alpha}}$ to denote the H{\"o}lder semi-norm. We use $|f|_U:=\sup_{y\in U}|f(y)|$ to denote the sup-norm. 

7. For $l\in \mathbb{N}$, $\alpha\in (0,1)$, and an open ball $B=B(0,d)$ with radius $d>0$, we use the following notation to denote the local weighted H{\"o}lder norms of a real valued function $f$ defined on $B$,
$$
|f|'_{l,\alpha,B}:=\sum_{i=0}^ld^i|D^if|_{B}+d^{l+\alpha}[D^lf]_{\alpha,B}.
$$

8. For a Banach space $(X,\norm{X}{\cdot})$, we use $(X^*,\norm{X^*}{\cdot})$ to denote the dual space. As usual, for any $\phi\in X^*$, $\norm{X^*}{\phi}:=\sup_{u\in X\setminus\{0\}}\frac{|\phi(u)|}{\norm{X}{u}}$.

9. For $l\in \mathbb{N}$, we use $C^l(\r^n)$ to denote the set of continuously differentiable up to order $l$ functions defined on $\r^n$. For $l\in \mathbb{N}$, $\alpha\in (0,1)$, we use
$$
C^{l,\alpha}_{\loc}(\r^n):= \{f\in C^l(\r^n)\mid \forall K\subseteq \r^n \text{ compact, } [D^lf]_{\alpha,K}<\infty \}.
$$

\noindent\textbf{Acknowledgements.}
L. Gong and Y.Y. Li are partially supported by NSF Grant DMS-2247410. L. Gong is also
partially supported by Professor Abbas Bahri Excellence Fellowship at Rutgers University. The authors thank the referee for carefully reading the paper and for helpful suggestions which significantly enhanced its exposition.

\section{Preliminaries and Outline}

\subsection{Background Metric}

First of all, to construct a background metric $\tilde{g}$ on $\s^n$, we let $\tilde{g}$ be the pullback metric by the stereographic projection $\s^n\rightarrow \r^n$ of a metric $\frac{4}{(1+|x|^2)^2}\hat{g}(x)$ on $x\in\r^n$. Consider
$$
\hat{g}=\exp (\epsilon \hat{h}),
$$
where $\epsilon>0$, $(\hat{h}_{\mu\nu})$ is an $n\times n$ symmetric matrix function on $\r^n$ with $\supp \hat{h}\subseteq B(0,1)$.

In this subsection, we will define matrix functions $\hat{H}$, $\hat{h}_k$, and $\hat{h}$ one by one.

Consider the following matrix function on $\r^n$, used in \cite{marques2009blow}. Similar examples can be found in \cite{brendle2009blow} and \cite{wei2013Qcurvature}. 
\begin{equation}\label{def of H}
\hat{H}_{\mu\nu}(\hat{x}):=(\tau_0+5|\hat{x}|^2-|\hat{x}|^4+\frac{1}{20}|\hat{x}|^6)W_{\mu\alpha\nu\beta}\hat{x}_{\alpha}\hat{x}_{\beta}, \quad \hat{x}\in \r^n,
\end{equation}
where $\tau_0\in \r$ is to be determined and $W: \r^n\times\r^n\times \r^n\times\r^n\rightarrow\r$ is a fixed multilinear form with all algebraic properties of Weyl tensors:
\begin{align*}
&W_{\mu\alpha\nu\beta}=W_{\nu\beta\mu\alpha}=-W_{\alpha\mu\nu\beta}=-W_{\mu\alpha\beta\nu};\\
&W_{\mu\alpha\nu\beta}+W_{\mu\nu\beta\alpha}+W_{\mu\beta\alpha\nu}=0;\\
&W_{\mu\alpha\mu\beta}=0.
\end{align*} and the following nontrivial condition 
$$
\sum_{\mu,\alpha,\nu,\beta=1}^{n}(W_{\mu\alpha\nu\beta}+W_{\mu\beta\nu\alpha})^2>0.
$$
Notice that $\hat{H}_{\mu\mu}(\hat{x})=0$, $\hat{x}_{\mu}\hat{H}_{\mu\nu}(\hat{x})=0$ and $D_{\mu}\hat{H}_{\mu\nu}(\hat{x})=0$. 

Consider a fixed non-negative cutoff function $\eta\in C_c^{\infty}((-1,1))$ satisfying $\eta(s)=1$ when $|s|<\frac{1}{2}$. For each integer $k\geq 3$, define
\begin{equation}\label{def of hath_k}
\hat{h}_k(x):=\sum_{j=1}^k\eta(k^4|x-\hat{P}^{k,j}|^2) \hat{H}(e^k (x-\hat{P}^{k,j})),
\end{equation}
where for each fixed $k$ and $j$ with $j\leq k$, the points
$$
\hat{P}^{k,j}:=O^{k,j}\hat{P}^{k} \text{ and } \hat{P}^k:=(k^{-1},0,0)\in \r\times \r \times \r^{n-2},
$$
and the rotation matrices
\begin{equation}\label{rotation matrices}
O^{k,j}:=\begin{bmatrix}
\cos\frac{2\pi(j-1)}{k} & \sin\frac{2\pi(j-1)}{k} & 0\\
-\sin\frac{2\pi(j-1)}{k} & \cos\frac{2\pi(j-1)}{k} & 0\\
0 & 0 & I_{n-2}
\end{bmatrix}.
\end{equation}
For $k\geq 3$, the support of $\hat{h}_k$ is contained in the union of $k$ disjoint open balls: $\{B(\hat{P}^{k,j},k^{-2})\}_{1\leq j\leq k}$. The open balls have the following property 
$$
\dist\left(B(\hat{P}^{k,j_1},k^{-2}),B(\hat{P}^{k,j_2},k^{-2})\right)\geq \frac{1}{3}k^{-1}, \text{ for any }j_1\neq j_2\in \{1,\cdots,k\}.
$$
\textbf{Remark.} For simplicity, we have arranged the $k$ points $\{\hat{P}^{k,j}\}_{1\leq j\leq k}$ with symmetry. Similar setup was used by Haskins and Kapouleas \cite{haskins2007special} and Wei and Yan \cite{wei2010infinitely}. In fact, the symmetry here does not play an essential role. The idea of this paper can still work if we arrange the points nicely but without symmetry.

Now, we are ready to give the definition of $\hat{h}$. For $c_0>0$,
\begin{equation}\label{def of hath}
\hat{h}(x):=\sum_{k=3}^{\infty}e^{-(8+c_0)k}\hat{h}_k(x),
\end{equation}
where $\hat{h}_k$ is defined by (\ref{def of hath_k}). Also, notice that $\supp \hat{h}\subseteq B(0,1)$,
$$
\hat{h}_{\mu\mu}(x)=0 \text{ and } D_{\mu}\hat{h}_{\mu\nu}(x)=0 \text{ for } x\in \r^n.
$$
Here, $\hat{h}$ is clearly smooth away from the origin. It is also smooth at the origin because for any integer $m\geq 0$ there exists a positive constant $C$ only depends on $n,m$, the multilinear form $W$ and the cutoff function $\eta$ such that
\begin{equation*}
|D^m\hat{h}|(x)\leq C k^{2m}e^{-c_0k}, \text{ for any } 0<|x|<2k^{-1}, k\geq 3.
\end{equation*}
Consequently, the Weyl tensor of $\hat{g}$ together with its derivatives of any order vanish at the origin. In this paper, we fix a multilinear form $W$ and a cutoff function $\eta$ such that $\norm{C^{3}(\r^n)}{\hat{h}}<1$.

\subsection{Approximate Solutions}

There are three numbers to be determined in the definition of $\hat{g}$. 

\noindent a. $\tau_0\in \r$ appears in (\ref{def of H}), the definition of $\hat{H}$;\\
b. $c_0>0$ appears in (\ref{def of hath}), the definition of $\hat{h}$;\\
c. $\epsilon>0$ appears in the definition of $\hat{g}$: $\hat{g}=\exp(\epsilon \hat{h})$.

By the stereographic projection, solving equation (\ref{yamabe eqn on sphere}) can be reduced to solving:
\begin{equation}\label{beforescale}
-L_{\hat{g}}\hat{v}_k=n(n-2)\hat{v}_k^{\frac{n+2}{n-2}}, \text{ in }\r^n.
\end{equation}

For each fixed $k\geq 3$, we pick $t_k=e^{-k}$ as the scaling parameter and do scaling $v_k(x):=t_k^{\frac{n-2}{2}}\hat{v}_k(t_kx)$. So, $\hat{v}_k$ solving equation (\ref{beforescale}) is equivalent to $v_k$ solving the following equation
\begin{equation}\label{afterscale}
-L_{g_{\epsilon, k}}v_k=n(n-2)v_k^{\frac{n+2}{n-2}}, \text{ in }\r^n,
\end{equation}
where $g_{\epsilon, k}(x):=\hat{g}(t_kx)=\exp (\epsilon \hat{h}(t_kx))$. 

In the following, we will give approximate solutions to equation (\ref{afterscale}). Firstly, denote the Aubin-Talenti bubbles by
$$
\sigma_{\xi,\lambda}(x):=\lambda^{\frac{n-2}{2}}(1 +\lambda^2|x-\xi|^2)^{-\frac{n-2}{2}}, \quad x\in \r^n,
$$
where $\xi\in \r^n$ and $\lambda\in \r^+$. According to Caffarelli, Gidas, and Spruck \cite{caffarelli1989asymptotic}, Aubin-Talenti bubbles are the only positive solutions of $-\Delta v=n(n-2)v^{\frac{n+2}{n-2}}$ in $\r^n$, where $n\geq 3$. 
For each integer $k\geq 3$, consider a parameter set
\begin{align*}
\mathcal{C}_k:=&\left\{\Xi=(\xi^1,\cdots,\xi^k)\text{ and }\Lambda=(\lambda_1,\cdots,\lambda_k)\mid \right.\\
&\left.(\xi^i,\lambda_i)\in\r^n\times \r^+, \frac{3}{4}<\lambda_i<\frac{4}{3},|\xi^i|<\frac{1}{2}, i=1,\cdots,k\right\}.
\end{align*}

For any $(\Xi,\Lambda)\in \mathcal{C}_k$, consider the approximate solutions to equation (\ref{afterscale}):
$$
u_{k,\Xi,\Lambda}(x)=\sum_{j=1}^k\sigma_{P^{k,j}+\xi^j,\lambda_j}(x),
$$
where $P^{k,j}:=t_k^{-1}\hat{P}^{k,j}$ and $\hat{P}^{k,j}:=O^{k,j}(k^{-1},0,0)$.

\begin{thm}\label{blowup3}
Let $n\geq 25$, $0<c_0<\frac{n-2}{2}-8$. There exist $\tau_0=\tau_0(n)\in \r$, $\epsilon=\epsilon(n, c_0)>0$, and $k_0=k_0(n, c_0)>0$ such that for any $k\geq k_0$, equation (\ref{afterscale}) has a positive smooth solution $v_{k}$ of the form
$$
v_k=u_{k,\bar{\Xi},\bar{\Lambda}}+\phi_k,
$$ 
where $(\bar{\Xi},\bar{\Lambda})\in \mathcal{C}_k$ and $\norm{L^{\frac{2n}{n-2}}(\r^n)}{\phi_k}\rightarrow 0$ as $k\rightarrow \infty$. Consequently, $\{v_k\}$ satisfies
\begin{equation*}
\int_{\r^n}v_k^{\frac{2n}{n-2}}dx\rightarrow \infty \text{ as } k\rightarrow \infty.
\end{equation*}
\end{thm}

After scaling back, $\hat{v}_k(x)=t_k^{-\frac{n-2}{2}}v_{k}(t_k^{-1}x)$ solves equation (\ref{beforescale}). Then composing with the stereographic projection, $\tilde{v}_k$ solves equation (\ref{yamabe eqn on sphere}) on $\s^n$. The volume blowup result in Theorem \ref{blowup2} comes from
$$
\Vol(\s^n,\tilde{g}_k)=\int_{\s^n}\tilde{v}_k^{\frac{2n}{n-2}}dV_{\tilde{g}}=\int_{\r^n}\hat{v}_k^{\frac{2n}{n-2}}dV_{\hat{g}}=\int_{\r^n}\hat{v}_k^{\frac{2n}{n-2}}dx\rightarrow\infty,
$$
where we have used $\det \hat{g}=\exp (\epsilon\trace \hat{h})=1$ in the last equality.

\subsection{Function Spaces and Operators}

Consider functions
\begin{equation}\label{def of tangent vectors}
\begin{aligned}
Z_{k,j,0}(x)&:=\partial_{\lambda_j}u_{k,\Xi,\Lambda}(x)=\partial_{\lambda_j}\sigma_{P^{k,j}+\xi^j,\lambda_j}(x),\text{ for any }j=1,\cdots,k;\\
Z_{k,j,\mu}(x)&:=\partial_{\xi^{j}_{\mu}}u_{k,\Xi,\Lambda}(x)=\partial_{\xi^j_{\mu}}\sigma_{P^{k,j}+\xi^j,\lambda_j}(x),\text{ for any }j=1,\cdots,k, \mu=1,\cdots,n,
\end{aligned}
\end{equation}
and two subspaces of $\mathcal{D}^{1,2}$,
\begin{equation*}
\begin{aligned}
T_{k,\Xi,\Lambda}&:=\spa\{Z_{k,j,\mu}\mid j=1,\cdots,k;\mu=0,\cdots,n\};\\
N_{k,\Xi,\Lambda}&:=\left\{\phi\in \mathcal{D}^{1,2} \mid \<Z_{k,j,\mu},\phi\>_{\mathcal{D}^{1,2} }=0, \forall j,\mu\right\}.
\end{aligned}
\end{equation*}

For each $k\geq 3$, consider $l\in \mathbb{N}, \alpha\in (0,1),s\in \r$, and we introduce the following weighted sup-norm and weighted H{\"o}lder semi-norm:
\begin{align}
\norm{X^{l}_{k,s}}{\phi}&=\sum_{i=0}^l\sup_{x\in\r^n}\left(\sum_{j=1}^k\<x-P^{k,j}\>^{-s-i}\right)^{-1}|D^i\phi|(x),\\
[\phi]_{X^{l,\alpha}_{k,s}}&=\sup_{x\in \r^n} \left(\sum_{j=1}^k\<x-P^{k,j}\>^{-s-l-\alpha}\right)^{-1}[D^l\phi]_{\alpha,B(x,d(x)/2)},
\end{align}
where
$$
\<x\>:=(1+|x|^2)^{\frac{1}{2}},\text{ and }d(x):=\min_{j=1,\cdots,k}\<x-P^{k,j}\>.
$$
Then, we define the weighted spaces
\begin{align*}
X^{l}_{k,s}&:=\left\{\phi\in C^{l}_{\loc}(\r^n)\mid \norm{X^{l}_{k,s}}{\phi}<\infty\right\},\\
X^{l,\alpha}_{k,s}&:=\left\{\phi\in C^{l,\alpha}_{\loc}(\r^n)\mid \norm{X^{l,\alpha}_{k,s}}{\phi}<\infty\right\}.
\end{align*}
where 
$$
\norm{X^{l,\alpha}_{k,s}}{\phi}:=\norm{X^{l}_{k,s}}{\phi}+[\phi]_{X^{l,\alpha}_{k,s}}.
$$

Now, we give the definitions of two operators. For any $\epsilon>0$, $k\geq 3$, consider a nonlinear map:
\begin{equation}\label{def of F}
F_{\epsilon,k}(v):=-L_{g_{\epsilon,k} }v-n(n-2)|v|^{\frac{4}{n-2}}v,
\end{equation}
where $g_{\epsilon,k}(x)=\exp (\epsilon \hat{h}(t_kx))$. For any $(\Xi,\Lambda)\in \mathcal{C}_k$, consider the linearized elliptic operator of $F_{\epsilon,k}$ at $u_{k,\Xi,\Lambda}$:
\begin{equation}\label{def of L}
L_{\epsilon,k,\Xi,\Lambda}\phi:=-L_{g_{\epsilon,k}}\phi-n(n+2)u_{k,\Xi,\Lambda}^{\frac{4}{n-2}}\phi.
\end{equation}

Solving equation (\ref{afterscale}) is equivalent to solving a positive zero point $F_{\epsilon,k}(v)=0$ of the nonlinear map $F_{\epsilon,k}$. It has a variational formulation: $F_{\epsilon,k}(v)=0$ is the Euler Lagrange equation of the following functional
\begin{equation}\label{def of energy I}
I_{g_{\epsilon,k}}(v):=\int_{\r^n}\frac{1}{2}|D_{g_{\epsilon,k}}v|^2+\frac{c(n)}{2}R_{g_{\epsilon,k} }v^2-\frac{(n-2)^2}{2}|v|^{\frac{2n}{n-2}}, \text{ for } v\in \mathcal{D}^{1,2}.
\end{equation}

\textbf{Simplified Notations.} In this paper, we work with Banach spaces $X^{2, \alpha}_{k,s}$ and $X^{0, \alpha}_{k,s+2}$ and for the simplicity of notations, we write 
\begin{equation}\label{def of X and Y}
X:=X^{2, \alpha}_{k,s} \text{ and } Y:=X^{0, \alpha}_{k,s+2}.
\end{equation}
We will drop the subscript $\epsilon, k$ and write 
$$
g:=g_{\epsilon,k}, \quad F:=F_{\epsilon,k}, \quad L_{\Xi,\Lambda}:=L_{\epsilon,k,\Xi,\Lambda}.
$$
Moreover, we will drop the subscript or superscript $k$ as follows:
$$
P^j:=P^{k,j}, \quad u_{\Xi,\Lambda}:=u_{k,\Xi,\Lambda},\quad Z_{j,\mu}:=Z_{k,j,\mu}, \quad T_{\Xi,\Lambda}:=T_{k,\Xi,\Lambda}, \quad  N_{\Xi,\Lambda}:=N_{k,\Xi,\Lambda}.
$$

\subsection{Outline}

In Section 3, we give some basic properties of the weighted H{\"o}lder norms. Among them, we prove that $ X \subseteq \mathcal{D}^{1,2} $ and the inclusion map is continuous. Thus, the following two subspaces are closed in $X$:
\begin{equation}\label{def of X_1,2}
X_1:=T_{\Xi,\Lambda} \text{ and } X_2:=N_{\Xi,\Lambda}\cap  X.
\end{equation}
Since $\mathcal{D}^{1,2} =T_{\Xi,\Lambda}\oplus N_{\Xi,\Lambda}$, we have a decomposition $X =X_1\oplus X_2$. 

In Section 4, we focus on the properties of the linearized operator $L_{\Xi,\Lambda}$. First, we will prove that for any $0<\epsilon<1$, $k\geq 3$, $\frac{n-2}{2}<s<n-2$ and $0<\alpha<\frac{4}{n-2}$, $F: X \rightarrow  Y $ is a $C^1$ map with $F'(u_{\Xi,\Lambda})=L_{\Xi,\Lambda}:  X \rightarrow  Y $. Moreover, define two subspaces of $Y$,
\begin{equation}\label{def of Y_1,2}
Y_1:=\Delta(X_1)\text { and }Y_2:= L_{\Xi,\Lambda}(X_2).
\end{equation}
We will prove $Y_1, Y_2$ are closed and $Y$ admits a decomposition: $Y =Y_1\oplus Y_2$. 

Second, we will prove a $k$-independent Schauder estimate of the linearized operator $L_{\Xi,\Lambda}$ and a $k$-independent invertibility of $\bar{L}:=\Proj_{Y_2}\circ L_{\Xi,\Lambda}\circ \iota: X_2\rightarrow Y_2$, where $\iota:X_2\rightarrow  X $ is the inclusion map. In fact, we will prove there exists a $k$-independent positive constant $C$ such that
\begin{align*}
\norm{ X }{\bar{L}^{-1}f_2}&\leq C\norm{ Y }{f_2},\text{ for any }f_2\in Y_2,\\
\norm{ Y }{\Proj_{Y_2}f}&\leq C\norm{ Y }{f}, \text{ for any }f\in  Y,\\
\norm{ X }{\phi}&\leq C(\norm{ Y }{L_{\Xi,\Lambda}\phi}+\norm{X^0_{k,s}}{\phi}), \text{ for any }\phi \in X.
\end{align*}

In Section 5.1, we divide solving the equation $F(v)=0$ into solving the following two equations:
\begin{align*}
\Proj_{Y_2}\circ F(v)&=0;\\
(\Id-\Proj_{Y_2})\circ F(v)&=0.
\end{align*}
For the first equation, we will use a version of implicit function theorem to prove that for any $(\Xi,\Lambda)\in \mathcal{C}_k$, there exists $\phi_{\Xi,\Lambda}\in X_2$ such that
$$
\Proj_{Y_2}\circ F(u_{\Xi,\Lambda}+\phi_{\Xi,\Lambda})=0.
$$
For the second equation, by the variational formulation of $F$, we reduce finding $(\bar{\Xi},\bar{\Lambda})\in\mathcal{C}_k$ such that
$$
(\Id-\Proj_{Y_2})\circ F(u_{\bar{\Xi},\bar{\Lambda}}+\phi_{\bar{\Xi},\bar{\Lambda}})=0,
$$
to finding a critical point of the following finite dimensional function $I:\mathcal{C}_k\rightarrow \r$:
\begin{equation*}
I(\Xi,\Lambda):=I_{g}(u_{\Xi,\Lambda}+\phi_{\Xi,\Lambda}).
\end{equation*}
Then, letting $v_k=u_{\bar{\Xi},\bar{\Lambda}}+\phi_{\bar{\Xi},\bar{\Lambda}}$, we will solve equation (\ref{afterscale}) and prove Theorem \ref{blowup3}.

In the rest of Section 5, we will find a critical point of $I(\Xi,\Lambda)$ by finding the leading term of $I(\Xi,\Lambda)-kI_0$, where $I_0:=I_{\delta}(\sigma_{\xi,\lambda})>0$ is a constant independent of $(\xi,\lambda)\in \r^n\times \r^+$, and then finding a local minimum point of the leading term. First, in Section 5.2, we will prove the following approximation: as $k\rightarrow \infty$, for any $(\Xi,\Lambda)\in \mathcal{C}_k$, ,
\begin{equation*}
I(\Xi,\Lambda)-kI_0=\epsilon^2 t_k^{16+2c_0}(G_k(\Xi,\Lambda)+o(k^{-2})).
\end{equation*}
Second, in Section 5.3, we will find a special open neighborhood $B((\Xi_0,\Lambda_0),k^{-1})\subseteq \mathcal{C}_k$ and prove that
$$
\min_{\partial B((\Xi_0,\Lambda_0),k^{-1})} G_k(\Xi,\Lambda)\geq G_k(\Xi_0,\Lambda_0)+C^{-1}k^{-2},
$$
where the $k$-independence of the positive constant $C$ is important.

\section{Properties of the Weighted Norms}

Recall that given $l\in \mathbb{N}$, $\alpha\in (0,1)$, and $s\in \r$, the modified weighted H{\"o}lder norms $\norm{X^{l,\alpha}_{k,s}}{\cdot}$ only depend on the choice of $k$ points in $\r^n$. So, in this section, we will work in a slightly more general setup:\\
For fixed $k\geq 3$ and $r>1$, pick $(r,0,0)\in \r\times \r \times \r^{n-2}$ and consider $k$ points 
$$
P^j=O^{k,j}(r,0,0), \text{ for } j=1,\cdots,k,
$$
where $O^{k,j}$ are the rotation matrices given by (\ref{rotation matrices}). 

In this setup, there are several simple facts. 

1. There exists $C=C(n)>0$ such that 
\begin{equation}\label{geometry}
|P^j-P^1|\geq \frac{1}{C}\left\{
\begin{aligned}
(j-1)\frac{r}{k} &\qquad\text{ for } j=1,\cdots, \left[\frac{k}{2}\right];\\
(k-j+1)\frac{r}{k} &\qquad\text{ for } j=\left[\frac{k}{2}\right]+1,\cdots, k.
\end{aligned}\right.  
\end{equation}

2. According to \cite[Lemma B.1]{wei2010infinitely} and \cite[Lemma A.1]{li2018multi}, for $a,b>0$ and $0<\tau\leq \min\{a,b\}$, there exists $C=C(n,\tau)>0$ such that for $i, j=1,\cdots, k$ and $i\neq j$, for any $x\in\r^n$,
\begin{equation}\label{fact 1}
\<x-P^i\>^{-a}\<x-P^j\>^{-b}\leq C|P^i-P^j|^{-\tau}\left(\<x-P^i\>^{-a-b+\tau}+\<x-P^j\>^{-a-b+\tau}\right).
\end{equation}
In fact, this inequality holds for any two distinct points $P^i, P^j\in \r^n$.

3. There exists $C=C(n)>0$ such that for any $x\in \r^n$ and $i=1,\cdots,k$,
\begin{equation}\label{B geometry}
C^{-1}\<x-P^i\>\leq \<y-P^i\>\leq C\<x-P^i\>, \text{ for any } y\in B(x,d(x)/2).
\end{equation}
Note that we can use ``triangle inequality'' and $y\in B(x,d(x)/2)$ to obtain that $\<x-P^i\>\leq C\<y-P^i\>+|x-y|\leq C\<y-P^i\>+d(x)/2$ and then obtain the first inequality in (\ref{B geometry}).

4. In some of the proofs of this paper, it is enough to just consider the following domain $\Omega_1$ rather than $\r^n$.
\begin{equation}\label{Omega_1}
\Omega_1:=\left\{(x_1,x_2,x'')\in \r^n\mid \frac{x_1}{\sqrt{x_1^2+x_2^2}}>\cos\frac{\pi}{k} \right\}.
\end{equation}
For $x\in \Omega_1$, there are several geometric relations
\begin{equation}\label{Omega_1 geometry}
\begin{aligned}
d(x)&=\<x-P^1\>;\\
\<x-P^i\>&\geq\<x-P^1\>, \text{ for } i=1,\cdots,k;\\
\<x-P^i\>&\geq \frac{1}{2}|P^i-P^1|, \text{ for } i=1,\cdots,k.
\end{aligned}
\end{equation}

Also, recall the classical weighted H{\"o}lder spaces $C^{l,\alpha}_s(\r^n)$ for $l\in \mathbb{N}, \alpha\in (0,1),s\in \r$, see e.g. \cite[Appendix A.2]{lee2021geometric}. For any $\phi\in C^{l,\alpha}_{\loc}(\r^n)$, define the following classical weighted sup-norm and classical weighted H{\"o}lder norm:
\begin{align*}
\norm{C^{l}_s}{\phi}&:=\sum_{j=0}^l\sup_{x\in\r^n}\<x\>^{s+j}|D^j\phi|(x),\\
[\phi]_{C^{l,\alpha}_s}&:=\sup_{x\in \r^n} \<x\>^{s+l+\alpha}[D^l\phi]_{\alpha,B(x,\<x\>/2)},\\
\norm{C^{l,\alpha}_s}{\phi}&:=\norm{C^{l}_s}{\phi}+[\phi]_{C^{l,\alpha}_s}.
\end{align*}
The classical weighted space is defined to be 
$$
C^{l,\alpha}_s(\r^n):=\left\{\phi\in C^{l,\alpha}_{\loc}(\r^n)\mid \norm{C^{l,\alpha}_s}{\phi}<\infty\right\}.
$$

\begin{prop}\label{modified norm}
Let $n\geq 5$, $\frac{n-2}{2}<s\leq n-2$,  $l\geq 1$ and $\alpha\in (0,1)$. There exists a positive constant $C=C(n,s,l,\alpha)$, independent of $k$ and $r$, such that for any $\phi \in X^{l,\alpha}_{k,s}$, $f\in  Y $, $(\Xi,\Lambda)\in \mathcal{C}_k$,
\begin{align}
C^{-1}k^{-1}r^{-l-\alpha-s}\norm{C^{l,\alpha}_s}{\phi}\leq\norm{X^{l,\alpha}_{k,s}}{\phi}&\leq Cr^{l+\alpha+s}\norm{C^{l,\alpha}_s}{\phi},\label{equivlent norm}\\
\norm{\mathcal{D}^{1,2} }{\phi}&\leq Ck(k/r)^{\frac{1}{2}(s-\frac{n-2}{2})}\norm{X^{l,\alpha}_{k,s}}{\phi},\label{D^12 embedding}\\
\norm{X^*}{f}&\leq Ck^{2}(k/r)^{s-\frac{n-2}{2}}\norm{ Y }{f},\label{dual space}\\
\norm{ X }{u_{\Xi,\Lambda}}&\leq C,\label{bound of u}\\
\norm{ X }{\sum_{j,\mu} Z_{j,\mu}}&\leq C.\label{bound of Z}
\end{align}
Note that $f\in  X^* $ is viewed as the linear functional $\phi \mapsto\int_{\r^n}f\phi$ where $\phi\in  X $.
\end{prop}

Using (\ref{equivlent norm}) and the Fredholmness of elliptic operators in classical weighted spaces, see e.g. \cite[Corollary A.42]{lee2021geometric}, we can obtain the following corollary.
\begin{coro}\label{fredholm}
Let $n\geq 5$, $1<s< n-2$, $0<\alpha<1$, $c_0>0$, and $\tau_0\in\r$. For any $\epsilon>0$, $k\geq 3$ and $(\Xi,\Lambda)\in \mathcal{C}_k$, the linearized operator
$$
L_{\Xi,\Lambda}:\left( X ,\norm{ X }{\cdot}\right)\rightarrow \left( Y ,\norm{ Y }{\cdot}\right) 
$$
is a Fredholm operator of index zero.
\end{coro}
\textbf{Remark.} In Proposition \ref{modified norm}, $s>\frac{n-2}{2}$ is used to prove estimates (\ref{D^12 embedding}) and (\ref{dual space}) and $s\leq n-2$ is used to prove estimates (\ref{bound of u}) and (\ref{bound of Z}). In Corollary \ref{fredholm}, $s>1$ is because of the usage of (\ref{equivlent norm}) and $s<n-2$ is to make sure the index of $L_{\Xi,\Lambda}$ is zero.

To prove Proposition \ref{modified norm}, we will give several lemmas.

\begin{lem}\label{lem3}
Let $s, s_1, s_2>0$ and $0< \tau\leq \min\{s, s_1, s_2\}$, there exists $C=C(n,\tau)>0$ such that for any $x\in \r^n$:
\begin{equation}\label{tau eqn1}
\sum_{i=1}^k\<x-P^i\>^{-s}-d(x)^{-s}\leq Ck\left(\frac{k}{r}\right)^{\tau}d(x)^{-s+\tau},
\end{equation}
and 
\begin{equation}\label{tau eqn2}
\sum_{i=1}^k\<x-P^i\>^{-s_1}\sum_{j=1}^k\<x-P^j\>^{-s_2}-\sum_{i=1}^k\<x-P^i\>^{-s_1-s_2}\leq Ck\left(\frac{k}{r}\right)^{\tau}\sum_{i=1}^k\<x-P^i\>^{-s_1-s_2+\tau}.
\end{equation}
\end{lem}
\textbf{Remark.} The $k\left(\frac{k}{r}\right)^{\tau}$ in this lemma can be replaced by $\left(\frac{k}{r}\right)^{\tau}$ when $\tau>1$. Similar idea can be seen in \cite[Lemma B.3]{wei2010infinitely}.

\begin{proof}
By (\ref{geometry}), for any $\tau>0$, we obtain
\begin{equation}\label{compute summation}
\sum_{j=2}^k|P^j-P^1|^{-\tau}\leq C\left(\frac{k}{r}\right)^{\tau}\sum_{j=2}^kj^{-\tau}\leq Ck\left(\frac{k}{r}\right)^{\tau}.
\end{equation}
Because of the symmetry of $d(x)$ and $\sum_{i}\<x-P^i\>^{-s}$, we just need to prove the inequalities (\ref{tau eqn1}) and (\ref{tau eqn2}) in $\Omega_1$ defined in (\ref{Omega_1}). 

To prove (\ref{tau eqn1}), we use the two inequalities in (\ref{Omega_1 geometry}) to obtain that
\begin{align*}
\sum_{i=1}^k\<x-P^i\>^{-s}-\<x-P^1\>^{-s}\leq&\<x-P^1\>^{-s+\tau}\sum_{i=2}^k\<x-P^i\>^{-\tau}\\
\leq&C\<x-P^1\>^{-s+\tau}\sum_{i=2}^k|P^i-P^1|^{-\tau}.
\end{align*}
Then by the equality in (\ref{Omega_1 geometry}) and (\ref{compute summation}), we can have inequality (\ref{tau eqn1}).

To prove (\ref{tau eqn2}), consider
$$
\sum_{i=1}^k\<x-P^i\>^{-s_1}\sum_{j=1}^k\<x-P^j\>^{-s_2}-\sum_{i=1}^k\<x-P^i\>^{-s_1-s_2}=\sum_{i=1}^k\sum_{j\neq i}\<x-P^i\>^{-s_1}\<x-P^j\>^{-s_2}.
$$
Then, we use inequality (\ref{fact 1}) for $a=s_1$ and $b=s_2$,
\begin{align*}
&\sum_{i=1}^k\sum_{j\neq i}\<x-P^i\>^{-s_1}\<x-P^j\>^{-s_2}\\
\leq &C\sum_{i=1}^k\sum_{j\neq i}(\<x-P^i\>^{-s_1-s_2+\tau}+\<x-P^j\>^{-s_1-s_2+\tau})|P^j-P^i|^{-\tau}\\
=&C\sum_{i=1}^k\<x-P^i\>^{-s_1-s_2+\tau}\sum_{j\neq i}|P^j-P^i|^{-\tau}.
\end{align*}
Then by (\ref{compute summation}), we can have inequality (\ref{tau eqn2}).
\end{proof}

The following lemma plays a key role to prove Proposition \ref{modified norm}. It is also used frequently in the rest of this paper. 
\begin{lem}\label{step function lemma}
For $n\geq 5$, $\tau>0$, $s\geq 1+\tau$, there exists a positive constant $C=C(n,\tau)$, independent of $k, r, s$, such that
\begin{equation}\label{step function gamma}
C^{-1}\gamma(d(x))d(x)^{-s}\leq \sum_{i=1}^k\<x-P^i\>^{-s}\leq C\gamma(d(x))d(x)^{-s}, \quad\text{ for } x\in\r^n,
\end{equation}
where $\gamma$ is defined on $\r^+$:
\begin{equation*}
\gamma(\rho):=\left\{\begin{aligned}
1, &\qquad 0<\rho\leq \frac{r}{k};\\
j+1, &\qquad j\frac{r}{k}<\rho\leq (j+1)\frac{r}{k}, j=1,\cdots,\left[\frac{k}{2}\right];\\
\left[\frac{k}{2}\right]+1, &\qquad \rho>\frac{r}{2},
\end{aligned}\right.
\end{equation*}
In particular, for $s_1,s_2\geq 1+\tau$,
\begin{equation}\label{step function equation}
C^{-1}d(x)^{s_1-s_2}\leq \left(\sum_{j=1}^k\<x-P^j\>^{-s_1}\right)^{-1}\sum_{i=1}^k\<x-P^i\>^{-s_2}\leq Cd(x)^{s_1-s_2}.
\end{equation}
\end{lem}

\textbf{Remark.} An immediate corollary of (\ref{step function equation}) is that for any $1+\tau\leq s_1\leq s_2$, there exists a positive constant $C=C(n,\tau)$, we have $\norm{X^{l,\alpha}_{k,s_1}}{\phi}\leq C\norm{X^{l,\alpha}_{k,s_2}}{\phi}$.

\begin{proof}
Because of the symmetry of $d(x)$ and $\sum_{i}\<x-P^i\>^{-s}$, we just need to prove the inequality (\ref{step function gamma}) in $\Omega_1$ defined as (\ref{Omega_1}).

We divide $\Omega_1$ into $\left[\frac{k}{2}\right]+2$ parts: 
for $j=0,\cdots,\left[\frac{k}{2}\right]$,
$$
j\frac{r}{k}\leq\<x-P^1\>\leq (j+1)\frac{r}{k},\text{ and }, \<x-P^1\>\geq\frac{r}{2}.
$$

Firstly, let us consider the last part: $\<x-P^1\>\geq\frac{r}{2}$. By this assumption, for $i=1,\cdots,k$,
$$
\<x-P^i\>\leq C\<x-P^1\>+|P^1-P^i|\leq C\<x-P^1\>+\frac{r}{2}\leq C\<x-P^1\>.
$$
Then, combining the above estimate with the first inequality in (\ref{Omega_1 geometry}), we have
$$
C^{-1}k\<x-P^1\>^{-s}\leq \sum_{i=1}^k\<x-P^i\>^{-s}\leq  k\<x-P^1\>^{-s}.
$$
So we have proved inequality (\ref{step function gamma}) for $\Omega_1\cap \{\<x-P^1\>\geq\frac{r}{2}\}$.

Secondly, consider the other parts $j\frac{r}{k}\leq\<x-P^1\>\leq (j+1)\frac{r}{k}$ where $j=0,\cdots,\left[\frac{k}{2}\right]$. By the first inequality in this assumption, for $i=1,\cdots, j+1$:
$$
\<x-P^i\>\leq C\<x-P^1\>+|P^1-P^i|\leq C\<x-P^1\>+C(i-1)\frac{r}{k} \leq C\<x-P^1\>.
$$
So we can obtain the lower bound of inequality (\ref{step function gamma}):
$$
\sum_{i=1}^k\<x-P^i\>^{-s}\geq 2\sum_{i=1}^{j+1}\<x-P^i\>^{-s}\geq C^{-1}(j+1)\<x-P^1\>^{-s}.
$$
To obtain the upper bound, we divide the summation into two parts and use the two inequalities from (\ref{Omega_1 geometry}) correspondingly:
\begin{align*}
\sum_{i=1}^k\<x-P^i\>^{-s}\leq &2\sum_{i=1}^{j+1}\<x-P^i\>^{-s}+2\sum_{i=j+2}^{\left[\frac{k}{2}\right]+1}\<x-P^i\>^{-s}\\
\leq &2(j+1)\<x-P^1\>^{-s}+C\sum_{i=j+2}^{\left[\frac{k}{2}\right]+1}|P^1-P^i|^{-s}.
\end{align*}
Then because of (\ref{geometry}) and the second inequality in the assumption, we have
$$
\sum_{i=j+2}^{\left[\frac{k}{2}\right]+1}|P^1-P^i|^{-s}\leq 
C\sum_{i=j+2}^{\left[\frac{k}{2}\right]+1}\left((i-1)\frac{r}{k}\right)^{-s}\leq C\sum_{i=j+1}^{\left[\frac{k}{2}\right]}i^{-s}\left(\frac{\<x-P^1\>}{j+1}\right)^{-s}.
$$
Combining the above three estimates, we have  
$$
C^{-1}(j+1)\<x-P^1\>^{-s}\leq \sum_{i=1}^k\<x-P^i\>^{-s}\leq  C(j+1)\<x-P^1\>^{-s}.
$$
So, we have proved inequality (\ref{step function gamma}) for each part $\Omega_1\cap \{j\frac{r}{k}\leq\<x-P^1\>\leq (j+1)\frac{r}{k}\}$ where $j=0,\cdots,\left[\frac{k}{2}\right]$. 

Combining all the $\left[\frac{k}{2}\right]+2$ parts together, we obtain (\ref{step function gamma}) for $x\in \Omega_1$. By symmetry, we obtain (\ref{step function gamma}) for $x\in\r^n$.
\end{proof}

Before the next lemma, consider another two weighted H{\"o}lder norms:
\begin{align*}
\norm{(X^{l,\alpha}_{k,s})'}{\phi}&:=\norm{(X^{l}_{k,s})'}{\phi}+\sup_{x\in \r^n} d(x)^{s+l+\alpha}[D^l\phi]_{\alpha,B(x,d(x)/2)},\\
\norm{(X^{l,\alpha}_{k,s})''}{\phi}&:=\norm{(X^{l}_{k,s})'}{\phi}+\sup_{x,y\in \r^n} \min\{d(x),d(y)\}^{s+l+\alpha}\frac{|D^l\phi(x)-D^l\phi(y)|}{|x-y|^{\alpha}},
\end{align*}
where $\norm{(X^{l}_{k,s})'}{\phi}:=\sum_{j=0}^l\sup_{x\in\r^n}d(x)^{s+j}|D^j\phi|(x)$.

\begin{lem}
Let $n\geq 5$, $\tau\geq 0$, $s\geq 1+\tau$,  $l\geq 1$ and $\alpha\in (0,1)$. There exists a positive constant $C=C(n,\tau,l,\alpha)$, independent of $k$ and $r$, such that for any $\phi \in X^{l,\alpha}_{k,s}$,
\begin{align}
C^{-1}k^{-1}\norm{(X^{l,\alpha}_{k,s})'}{\phi}&\leq\norm{X^{l,\alpha}_{k,s}}{\phi}\leq C\norm{(X^{l,\alpha}_{k,s})'}{\phi}, \label{equi norm 1}\\
C^{-1}\norm{(X^{l,\alpha}_{k,s})'}{\phi}&\leq \norm{(X^{l,\alpha}_{k,s})''}{\phi}\leq C\norm{(X^{l,\alpha}_{k,s})'}{\phi}\label{equi norm 2}. 
\end{align}
\end{lem} 
\begin{proof}
First, (\ref{equi norm 1}) is a direct corollary of estimate (\ref{step function gamma}).

To prove the first inequality in (\ref{equi norm 2}), consider $y_1,y_2\in B(x,d(x)/2)$ and $d(y_1)\geq d(y_2)$. Then for any $j=1,\cdots,k$, by (\ref{B geometry}), we can obtain
$$
d(x)=\min_{j}\<x-P^j\>\leq  C\min_{j}\<y_2-P^{j}\>=Cd(y_2).
$$
So, there is
$$
d(x)^{s+l+\alpha}\frac{|D^l\phi(y_1)-D^l\phi(y_2)|}{|y_1-y_2|^{\alpha}}\leq Cd(y_2)^{s+l+\alpha}\frac{|D^l\phi(y_1)-D^l\phi(y_2)|}{|y_1-y_2|^{\alpha}}.
$$

To prove the second inequality in (\ref{equi norm 2}), consider $d(y)\geq d(x)$. Then when $y\in B(x,d(x)/2)$, it is trivial. When $y\notin B(x,d(x)/2)$, we have
\begin{align*}
d(x)^{s+l+\alpha}\frac{|D^l\phi(x)-D^l\phi(y)|}{|x-y|^{\alpha}}&\leq C d(x)^{s+l}(|D^l\phi|(x)+|D^l\phi|(y))\\
&\leq Cd(x)^{s+l}|D^l\phi|(x)+d(y)^{s+l}|D^l\phi|(y).
\end{align*}
\end{proof}

Proposition \ref{modified norm} follows from our lemmas.
\begin{proof}
(a). Let us prove (\ref{equivlent norm}). Similar with (\ref{equi norm 2}), the classic weighted H{\"o}lder norm $\norm{C^{l,\alpha}_s}{\cdot}$ has an equivalent expression:
\begin{align*}
\norm{(C^{l,\alpha}_{s})''}{\phi}:=&\sum_{j=0}^l\sup_{x\in\r^n}\<x\>^{s+j}|D^j\phi|(x)\\
&+\sup_{x,y\in \r^n} \min\{\<x\>,\<y\>\}^{s+l+\alpha}\frac{|D^l\phi(x)-D^l\phi(y)|}{|x-y|^{\alpha}}.
\end{align*}
In other words, we have
$$
C^{-1}\norm{(C^{l,\alpha}_{s})''}{\phi}\leq \norm{C^{l,\alpha}_s}{\phi}\leq C\norm{(C^{l,\alpha}_{s})''}{\phi}.
$$

Combining (\ref{equi norm 1}) and (\ref{equi norm 2}), we have
$$
C^{-1}k^{-1}\norm{(X^{l,\alpha}_{k,s})''}{\phi}\leq\norm{X^{l,\alpha}_{k,s}}{\phi}\leq C\norm{(X^{l,\alpha}_{k,s})''}{\phi}.
$$

Since $r,\<x\>,d(x)\geq 1$, there are $d(x)\leq C\<x\>+r \leq Cr\<x\>$ and $\<x\>\leq Cd(x)+r\leq Crd(x)$. So, we can obtain that $C^{-1}r^{-1}\<x\>\leq d(x)\leq Cr\<x\>$. Then we have
$$
C^{-1}r^{-l-\alpha-s}\norm{(C^{l,\alpha}_{s})''}{\phi} \leq \norm{(X^{l,\alpha}_{k,s})''}{\phi}\leq C r^{l+\alpha+s}\norm{(C^{l,\alpha}_{s})''}{\phi}.
$$

Combining the above three equivalent relations, we can obtain (\ref{equivlent norm}). 

(b). Let us prove (\ref{D^12 embedding}). Use inequality (\ref{tau eqn2}), picking $s_1=s_2=s+1$ and $\tau=s-\frac{n-2}{2}>0$,
\begin{align*}
\int_{\r^n}|D\phi|^2&\leq C\norm{X^l_{k,s}}{\phi}^2\int_{\r^n}\sum_{i=1}^k\<x-P^i\>^{-s-1}\sum_{j=1}^k\<x-P^j\>^{-s-1}\\
&\leq Ck\left(\frac{k}{r}\right)^{s-\frac{n-2}{2}}\norm{X^l_{k,s}}{\phi}^2\int_{\r^n}\sum_{i=1}^k\<x-P^i\>^{-s-\frac{n-2}{2}-2}.
\end{align*}
Then we obtain $\norm{\mathcal{D}^{1,2} }{\phi}\leq Ck (k/r)^{\frac{1}{2}(s-\frac{n-2}{2})}\norm{X^{l,\alpha}_{k,s}}{\phi}$.

(c). Let us prove (\ref{dual space}). Use inequality (\ref{tau eqn2}), picking $s_1=s, s_2=s+2$ and $\tau=s-\frac{n-2}{2}>0$,
\begin{align*}
\left|\int_{\r^n}f\phi\right|&\leq C\norm{X^0_{k,s+2}}{f}\norm{X^0_{k,s}}{\phi}\int_{\r^n}\sum_{i=1}^k\<x-P^i\>^{-s}\sum_{j=1}^k\<x-P^j\>^{-s-2}\\
&\leq Ck\left(\frac{k}{r}\right)^{s-\frac{n-2}{2}}\norm{X^0_{k,s+2}}{f}\norm{X^0_{k,s}}{\phi}\int_{\r^n}\sum_{i=1}^k\<x-P^i\>^{-s-\frac{n-2}{2}-2}.
\end{align*}
Then we obtain $\norm{X^*}{f}\leq Ck^{2}(k/r)^{s-\frac{n-2}{2}}\norm{ Y }{f}$.

(d). Let us prove (\ref{bound of u}). For $l=0,1,2,3$, we have
\begin{equation}\label{u upper bound info}
|D^lu_{\Xi,\Lambda}|(x)\leq C\sum_i\<x-P^i\>^{-(n-2)-l}, \quad x\in \r^n.
\end{equation}

First, let us deal with the weighted sup-norm. For $l=0,1,2$, because of $s\leq n-2$, $d(x)\geq 1$ and Lemma \ref{step function lemma} with $s_1=s+l$, $s_2=n-2+l$, we have
$$
\sum_{i}\<x-P^i\>^{-(n-2)-l}\leq C d(x)^{s-(n-2)}\sum_{i}\<x-P^i\>^{-s-l}\leq C\sum_{i}\<x-P^i\>^{-s-l}.
$$
So we have $\norm{X^2_{k,s}}{u_{\Lambda,\Xi}}\leq C$.

Second, deal with the weighted H{\"o}lder semi-norm. For simplicity, write $B(x)=B(x,d(x)/2)$. For $\alpha\in (0,1)$, by the interpolation inequality 1 from Lemma \ref{holder semi norm lemma} with $\delta=1$ and (\ref{u upper bound info}), (\ref{B geometry}), we have
\begin{align*}
[D^2u_{\Xi,\Lambda}]_{\alpha,B(x)}&\leq Cd(x)^{1-\alpha}|D^3u_{\Xi,\Lambda}(y)|_{B(x)}+Cd(x)^{-\alpha}|D^2u_{\Xi,\Lambda}(y)|_{B(x)}\\
&\leq Cd(x)^{-\alpha}\sum_{i}\<x-P^i\>^{-(n-2)-2}.
\end{align*}
Then we can use Lemma \ref{step function lemma} with $s_1=s+2+\alpha$ and $s_2=(n-2)+2$ to obtain
$$
[D^2u_{\Xi,\Lambda}]_{\alpha,B(x)}\leq Cd(x)^{s-(n-2)}\sum_{i}\<x-P^i\>^{-s-2-\alpha}\leq C\sum_{i}\<x-P^i\>^{-s-2-\alpha}.
$$
So we have $[u_{\Lambda,\Xi}]_{X^{2,\alpha}_{k,s}}\leq C$.

Combining the estimates of weighted sup-norm and H{\"o}lder semi-norm, we obtain $\norm{ X }{u_{\Xi,\Lambda}}\leq C$.

(e). The proof of (\ref{bound of Z}) is the same with that of (\ref{bound of u}) because $\sum_{j,\mu} Z_{j,\mu}$ has similar estimates as (\ref{u upper bound info}). For $l=0,1,2,3$, we have
\begin{equation*}
\left|D^l\sum_{j,\mu} Z_{j,\mu}\right|(x)\leq C\sum_i\<x-P^i\>^{-(n-2)-l}, \quad x\in \r^n.
\end{equation*}
\end{proof}

\section{Linearized Elliptic Operator $F'(u_{\Xi,\Lambda})$}

Recall that given $c_0>0$, $\tau_0\in\r$, $\hat{h}$ is defined by (\ref{def of hath}). And then for any $\epsilon>0$, $k\geq 3$, recall that $t_k=e^{-k}$, and
$$
g(x)=\exp(\epsilon \hat{h}(t_kx)).
$$
Since we do estimates for a fixed integer $k$ in most of the proofs of this paper, for simplicity, we drop the subscript $k$ for the scaling parameter $t_k$ and write 
$$
t:=t_k=e^{-k}.
$$
In this case, for $j=1,\cdots,k$, the scaled points $P^{k,j}=t^{-1}\hat{P}^{k,j}$ have norm
$$
r:=|P^{k,1}|=(tk)^{-1}=e^{k}/k.
$$

Recall the nonlinear elliptic operator $F$ defined by (\ref{def of F}) and the linear operator $L_{\Xi,\Lambda}$ defined by (\ref{def of L}), for any $(\Xi,\Lambda)\in \mathcal{C}_k$:
\begin{align*}
F(v)&=-L_{g}v-n(n-2)|v|^{\frac{4}{n-2}}v\\ 
L_{\Xi,\Lambda}\phi&=-L_{g}\phi-n(n+2)u_{\Xi,\Lambda}^{\frac{4}{n-2}}\phi.
\end{align*}
In the end, recall the two Banach spaces $X, Y$ we work with from (\ref{def of X and Y}), the two subspaces $X_1,X_2\subseteq X$ from (\ref{def of X_1,2}) and the two subspaces $Y_1,Y_2\subseteq Y$ from (\ref{def of Y_1,2}).

\subsection{Boundedness}

In some proofs, we will use the following elementary Lemma, see e.g. \cite[Lemma 1.4]{li1997singularly}.
\begin{lem}\label{lem1}
Let $a>0$, $b\in \r$, $\beta>0$. There is a positive constant $C=C(\beta)$, such that\\when $0<\beta\leq 1$, we have 
$$
||a+b|^{\beta}-a^{\beta}|\leq C\min\{|b|^{\beta}, a^{\beta-1}|b|\};
$$
when $1<\beta\leq 2$, we have
$$
\left||a+b|^{\beta-1}(a+b)-a^{\beta}-\beta a^{\beta-1}b\right|\leq C\min\{|b|^{\beta}, a^{\beta-2}|b|^2\};
$$
when $2<\beta\leq 3$, we have
$$
\left||a+b|^{\beta}-a^{\beta}-\beta a^{\beta-1}b-\frac{1}{2}\beta(\beta-1)a^{\beta-2}b^2\right|\leq C\min\{|b|^{\beta}, a^{\beta-3}|b|^3\}.
$$
\end{lem}  

We have a H{\"o}lder semi-norm analogue of Lemma \ref{lem1}.
\begin{lem}\label{lem1-holder}
Let $0<\alpha\leq \beta< 1$, $\tau\geq 1$, $B\subseteq \r^n$ be an open ball with radius $R>0$, $u, \phi\in C^1(B)$, $u>0$. There is a positive constant $C=C(\alpha,\beta)$ such that 
\begin{align}
\left[ |\phi|^{\beta}\right]_{\beta, B}&\leq C|D\phi|^{\beta}_B\label{alpha=beta}\\
\left[ |\phi|^{\beta}\right]_{\alpha, B}&\leq C\left(R^{-\alpha}|\phi|^{\beta}_B+R^{\beta-\alpha}|D\phi|^{\beta}_B\right)\label{alpha<beta}\\
\left[ |u+\phi|^{\beta}-u^{\beta}\right]_{\alpha, B}&\leq C\left(R^{-\tau\alpha}+R^{\tau(\beta-\alpha)}(|Du|^{\beta}_B+|D\phi|^{\beta}_B)\right)|\phi|_B^{\beta-\alpha}\label{tricky}.
\end{align} 
\end{lem}  
\textbf{Remark}: In (\ref{alpha=beta}), $|D\phi|^{\beta}_B$ can not be replaced by $\left[\phi\right]_{\gamma,B}^{\beta}$ for any $\gamma\in (0,1)$. We can pick a sequence of functions $\phi_{\epsilon}(x)=(\epsilon+|x|^2)^{\gamma/2}$, for $\epsilon>0$ as a counterexample.

\begin{proof}
By Lemma \ref{lem1}, we know that for any $x,y\in B$, there is
$$
\frac{|\phi^{\beta}(x)-\phi^{\beta}(y)|}{|x-y|^{\beta}}\leq C\frac{|\phi(x)-\phi(y)|^{\beta}}{|x-y|^{\beta}}\leq C|D\phi|_{B}^{\beta}.
$$
So we obtain (\ref{alpha=beta}).

Using the interpolation inequality 2 from Lemma \ref{holder semi norm lemma} with $\delta=1$, we have
$$
\left[ |\phi|^{\beta}\right]_{\alpha, B}\leq C\left(R^{-\alpha}|\phi|^{\beta}_B+R^{\beta-\alpha}[|\phi|^{\beta}]_{\beta,B}\right).
$$
Combining the above estimate with (\ref{alpha=beta}), we obtain (\ref{alpha<beta}).

Using the interpolation inequality 2 from Lemma \ref{holder semi norm lemma} with $\delta>0$, we have 
$$
\left[ |u+\phi|^{\beta}-u^{\beta}\right]_{\alpha, B}\leq C\left((\delta R)^{-\alpha}\left| |u+\phi|^{\beta}-u^{\beta}\right|_{B}+(\delta R)^{\beta-\alpha}\left[ |u+\phi|^{\beta}-u^{\beta}\right]_{\beta, B}\right).
$$
Then we use Lemma \ref{lem1} for the first term and triangle inequality from Lemma \ref{holder semi norm lemma} for the second term to obtain that
\begin{equation*}
\left[ |u+\phi|^{\beta}-u^{\beta}\right]_{\alpha, B}\leq C\left((\delta R)^{-\alpha}\left| \phi\right|_{B}^{\beta}+(\delta R)^{\beta-\alpha} \left(\left[|u+\phi|^{\beta}\right]_{\beta, B}+\left[u^{\beta}\right]_{\beta, B}\right)\right).  
\end{equation*} 
Then combining the above estimate with (\ref{alpha=beta}), we have 
\begin{equation*}
\left[ |u+\phi|^{\beta}-u^{\beta}\right]_{\alpha, B}\leq C\left((\delta R)^{-\alpha}\left| \phi\right|_{B}^{\beta}+(\delta R)^{\beta-\alpha} (|Du|^{\beta}_B+|D\phi|^{\beta}_B)\right).
\end{equation*}
Pick $\delta=R^{\tau-1}|\phi|_{B}$. We can obtain (\ref{tricky}).
\end{proof}

\begin{lem}
Let $k\geq 3$, $0<\tau<1$, $t=e^{-k}$ and $r=(tk)^{-1}$. There exists a positive constant $C=C(\tau)$ such that 
\begin{equation}\label{B inside outside argument}
td(x)\leq C \text{ when } x\in \{B(0,t^{-1})\cap B(x,\tau d(x))\neq \emptyset\}.
\end{equation}
\end{lem}
\begin{proof}
By the definition of $d(x)$, we know that for any $x\in \r^n$,
$$
td(x)\leq t(1+|x|+r)= t+k^{-1}+t|x|.
$$
For $x\in \{B(0,t^{-1})\cap B(x,\tau d(x))\neq \emptyset\}$, there is
$$
|x|<t^{-1}+\tau d(x).
$$
Combining the above two estimates together, we have
$$
td(x)\leq t+k^{-1}+1+\tau td(x).
$$
Then we can obtain (\ref{B inside outside argument}) since $\tau <1$.
\end{proof}

Now, we can give the main proposition of this subsection.
\begin{prop}\label{Welldefined}
Let $n\geq 5$, $0<\alpha<\frac{4}{n-2}$, $\frac{n-2}{2}<s<n-2$, $c_0>0$, and $\tau_0\in \r$. For any $0<\epsilon<1$, $k\geq 3$, $F$ is a $C^1$ map from $ X $ into $ Y $. Moreover, there exists a positive constant $C=C(n,s,\alpha,c_0,\tau_0)$, such that for any $(\Xi,\Lambda)\in \mathcal{C}_k$, any $v, \phi_1,\phi_2\in  X $ with $\norm{ X }{\phi_1}, \norm{ X }{\phi_2}\leq 1$,
\begin{align}
\norm{ Y }{F(v)}&\leq C\left(1+\epsilon+\norm{ X }{v}^{\frac{4}{n-2}}\right)\norm{ X }{v};\label{Welldef of F}\\
\norm{ Y }{L_{\Xi,\Lambda}v}&\leq C\left(1+\epsilon+\norm{ X }{u_{\Xi,\Lambda}}^{\frac{4}{n-2}}\right)\norm{ X }{v};\label{boundedness of L}
\end{align}
and 
\begin{equation}\label{Lip of F}
\begin{aligned}
    &\norm{ Y }{F(u_{\Xi,\Lambda}+\phi_1)-F(u_{\Xi,\Lambda}+\phi_2)-L_{\Xi,\Lambda}(\phi_1-\phi_2)}\\
    \leq &C\left(\norm{ X }{\phi_1}+\norm{ X }{\phi_2}\right)^{\frac{4}{n-2}-\alpha}\norm{ X }{\phi_1-\phi_2}.
\end{aligned}
\end{equation}
Moreover, $L_{\Xi,\Lambda}$ is the linearized operator $F'(u_{\Xi,\Lambda})$.
\end{prop}
\textbf{Remark.} In this proposition, $s>\frac{n-2}{2}$ is used to prove (\ref{Welldef of F}) while both $s>\frac{n-2}{2}$ and $s<n-2$ are used to prove (\ref{boundedness of L}) and (\ref{Lip of F}).

\begin{proof}
(a). Let us prove $\norm{ Y }{F(v)}<\infty$ for any $v\in  X $. In this proof, we simply write $B(x)=B(x,d(x)/2)$. We write $F(v)$ as:
$$
F(v)=-(\RN{1}+\RN{2}+\RN{3}), 
$$
where $\RN{1}=\Delta v$, $\RN{2}=(L_{g }-\Delta) v$, and $\RN{3}=n(n-2)|v|^{\frac{4}{n-2}}v$. Recalling that $Y=X^{0, \alpha}_{k,s+2}$, we need to estimate both the weighted sup-norm $\norm{X^{0}_{k,s+2}}{\cdot}$ and weighted H{\"o}lder semi-norm $[\cdot]_{X^{0,\alpha}_{k,s+2}}$ of $\RN{1}, \RN{2}, \RN{3}$.

We obtain $\norm{ Y }{\RN{1}}\leq \norm{ X }{v}$ directly by the definition of $v\in X$.

Next, we do the estimates of $|\RN{2}|(x)$ and $[\RN{2}]_{\alpha,B(x)}$. Because $\supp h=B(0,t^{-1})$, we know that $|\RN{2}|(x)=[\RN{2}]_{\alpha,B(x)}=0$ when $x\in \{B(0,t^{-1})\cap B(x)= \emptyset\}$. So, we just need work with $x\in \{B(0,t^{-1})\cap B(x)\neq \emptyset\}$. By (\ref{B inside outside argument}), we know that $td(x)\leq C$.

Write 
\begin{equation}\label{coefficient}
\RN{2}=a^{\mu\nu}D_{\mu\nu}v+b^{\mu}D_{\mu}v+cv,    
\end{equation}
where $a^{\mu\nu}=g^{\mu\nu}-\delta^{\mu\nu}$, $b^{\mu}=D_{\nu}g^{\mu\nu}+\frac{1}{2}g^{\alpha\beta}D_{\nu}g_{\alpha\beta}g^{\mu\nu}$, and $c=-c(n)R_g$. 

For $x\in \r^n$, $l=0,1,2$, because $v\in X$ and (\ref{B geometry}), we know that  
\begin{equation}\label{translation of the sup norm}
|D^lv|_{B(x)}\leq \norm{ X }{v}\left|\sum_{i=1}^k\<y-P^i\>^{-s-l}\right|_{B(x)}\leq C\norm{ X }{v}\sum_{i=1}^k\<x-P^i\>^{-s-l}.    
\end{equation}
Using the interpolation inequality 1 from Lemma \ref{holder semi norm lemma} with $\delta=1$, combining with (\ref{translation of the sup norm}) and $d(x)\leq \<x-P^i\>$, we obtain the following summary: for $x\in \r^n$, $l=0,1,2$,
\begin{equation}\label{translate of the holder norm}
[D^lv]_{\alpha,B(x)}\leq C\norm{ X }{v}d(x)^{-\alpha}\sum_{i=1}^k\<x-P^i\>^{-s-l}.
\end{equation}

Recalling the $g(x)=\hat{g}(tx)$, we can similarly define $\hat{a}^{\mu\nu}, \hat{b}^{\mu},\hat{c}$. Then by scaling, $\hat{g}=\exp(\epsilon\hat{h})$ and $\norm{C^3(B(0,1))}{\hat{h}}<1$, we know that
\begin{equation}\label{scaling coeifficients}
\begin{aligned}
|a^{\mu\nu}|_{B(x)}&=|\hat{a}^{\mu\nu}|_{tB(x)}\leq C\epsilon,\\
|b^{\mu}|_{B(x)}&=t|\hat{b}^{\mu}|_{tB(x)}\leq C\epsilon t,\\
|c|_{B(x)}&=t^2|\hat{c}|_{tB(x)}\leq C\epsilon t^2, 
\end{aligned}
\qquad
\begin{aligned}
\relax[a^{\mu\nu}]_{\alpha,B(x)}&=t^{\alpha}[\hat{a}^{\mu\nu}]_{\alpha,tB(x)}\leq C\epsilon t^{\alpha};\\
[b^{\mu}]_{\alpha,B(x)}&=t^{1+\alpha}[\hat{b}^{\mu}]_{\alpha,tB(x)}\leq C\epsilon t^{1+\alpha};\\
[c]_{\alpha,B(x)}&=t^{2+\alpha}[\hat{c}]_{\alpha,tB(x)}\leq C\epsilon t^{2+\alpha}.
\end{aligned}
\end{equation}

Now, we are ready to do the estimates of $|\RN{2}|(x)$ and $[\RN{2}]_{\alpha,B(x)}$ for $x\in \{x\in\r^n\mid B(x)\cap B(0,t^{-1})\neq \emptyset\}$. \\
For $|\RN{2}|(x)$, according to (\ref{translation of the sup norm}), (\ref{scaling coeifficients}), we can obtain
$$
|\RN{2}|(x)\leq C\sum_{l=0}^2\epsilon t^l\norm{ X }{v}\sum_{i=1}^k\<x-P^i\>^{-s-2+l}.
$$
Then by Lemma \ref{step function lemma} and $td(x)\leq C$, we have that
$$
\left(\sum_{i=1}^k\<x-P^i\>^{-s-2}\right)^{-1}|\RN{2}|(x)\leq C\sum_{l=0}^2\epsilon t^ld(x)^l\norm{ X }{v}\leq C\epsilon\norm{ X }{v}.
$$
For $[\RN{2}]_{\alpha,B(x)}$, using the product rule from Lemma \ref{holder semi norm lemma} and (\ref{translation of the sup norm}), (\ref{translate of the holder norm}), (\ref{scaling coeifficients}), we can obtain
$$
[\RN{2}]_{\alpha,B(x)}\leq C\sum_{l=0}^2\epsilon (t^ld(x)^{-\alpha}+t^{l+\alpha})\norm{ X }{v}\sum_{i=1}^k\<x-P^i\>^{-s-2+l}.
$$
Then by Lemma \ref{step function lemma} and $td(x)\leq C$, we have that
$$
\bigg(\sum_{i=1}^k\<x-P^i\>^{-s-2-\alpha}\bigg)^{-1}[\RN{2}]_{\alpha,B(x)}\leq C\sum_{l=0}^2\epsilon (t^ld(x)^{l}+t^{l+\alpha}d(x)^{l+\alpha})\norm{ X }{v}\leq C\epsilon\norm{ X }{v}.
$$
So, we obtain $\norm{ Y }{\RN{2}}\leq C\epsilon\norm{ X }{v}$.

The estimates of the nonlinear term $\RN{3}$ are tricky. For $|\RN{3}|(x)$, by (\ref{translation of the sup norm}) and Lemma \ref{step function lemma}, we obtain
\begin{equation}\label{nonlinear proof}
\begin{aligned}
\norm{X^{0}_{k,s+2}}{\RN{3}}&\leq C\norm{ X }{v}^{\frac{n+2}{n-2}}\sup_{x\in\r^n}\left(\sum_{i=1}^k\<x-P^i\>^{-s-2}\right)^{-1}\left(\sum_{i=1}^k\<x-P^i\>^{-s}\right)^{1+\frac{4}{n-2}},\\
&\leq C\norm{ X }{v}^{\frac{n+2}{n-2}}\sup_{x\in\r^n}d(x)^2\left(\sum_{i=1}^k\<x-P^i\>^{-s}\right)^{\frac{4}{n-2}}.
\end{aligned}
\end{equation}
Also, by Lemma \ref{lem3} with $\tau=s-\frac{n-2}{2}>0$ and the relation $r=e^{k}/k$, we obtain
\begin{equation}\label{nonlinear proof 2}
\sum_{i=1}^k\<x-P^i\>^{-s}\leq Ck(k/r)^{s-\frac{n-2}{2}}d(x)^{-\frac{n-2}{2}}\leq Cd(x)^{-\frac{n-2}{2}}, \text{ for } x\in \r^n.
\end{equation}
Combining (\ref{nonlinear proof}) and (\ref{nonlinear proof 2}), we have that $\norm{X^{0}_{k,s+2}}{\RN{3}}\leq C\norm{ X }{v}^{\frac{n+2}{n-2}}$.

For $[\RN{3}]_{\alpha,B(x)}$, using the interpolation inequality 1 with $\delta=1$, we have
$$
[\RN{3}]_{\alpha,B(x)}\leq C\left(d(x)^{-\alpha}|v|_{B(x)}+d(x)^{1-\alpha}|Dv|_{B(x)}\right)|v|^{\frac{4}{n-2}}_{B(x)}.
$$
Then according to (\ref{translate of the holder norm}) and Lemma \ref{step function lemma}, we obtain that
\begin{equation}\label{nonlinear proof holder}
\begin{aligned}
\relax[\RN{3}]_{X^{0,\alpha}_{k,s+2}}&\leq C\norm{ X }{v}^{\frac{n+2}{n-2}}\sup_{x\in\r^n}d(x)^{-\alpha}\bigg(\sum_{i=1}^k\<x-P^i\>^{-s-2-\alpha}\bigg)^{-1}\bigg(\sum_{i=1}^k\<x-P^i\>^{-s}\bigg)^{\frac{n+2}{n-2}},\\
&\leq C\norm{ X }{v}^{\frac{n+2}{n-2}}\sup_{x\in\r^n}d(x)^2\bigg(\sum_{i=1}^k\<x-P^i\>^{-s}\bigg)^{\frac{4}{n-2}}.
\end{aligned}
\end{equation}
Then combining (\ref{nonlinear proof holder}) and (\ref{nonlinear proof 2}), we obtain $[\RN{3}]_{X^{0,\alpha}_{k,s+2}}\leq C\norm{ X }{v}^{\frac{n+2}{n-2}}$.

Combining the estimates of $\RN{1}, \RN{2}, \RN{3}$ together, we obtain the estimate (\ref{Welldef of F}). Consequently, for any $v\in X$, $F(v)\in Y$.

Moreover, recall that for any $v\in X$, $L_{\Xi,\Lambda}v=-(\RN{1}+\RN{2})-n(n+2)u_{\Xi,\Lambda}^{\frac{4}{n-2}}v$. So the estimates of $\RN{1}, \RN{2}$ still works. By (\ref{bound of u}), we know that $u_{\Xi,\Lambda}\in X$. We can do similar arguments as in (\ref{nonlinear proof}) and (\ref{nonlinear proof holder}) to obtain
$$
\norm{ Y }{u_{\Xi,\Lambda}^{\frac{4}{n-2}}v}\leq C\norm{ X }{u_{\Xi,\Lambda}}^{\frac{4}{n-2}}\norm{ X }{v}.
$$
Combining the above estimate with the estimates of $\RN{1}, \RN{2}$, we can obtain (\ref{boundedness of L}).

(b). We will prove the estimate (\ref{Lip of F}). For simplicity, in the rest of this proof, we denote $u=u_{\Xi,\Lambda}$. First, we can use the mean value theorem to obtain
\begin{align*}
&F(u+\phi_1)-F(u+\phi_2)-L_{\Xi,\Lambda}(\phi_1-\phi_2)\\
=&\frac{n+2}{n-2}\left(\int_{0}^{1}|u+\theta\phi_1+(1-\theta)\phi_2|^{\frac{4}{n-2}}-u^{\frac{4}{n-2}}d\theta\right)(\phi_1-\phi_2).
\end{align*}
Since Lemma \ref{lem1}, we know that for any $\theta\in (0,1)$,
$$
\left||u+\theta\phi_1+(1-\theta)\phi_2|^{\frac{4}{n-2}}-u^{\frac{4}{n-2}}\right|\leq C\left(|\phi_1|^{\frac{4}{n-2}}+|\phi_2|^{\frac{4}{n-2}}\right).
$$
Because $u,\phi_1,\phi_2\in X$, the estimates (\ref{translation of the sup norm}) also work for $u,\phi_1,\phi_2$. So, the proof of sup-norm part is similar to (\ref{nonlinear proof}) and (\ref{nonlinear proof 2}). We can obtain
\begin{equation}\label{nonlinear 1}
\norm{X^0_{k,s+2}}{F(u+\phi_1)-F(u+\phi_2)-L_{\Xi,\Lambda}(\phi_1-\phi_2)}\leq C\left(\norm{ X }{\phi_1}+\norm{ X }{\phi_2}\right)^{\frac{4}{n-2}}\norm{ X }{\phi_1-\phi_2}.    
\end{equation}

By mean value theorem and the product rule from Lemma \ref{holder semi norm lemma}, we have
$$
[F(u+\phi_1)-F(u+\phi_2)-L_{\Xi,\Lambda}(\phi_1-\phi_2)]_{X^{0,\alpha}_{k,s+2}}\leq C\sup_{x\in\r^n}(\RN{4}+\RN{5}),
$$
where
\begin{align*}
\RN{4}&:=\left(\sum_{i=1}^k\<x-P^i\>^{-s-2-\alpha}\right)^{-1}\left||u+\theta\phi_1+(1-\theta)\phi_2|^{\frac{4}{n-2}}-u^{\frac{4}{n-2}}\right|_{B(x)}[\phi_1-\phi_2]_{\alpha,B(x)},\\
\RN{5}&:=\left(\sum_{i=1}^k\<x-P^i\>^{-s-2-\alpha}\right)^{-1}\left[|u+\theta\phi_1+(1-\theta)\phi_2|^{\frac{4}{n-2}}-u^{\frac{4}{n-2}}\right]_{\alpha,B(x)}|\phi_1-\phi_2|_{B(x)}.
\end{align*}

Again, because of Lemma \ref{lem1}, the estimate of $\RN{4}$ is similar to (\ref{nonlinear proof}) and (\ref{nonlinear proof 2}). We can obtain 
\begin{equation}\label{nonlinear 2}
\RN{4}\leq C\left(\norm{ X }{\phi_1}+\norm{ X }{\phi_2}\right)^{\frac{4}{n-2}}\norm{ X }{\phi_1-\phi_2}.    
\end{equation}

We use (\ref{tricky}) from Lemma \ref{lem1-holder} with $\tau=\frac{n-2}{2}+1$, (\ref{translate of the holder norm}) and $\norm{ X }{\phi_1}, \norm{ X }{\phi_1}, \norm{ X }{ u}\leq C$ to obtain
\begin{align*}
&\left[ |u+\theta\phi_1+(1-\theta)\phi_2|^{\frac{4}{n-2}}-u^{\frac{4}{n-2}}\right]_{\alpha, B(x)}\\
\leq &Cd(x)^{-\alpha}\left((\norm{ X }{\phi_1}+\norm{ X }{\phi_2})\left(\sum_{i=1}^k\<x-P^i\>^{-s}\right)\right)^{\frac{4}{n-2}-\alpha}\\
&\cdot\left(d(x)^{-\alpha\frac{n-2}{2}}+d(x)^{(\frac{4}{n-2}-\alpha)\frac{n-2}{2}}\left(d(x)\sum_{i=1}^k\<x-P^i\>^{-s-1}\right)^{\frac{4}{n-2}}\right).
\end{align*}
Combining the above estimate with $d(x)\leq \<x-P^i\>$ and (\ref{nonlinear proof 2}), we have
$$
\left[ |u+\theta\phi_1+(1-\theta)\phi_2|^{\frac{4}{n-2}}-u^{\frac{4}{n-2}}\right]_{\alpha, B(x)}\leq C\left(\norm{ X }{\phi_1}+\norm{ X }{\phi_2}\right)^{\frac{4}{n-2}-\alpha}d(x)^{-2-\alpha}.
$$
Combining the above estimate with (\ref{translate of the holder norm}) and Lemma \ref{step function lemma}, we obtain
\begin{equation}\label{nonlinear 3}
\RN{5}\leq C\left(\norm{ X }{\phi_1}+\norm{ X }{\phi_2}\right)^{\frac{4}{n-2}-\alpha}\norm{ X }{\phi_1-\phi_2}.
\end{equation}

Combining (\ref{nonlinear 1}), (\ref{nonlinear 2}), and (\ref{nonlinear 3}), we obtain estimate (\ref{Lip of F}). Consequently, $F: X\rightarrow Y$ is $C^1$ map with $F'(u)=L_{\Xi,\Lambda}$.

\end{proof}

\subsection{Schauder Estimates}

As the second part of this section, we give the following $k$-independent Schauder estimate of $L_{\Xi,\Lambda}$.
\begin{prop}
Let $n\geq 5$, $0<\alpha<\frac{4}{n-2}$, $\frac{n-2}{2}<s<n-2$, $c_0>0$, and $\tau_0\in \r$. Then there exists $C=C(n,s,\alpha,c_0,\tau_0)$, such that for $0<\epsilon<1$, $k\geq 3$, $(\Xi,\Lambda)\in \mathcal{C}_k$ and $\phi\in  X $,
\begin{equation}\label{schauder of Laplace}
\norm{ X }{\phi}\leq C(\norm{ Y }{\Delta\phi}+\norm{X^0_{k,s}}{\phi})
\end{equation}
and
\begin{equation}\label{schauder of L}
\norm{ X }{\phi}\leq C(\norm{ Y }{L_{\Xi,\Lambda}\phi}+\norm{X^0_{k,s}}{\phi})
\end{equation}
\end{prop}
\textbf{Remark.} In this proposition, $\frac{n-2}{2}<s<n-2$ is needed for estimate (\ref{schauder of L}). $s>1$ is enough for (\ref{schauder of Laplace}).

\begin{proof}
The proof of estimate (\ref{schauder of Laplace}) is similar and simpler than that of (\ref{schauder of L}). So here we just give the proof of (\ref{schauder of L}). We write $B_1(x):=B(x,d(x)/2)$ and $B_2(x):=B(x,3d(x)/4)$. Recall 
$$
L_{\Xi,\Lambda}=-\Delta-(a^{\mu\nu}D_{\mu\nu}+b^{\mu}D_{\mu}+c)-n(n+2)u_{\Xi,\Lambda}^{\frac{4}{n-2}},
$$
where $a^{\mu\nu}=g^{\mu\nu}-\delta^{\mu\nu}$, $b^{\mu}=D_{\nu}g^{\mu\nu}+\frac{1}{2}g^{\alpha\beta}D_{\nu}g_{\alpha\beta}g^{\mu\nu}$, and $c=c(n)R_g$.

In order to use interior Schauder estimate, we need to estimate the local H{\"o}lder norm $|\cdot|'_{\alpha,B_2(x)}$ of the coefficients of $L_{\Xi,\Lambda}$. Because $\supp h=B(0,t^{-1})$, we know that $a^{\mu\nu},b^{\mu},c=0$ in $B_2(x)$ when $x\in \{B(0,t^{-1})\cap B_2(x)= \emptyset\}$. So, next we just need work with $x\in\{B(0,t^{-1})\cap B_2(x)\neq \emptyset\}$. By (\ref{scaling coeifficients}) and (\ref{B inside outside argument}), we obtain
\begin{equation}\label{schauder 2}
\begin{aligned}
&|a^{\mu\nu}|'_{\alpha,B_2(x)}+d(x)|b^{\mu}|'_{\alpha,B_2(x)}+d(x)^2|c|'_{\alpha,B_2(x)}\\
\leq& C\sum_{l=0}^2t^ld(x)^l+C\sum_{l=0}^2t^{l+\alpha}d(x)^{l+\alpha}\leq C.  
\end{aligned}
\end{equation}

Second, let us control the last term $u_{\Xi,\Lambda}^{\frac{4}{n-2}}$. We use (\ref{alpha<beta}) from Lemma \ref{lem1-holder} to obtain
$$
|u^{\frac{4}{n-2}}_{\Xi,\Lambda}|'_{\alpha,B_2(x)}\leq C|u_{\Xi,\Lambda}|^{\frac{4}{n-2}}_{B_2(x)}+Cd(x)^{\frac{4}{n-2}}|Du_{\Xi,\Lambda}|^{\frac{4}{n-2}}_{B_2(x)}.
$$
Then because $u_{\Xi,\Lambda}\in X$ and (\ref{B geometry}), we have
$$
|u^{\frac{4}{n-2}}_{\Xi,\Lambda}|'_{\alpha,B_2(x)}\leq C\left(\sum_{i=1}^k\<x-P^i\>^{-s}\right)^{\frac{4}{n-2}}+C\left(d(x)\sum_{i=1}^k\<x-P^i\>^{-s-1}\right)^{\frac{4}{n-2}}.
$$
Then because $d(x)\leq \<x-P^i\>$ and (\ref{nonlinear proof 2}), we can obtain
\begin{equation}\label{schauder 1}
d(x)^2|u^{\frac{4}{n-2}}_{\Xi,\Lambda}|'_{\alpha,B_2(x)}\leq C.
\end{equation}

So, combining estimates (\ref{schauder 1}), (\ref{schauder 2}) and $|\delta^{\mu\nu}|'_{\alpha,B_2(x)}=1$, we can obtain the interior Schauder estimate: there exists a positive constant $C$ independent of the position $x\in \r^n$, such that:
\begin{equation}\label{interior schauder}
|\phi|'_{2,\alpha,B_1(x)}\leq C(d(x)^2|L_{\Xi,\Lambda}\phi|'_{\alpha,B_2(x)}+|\phi|_{B_2(x)}).
\end{equation}

We use the first inequality in (\ref{step function gamma}) to obtain that for any $x\in \r^n$, 
\begin{equation}\label{lower bound}
\begin{aligned}
&\sum_{l=0}^2\left(\sum_{i=1}^k\<x-P^i\>^{-s-l}\right)^{-1}|D^l\phi(x)|+\left(\sum_{i=1}^k\<x-P^i\>^{-s-2-\alpha}\right)^{-1}[D^2\phi]_{\alpha,B_1(x)}\\
\leq &C\gamma(d(x))^{-1}d(x)^{s}|\phi|'_{2,\alpha,B_1(x)}.
\end{aligned}
\end{equation}
We use the second inequality in (\ref{step function gamma}) to obtain that for any $x\in \r^n$,
\begin{equation}\label{upper bound}
\begin{aligned}
&\gamma(d(x))^{-1}d(x)^{s}\left(d(x)^2|L_{\Xi,\Lambda}\phi|'_{\alpha,B_2(x)}+|\phi|_{B_2(x)}\right)
\\
\leq &C\bigg(\sum_{i=1}^k\<x-P^i\>^{-s-2}\bigg)^{-1}|L_{\Xi,\Lambda}\phi|_{B_2(x)}+C\bigg(\sum_{i=1}^k\<x-P^i\>^{-s-2-\alpha}\bigg)^{-1}[L_{\Xi,\Lambda}\phi]_{\alpha,B_2(x)}\\
+&C\bigg(\sum_{i=1}^k\<x-P^i\>^{-s}\bigg)^{-1}|\phi|_{B_2(x)}.
\end{aligned}
\end{equation}

Assume that $y\in \overline{B_2(x)}$ such that $|\phi|(y)=|\phi|_{B_2(x)}$. Since (\ref{B geometry}), we have 
$$
\left(\sum_{i=1}^k\<x-P^i\>^{-s}\right)^{-1}|\phi|_{B_2(x)}\leq C\left(\sum_{i=1}^k\<y-P^i\>^{-s}\right)^{-1}|\phi|(y)\leq C\norm{X^0_{k,s}}{\phi}.
$$
Similarly, we have 
$$
\left(\sum_{i=1}^k\<x-P^i\>^{-s-2}\right)^{-1}|L_{\Xi,\Lambda}\phi|_{B_2(x)}\leq C\norm{X^{0}_{k,s+2}}{L_{\Xi,\Lambda}\phi}.
$$

Also, assume that $y,\bar{y}\in \overline{B_2(x)}$ such that $\frac{|L_{\Xi,\Lambda}\phi(y)-L_{\Xi,\Lambda}\phi(\bar{y})|}{|y-\bar{y}|^{\alpha}}=[L_{\Xi,\Lambda}\phi]_{\alpha,B_2(x)}$. Then we can pick a sequence of points $\{y_j\}$ on the line segment connecting $y, \bar{y}$ such that $y_1=y$ and $|y_j-y_{j+1}|=d(y_j)/2$.
Then because for any $y_j\in \overline{B_2(x)}$, by (\ref{B geometry}), there is $d(y_j)/2\geq d(x)/16$. So, $|y_j-y_{j+1}|\geq d(x)/16$. Combining with $|y-\bar{y}|\leq 3d(x)/2$, there exists $l\leq 24$ such that, $\bar{y}\in B_1(y_l)$. Then we can pick $1\leq \bar{j}\leq l$ such that $[L_{\Xi,\Lambda}\phi]_{\alpha,B_1(y_{\bar{j}})}=\max_{1\leq j\leq l}[L_{\Xi,\Lambda}\phi]_{\alpha,B_1(y_{j})}$. So, combining with (\ref{B geometry}), the H{\"o}lder semi-norm term
\begin{align*}
&\left(\sum_{i=1}^k\<x-P^i\>^{-s-2-\alpha}\right)^{-1}[L_{\Xi,\Lambda}\phi]_{\alpha,B_2(x)}\\
\leq &C\left(\sum_{i=1}^k\<y_{\bar{j}}-P^i\>^{-s-2-\alpha}\right)^{-1}[L_{\Xi,\Lambda}\phi]_{\alpha,B_1(y_{\bar{j}})}\leq C[L_{\Xi,\Lambda}\phi]_{X^{0,\alpha}_{k,s+2}}.
\end{align*}

So, combining the above three estimates with (\ref{upper bound}), we obtain that for any $x\in \r^n$,
\begin{equation}\label{upper bound again}
\gamma(d(x))^{-1}d(x)^{s}\left(d(x)^2|L_{\Xi,\Lambda}\phi|'_{\alpha,B_2(x)}+|\phi|_{B_2(x)}\right)
\leq C(\norm{Y}{L_{\Xi,\Lambda}\phi}+\norm{X^0_{k,s}}{\phi}).
\end{equation}

Combining (\ref{interior schauder}), (\ref{lower bound}), and (\ref{upper bound again}), we obtain that for any $x\in \r^n$,
\begin{align*}
&\sum_{l=0}^2\left(\sum_{i=1}^k\<x-P^i\>^{-s-l}\right)^{-1}|D^l\phi(x)|+\left(\sum_{i=1}^k\<x-P^i\>^{-s-2-\alpha}\right)^{-1}[D^2\phi]_{\alpha,B_1(x)}\\
\leq &C(\norm{Y}{L_{\Xi,\Lambda}\phi}+\norm{X^0_{k,s}}{\phi}).
\end{align*}
In the end, we can obtain (\ref{schauder of L}) because $x\in \r^n$ is arbitrary and the positive constant $C$ above does not depend on $x$. 
\end{proof}

The following Lemma is well known, see e.g. \cite[Lemma B.2]{wei2010infinitely}, 
\begin{lem}\label{lem4}
Let $n\geq 3$, $0<s<n-2$, there exists $C=C(n,s)>0$ such that 
$$
\int_{\r^n}\frac{\<x\>^{-s-2}}{|y-x|^{n-2}}dx\leq C\<y\>^{-s}.
$$
\end{lem}

A corollary of the $k$-independent Schauder estimate and the Lemma above is the boundedness of the Newtonian potential $\Gamma(f):=\int_{\r^n}|x-y|^{-(n-2)}f(y)dy$.
\begin{coro}\label{Gammaboundedness}
Let $n\geq 5$, $1<s<n-2$, $\alpha\in(0,1)$. Then $\Gamma$ is a bounded operator from $ Y $ to $ X $. Moreover, there exists $C=C(n,s,\alpha)>0$, independent of $k$, such that for any $f\in  Y $,
$$
\norm{ X }{\Gamma(f)}\leq C\norm{ Y }{f}.
$$
\end{coro}
\begin{proof}
According to Schauder estimate (\ref{schauder of Laplace}), we only need to prove that $$
\norm{X^0_{k,s}}{\Gamma(f)}\leq C\norm{ Y }{f}.
$$
The proof is straightforward with the help of Lemma \ref{lem4}. 
\begin{align*}
\left|\int_{\r^n}\frac{f(y)}{|x-y|^{n-2}}dy\right|
&=\norm{X^0_{k,s+2}}{f}\sum_{i=1}^k\int_{\r^n}\frac{\<y-P^i\>^{-s-2}}{|x-y|^{n-2}}dy\\
&\leq C\norm{X^0_{k,s+2}}{f}\sum_{i=1}^k\<x-P^i\>^{-s}.
\end{align*}
\end{proof}

\subsection{Invertibility}

Recall the approximate solutions $u_{\Xi,\Lambda}(x)=\sum_{i=1}^k\sigma_{P^i+\xi^i,\lambda_i}(x)$. To simplify the notation, we denote
$$
\sigma_i(x):=\sigma_{P^i+\xi^i,\lambda_i}(x).
$$ 
\begin{prop}\label{Surjectivity}
Let $n\geq 6$, $0<\alpha<\frac{4}{n-2}$, $\frac{n-2}{2}<s<n-2$, $c_0>0$, and $\tau_0\in \r$. There exist positive constants $\epsilon_0, k_0, C$ depending on $ n,s,\alpha, c_0, \tau_0$ such that for any $0<\epsilon\leq \epsilon_0$, $k\geq k_0$, $(\Xi,\Lambda)\in \mathcal{C}_k$, we have that $Y_1=\Delta(X_1)$ and $Y_2=L_{\Xi,\Lambda}(X_2)$ are closed in $ Y $ and there is a decomposition 
\begin{equation*}
Y =Y_1\oplus Y_2.
\end{equation*}
Also, $\bar{L}: X_2\rightarrow Y_2$ is invertible. Moreover,
\begin{align}
\norm{ X }{\bar{L}^{-1}f_2}&\leq C\norm{ Y }{f_2},\text{ for any }f_2\in Y_2, \label{invertibility}\\
\norm{ Y }{\Proj_{Y_2}f}&\leq C\norm{ Y }{f}, \text{ for any }f\in  Y \label{bound of proj}.
\end{align}
\end{prop}

Before the proof of Proposition \ref{Surjectivity}, we give the following three lemmas. 

\textbf{Remark.} $\frac{n-2}{2}<s<n-2$ will be used to prove Lemmas \ref{smallnessfu} and \ref{smallnessZ}. In addition, Schauder estimate (\ref{schauder of L}) and Corollary \ref{Gammaboundedness}, which need the condition on $s$, will be used in the proof of Proposition \ref{Surjectivity} as well. 

The first lemma is directly from Bahri's book \cite[Section 2.1, F1-F6]{bahri1989critical}.
\begin{lem}
We have $C=C(n)$, such that for $i,j=1,\cdots,k$ and $\mu,\nu=1,\cdots,n$, 
\begin{equation}\label{orthognality of Z}
\begin{aligned}
&\left|\<Z_{i,0},Z_{j,0}\>_{\mathcal{D}^{1,2} }-a_1\lambda_i^{-2}\delta_{ij}\right|\leq C(k/r)^{n-2};\\
&\left|\<Z_{i,\mu},Z_{j,\nu}\>_{\mathcal{D}^{1,2} }-a_2\lambda_i^{2}\delta_{ij}\delta_{\mu\nu}\right|\leq C(k/r)^{n-2};\\
&\left|\<Z_{i,\mu},Z_{j,0}\>_{\mathcal{D}^{1,2} }\right|\leq C(k/r)^{n-2},
\end{aligned}
\end{equation}
where $a_1=n(n+2)(n-2)^2\int_{\r^n}\frac{1+(1+|x|^2)^{2}}{4(1+|x|^2)^{n+2}}$ and $a_2=n(n+2)(n-2)^2\int_{\r^n}\frac{x_1^2}{(1+|x|^2)^{n+2}}$.
\end{lem}

Next, we will prove one important ingredient of Lyapunov-Schmidt reduction: the smallness. 
\begin{lem}\label{smallnessfu}
Let $n\geq 6$, $\frac{n-2}{2}<s< n-2$, $0<\alpha<\frac{4}{n-2}$, $c_0>0$, $\tau_0\in\r$. Then there exists $C=C(n,s,\alpha,c_0,\tau_0)$, such that for any $0<\epsilon<1$, $k\geq 3$, $(\Xi,\Lambda)\in\mathcal{C}_k$, 
\begin{equation}\label{smallness1}
\norm{ Y }{F(u_{\Xi,\Lambda})}\leq C\left(\epsilon+(k/r)^{n-s}\right).
\end{equation}
\end{lem}
\textbf{Remark.} The Lemma above is enough to be used to prove Proposition \ref{Surjectivity}. A stronger estimate of $\norm{ Y }{F(u_{\Xi,\Lambda})}$ can be seen in Lemma \ref{basicproperty2}.

\begin{proof}
For simplicity, we will drop the subscript $\Xi,\Lambda$ of $u_{\Xi,\Lambda}$ and write $u=u_{\Xi,\Lambda}$ in this proof. We use (\ref{def of F}) and $-\Delta \sigma_i=n(n-2)\sigma_i^{\frac{n+2}{n-2}}$ for any $i=1,\cdots,k$ to obtain 
$$
F(u)=-(L_{g }-\Delta)u+n(n-2)\left(\sum_{i=1}^k\sigma_i^{\frac{n+2}{n-2}}-(\sum_{i=1}^k\sigma_i)^{\frac{n+2}{n-2}}\right).
$$
We can use Proposition \ref{Welldefined} and (\ref{bound of u}) to deal with the linear part
\begin{equation}\label{smallnesslinear}
\norm{ Y }{(L_{g }-\Delta)u}\leq C\epsilon\norm{ X }{u}\leq C\epsilon.
\end{equation}

For the nonlinear part, because of the symmetry of $d(x)$ and $\sum_{i}\<x-P^i\>^{-s}$, we just need to assume that $x\in \Omega_1$, defined as (\ref{Omega_1}). We will do the sup-norm first. Because of  triangle inequality and Lemma \ref{lem1},
\begin{equation}\label{smallness00}
\begin{aligned}
\left|\sum_{i=1}^k\sigma_i^{\frac{n+2}{n-2}}-(\sum_{i=1}^k\sigma_i)^{\frac{n+2}{n-2}}\right|&\leq \sum_{i=2}^k\sigma_i^{\frac{n+2}{n-2}}+
\left|\sigma_1^{\frac{n+2}{n-2}}-(\sigma_1+\sum_{i=2}^k\sigma_i)^{\frac{n+2}{n-2}}\right|\\
&\leq C\left(\sum_{i=2}^k\sigma_i^{\frac{n+2}{n-2}}+\sigma_1^{\frac{4}{n-2}}\sum_{i=2}^k\sigma_i+(\sum_{i=2}^k\sigma_i)^{\frac{n+2}{n-2}}\right).
\end{aligned}    
\end{equation}

The estimates of the three terms on the right hand side of (\ref{smallness00}) are similar. Here we just give the proof of the estimate for the third term. Using (\ref{Omega_1 geometry}) and (\ref{tau eqn1}) with $\tau=n-2-(s+2)\frac{n-2}{n+2}>1$, we have
\begin{equation}\label{smallnesssup1}
\begin{aligned}
&\sup_{x\in\Omega_1}\left(\sum_{i=1}^{k}\<x-P^i\>^{-s-2}\right)^{-1}\left|(\sum_{i=2}^k\sigma_i)^{\frac{n+2}{n-2}}\right|\\
\leq& C\sup_{x\in\Omega_1}d(x)^{s+2}\left(\sum_{i=2}^k\<x-P^i\>^{-(n-2)}\right)^{\frac{n+2}{n-2}}\leq C(k/r)^{n-s}.
\end{aligned}
\end{equation}
Combining the three terms together, we can get
\begin{equation}\label{smallnesssup2}
\norm{X^{0}_{k,s+2}}{\sum_{i=1}^k\sigma_i^{\frac{n+2}{n-2}}-(\sum_{i=1}^k\sigma_i)^{\frac{n+2}{n-2}}}\leq C(k/r)^{n-s}.
\end{equation}

For the H{\"o}lder semi-norm, fix $x\in \Omega_1$, $B(x)=B(x,d(x)/2)$.  We use triangle inequality from Lemma \ref{holder semi norm lemma} and mean value theorem to obtain 
\begin{align*}
&\left[\sum_{i=1}^k\sigma_i^{\frac{n+2}{n-2}}-(\sum_{i=1}^k\sigma_i)^{\frac{n+2}{n-2}}\right]_{\alpha,B(x)}\\
\leq& \left[\sigma_1^{\frac{n+2}{n-2}}-(\sigma_1+\sum_{i=2}^k\sigma_i)^{\frac{n+2}{n-2}}\right]_{\alpha,B(x)}+\left[\sum_{i=2}^k\sigma_i^{\frac{n+2}{n-2}}\right]_{\alpha,B(x)}\\
= &\frac{n+2}{n-2}\left[\int_{0}^{1}|\sigma_1+\theta\sum_{i=2}^k\sigma_i|^{\frac{4}{n-2}}d\theta\sum_{i=2}^k\sigma_i\right]_{\alpha,B(x)} +\left[\sum_{i=2}^k\sigma_i^{\frac{n+2}{n-2}}\right]_{\alpha,B(x)}.
\end{align*}
We keep the second term and use the product rule from Lemma \ref{holder semi norm lemma} and the interpolation inequality 2 from Lemma \ref{holder semi norm lemma} with $\delta=1$ for the first term. So,
$$
\left[\sum_{i=1}^k\sigma_i^{\frac{n+2}{n-2}}-(\sum_{i=1}^k\sigma_i)^{\frac{n+2}{n-2}}\right]_{\alpha,B(x)}\leq C(\RN{1}+\RN{2}+\RN{3})+\RN{4},
$$
where $\RN{4}:=\left[\sum_{i=2}^k\sigma_i^{\frac{n+2}{n-2}}\right]_{\alpha,B(x)}$, and 
\begin{align*}
\RN{1}&:=\left|\int_{0}^{1}|\sigma_1+\theta\sum_{i=2}^k\sigma_i|^{\frac{4}{n-2}}\right|_{B(x)}\left[\sum_{i=2}^k\sigma_i\right]_{\alpha,B(x)},\\
\RN{2}&:=d(x)^{-\alpha}\left|\int_{0}^{1}|\sigma_1+\theta\sum_{i=2}^k\sigma_i|^{\frac{4}{n-2}}\right|_{B(x)}\left|\sum_{i=2}^k\sigma_i\right|_{B(x)},\\
\RN{3}&:=
d(x)^{\frac{4}{n-2}-\alpha}\left[\int_{0}^{1}|\sigma_1+\theta\sum_{i=2}^k\sigma_i|^{\frac{4}{n-2}}\right]_{\frac{4}{n-2},B(x)}\left|\sum_{i=2}^k\sigma_i\right|_{B(x)}.
\end{align*}

Again, the estimates of \RN{1}, \RN{2}, \RN{4} are similar. Here we just give the estimate of \RN{4}. Use the interpolation inequality 1 from Lemma \ref{holder semi norm lemma} with $\delta=1$:
\begin{equation}\label{interpolation}
\begin{aligned}
\RN{4}&\leq C\left(d(x)^{-\alpha}\left|\sum_{i=2}^k\sigma_i^{\frac{n+2}{n-2}}\right|_{B(x)}+d(x)^{1-\alpha}\left|D(\sum_{i=2}^k\sigma_i^{\frac{n+2}{n-2}})\right|_{B(x)}\right)\\
&\leq Cd(x)^{-\alpha}\sum_{i=2}^k\<x-P^i\>^{-(n+2)}.
\end{aligned}
\end{equation}
Then we can use (\ref{interpolation}), (\ref{Omega_1 geometry}) and (\ref{tau eqn1}) with $\tau=n-s>1$ to obtain
\begin{equation}\label{smallnessHolder1}
\begin{aligned}
\sup_{x\in\Omega_1}\left(\sum_{i=1}^{k}\<x-P^i\>^{-s-2-\alpha}\right)^{-1}\RN{4}&\leq C\sup_{x\in\Omega_1}d(x)^{s-2}\sum_{i=2}^k\<x-P^i\>^{-(n-2)}\\
&\leq C (k/r)^{n-s}.
\end{aligned}
\end{equation}

For \RN{3}, we use (\ref{alpha=beta}) from Lemma \ref{lem1-holder}:
\begin{align*}
\RN{3}\leq &Cd(x)^{\frac{4}{n-2}-\alpha}\left(|D\sigma_1|^{\frac{4}{n-2}}_{B(x)}+\left|D\sum_{i=2}^k\sigma_i\right|^{\frac{4}{n-2}}_{B(x)}\right)\left|\sum_{i=2}^k\sigma_i\right|_{B(x)}\\
\leq &Cd(x)^{-4-\alpha}\left(\sum_{i=2}^{k}\<x-P^i\>^{-(n-2)}\right)+Cd(x)^{-\alpha}\left(\sum_{i=2}^k\<x-P^i\>^{-(n-2)}\right)^{\frac{n+2}{n-2}}.
\end{align*}
Similar with (\ref{smallnesssup1}) and (\ref{smallnessHolder1}), we have
\begin{equation*}
\sup_{x\in\Omega_1}\left(\sum_{i=1}^{k}\<x-P^i\>^{-s-2-\alpha}\right)^{-1}\RN{3} \leq C (k/r)^{n-s}.
\end{equation*}
Combining the terms $\RN{1}$, $\RN{2}$, $\RN{3}$, and \RN{4} together, we can obtain
\begin{equation}\label{smallnessHolder2}
\left[\sum_{i=1}^k\sigma_i^{\frac{n+2}{n-2}}-(\sum_{i=1}^k\sigma_i)^{\frac{n+2}{n-2}}\right]_{X^{0,\alpha}_{k,s+2}}\leq C (k/r)^{n-s}.
\end{equation}
Combining estimates (\ref{smallnesslinear}), (\ref{smallnesssup2}), (\ref{smallnessHolder2}) together, we can obtain the estimate (\ref{smallness1}).
\end{proof}

\begin{lem}\label{smallnessZ}
Let $n\geq 6$, $\frac{n-2}{2}<s< n-2$, $0<\alpha<\frac{4}{n-2}$, $c_0>0$, and $\tau_0\in \r$. Then there exists $C=C(n,s,\alpha,c_0,\tau_0)$, such that for any $0<\epsilon<1$, $k\geq 3$, $(\Xi,\Lambda)\in\mathcal{C}_k$, $j=1,\cdots,k$, $\mu=0,\cdots,n$,
\begin{equation*}
\norm{X^*}{L_{\Xi,\Lambda} Z_{j,\mu}}\leq C\left(\epsilon+(k/r)^{2}\right).
\end{equation*}
\end{lem}
\begin{proof}
Without loss of generality, pick $j=1$, $\mu=0$. First, by the definition of $L_{\Xi,\Lambda}$ (\ref{def of L}) and the equation $-\Delta Z_{1,0}=n(n+2)\sigma_1^{\frac{4}{n-2}}Z_{1,0}$, we have
$$
L_{\Xi,\Lambda}Z_{1,0}=-(L_{g }-\Delta)Z_{1,0}+n(n+2)\left(\sigma_1^{\frac{4}{n-2}}-(\sum_{i=1}^k\sigma_i)^{\frac{4}{n-2}}\right)Z_{1,0}.
$$

We use the same method of Proposition \ref{Welldefined} to deal with the linear part:
$$
|(L_g-\Delta) Z_{1,0}|(x)\leq C\epsilon \<x-P^1\>^{-n}, \text{ for any } x\in \r^n.
$$
Then for any test function $\phi\in  X $, we have
\begin{equation}\label{use symmetry}
\begin{aligned}
\left|\int_{\r^n}\phi |(L_g-\Delta )Z_{1,0}|\right|&\leq C\epsilon\norm{X^0_{k,s}}{\phi}\int_{\r^n}\<x-P^1\>^{-n}\sum_i\<x-P^i\>^{-s}\\
&\leq C\epsilon\norm{X^0_{k,s}}{\phi}\frac{1}{k}\int_{\r^n}\sum_{j}\<x-P^j\>^{-n}\sum_{i}\<x-P^i\>^{-s}\\
&\leq C\epsilon\norm{X^0_{k,s}}{\phi}\frac{1}{k}\int_{\r^n}\sum_{i}\<x-P^i\>^{-s-n+\tau}\leq C\epsilon\norm{X^0_{k,s}}{\phi},
\end{aligned}
\end{equation}
where we have used (\ref{tau eqn2}) with $\tau<s$ in the third inequality.

Next, we focus on the nonlinear term. Because of the Lemma \ref{lem1},
$$
\left|\sigma_1^{\frac{4}{n-2}}-(\sum_{i=1}^k\sigma_i)^{\frac{4}{n-2}}\right|\leq C \sigma_1^{\frac{4}{n-2}-1}\sum_{i=2}^k\sigma_i.
$$
Then for any test function $\phi\in  X $, we have
\begin{align*}
&\int_{\r^n}\left|(\sigma_1^{\frac{4}{n-2}}-(\sum_{i=1}^k\sigma_i)^{\frac{4}{n-2}})Z_{1,0}\phi\right|\\
\leq& C\norm{X^0_{k,s}}{\phi} \int_{\r^n}\<x-P^1\>^{-4}\sum_{i=2}^k\<x-P^i\>^{-(n-2)}\sum_{j=1}^k\<x-P^j\>^{-s}.
\end{align*}
Because of the periodicity, we know that
\begin{align*}
&\int_{\r^n}\<x-P^1\>^{-4}\sum_{i=2}^k\<x-P^i\>^{-(n-2)}\sum_{j=1}^k\<x-P^j\>^{-s} \\
=&\frac{1}{k}\int_{\r^n}\left(\sum_{l}\<x-P^l\>^{-4}\sum_{i\neq l}\<x-P^i\>^{-(n-2)}\right)\sum_j\<x-P^j\>^{-s}.
\end{align*}
Use (\ref{tau eqn2}) with $\tau=2$, 
\begin{align*}
\sum_{l}\<x-P^l\>^{-4}\sum_{i\neq l}\<x-P^i\>^{-(n-2)}&\leq C\sum_{l=1}^k\<x-P^l\>^{-n}\sum_{i\neq l}|P^i-P^l|^{-2}\\
&\leq C(k/r)^{2}\sum_{i=1}^k\<x-P^i\>^{-n}.
\end{align*}
Using (\ref{tau eqn2}) with $\tau<s$ again, we have 
\begin{align*}
&\int_{\r^n}\<x-P^1\>^{-4}\sum_{i=2}^k\<x-P^i\>^{-(n-2)}\sum_{j=1}^k\<x-P^j\>^{-s}\\
\leq& C(k/r)^{2}\frac{1}{k}\int_{\r^n}\sum_{i=1}^k\<x-P^i\>^{-n-s+\tau}\leq C(k/r)^{2}.
\end{align*}
So, for any test function $\phi\in X$, there is
\begin{equation}\label{use symmetry nonlinear}
\int_{\r^n}\left|(\sigma_1^{\frac{4}{n-2}}-(\sum_{i=1}^k\sigma_i)^{\frac{4}{n-2}})Z_{1,0}\phi\right|\leq C(k/r)^{2}\norm{X^0_{k,s}}{\phi}
\end{equation}
Combining (\ref{use symmetry}) and (\ref{use symmetry nonlinear}), we can obtain Lemma \ref{smallnessZ}.
\end{proof}

Now, we are ready to prove Proposition \ref{Surjectivity}.
\begin{proof}
First, we prove that there exist $\epsilon_0, k_0, C$ such that for any $0<\epsilon\leq \epsilon_0$, any $k\geq k_0$, $(\Xi,\Lambda)\in \mathcal{C}_k$, we have the following estimate
\begin{equation}\label{a prior estimate}
\norm{ X }{\phi}\leq C\norm{ Y }{f}, \text{ for any }f\in Y_1+Y_2,
\end{equation}
where $\phi\in X_2$, $f_1\in Y_1$ such that $f=L_{\Xi,\Lambda}\phi+f_1$.

We prove it by contradiction argument. If there are not such $\epsilon_0$, $k_0$ and $C$ for fixed $n,s,\alpha, c_0, \tau_0$, consider a sequence $\epsilon_l\rightarrow 0$ and a sequence $k_l\rightarrow \infty$, $(\Xi_l,\Lambda_l)\in \mathcal{C}_{k_l}$ and $f_l\in (Y_1)_l+(Y_2)_l$ satisfying
\begin{equation}\label{linearized eqn in contradiction arg}
f_l=L_{\Xi_l,\Lambda_l}\phi_l+(f_1)_l,
\end{equation}
where $\phi_l \in (X_2)_l$, $(f_1)_l\in (Y_1)_l$ with $\norm{X_l}{\phi_l}=1$ but
$$
\norm{Y_l}{f_l}\rightarrow 0, \text{ when }l\rightarrow\infty.
$$
For simplicity, we drop the subscript $l$ in the following proof. Applying the $\Gamma$ operator on both side of the equation (\ref{linearized eqn in contradiction arg}), we have
\begin{align*}
\phi(x)=&C(n)\int_{\r^n}\frac{(L_{g }-\Delta)\phi(y)}{|x-y|^{n-2}}+\frac{n(n+2)u_{\Xi,\Lambda}(y)^{\frac{4}{n-2}}\phi(y)}{|x-y|^{n-2}}dy
\\&+C(n)\int_{\r^n}\frac{f(y)}{|x-y|^{n-2}}
-\frac{f_1(y)}{|x-y|^{n-2}}dy.
\end{align*}

Let us deal with the fourth term. We write $f_1=\sum_{j,\mu}c_{j,\mu}\Delta Z_{j,\mu}$ because $f_1\in Y_1$. Because of estimate (\ref{bound of Z}), we know that $\norm{ Y }{\sum_{j,\mu}\Delta Z_{j,\mu}}\leq C$. Then, 
\begin{equation}\label{f_1 bound}
\norm{ Y }{f_1}\leq C\max_{j,\mu}|c_{j,\mu}|.
\end{equation}
Without loss of generality, assume that $|c_{1,0}|=\max_{j,\mu}|c_{j,\mu}|$. Multiply the linearized equation (\ref{linearized eqn in contradiction arg}) by $Z_{1,0}$ and integrate over $\r^n$,
$$
\int_{\r^n}L_{\Xi,\Lambda} \phi Z_{1,0}=\int_{\r^n}fZ_{1,0}-\sum_{j,\mu}c_{j,\mu}\int_{\r^n}Z_{1,0}\Delta Z_{j,\mu}.
$$
Then do integration by parts on both sides, use (\ref{orthognality of Z}),
$$
\left|\int_{\r^n}\phi L_{\Xi,\Lambda}Z_{1,0}-\int_{\r^n}fZ_{1,0}+a_1c_{1,0}\lambda_1^{-2}\right|\leq Ck(k/r)^{n-2}|c_{1,0}|.
$$
Because $a_1, \lambda_1$ are bounded and $k(k/r)^{n-2}\rightarrow 0$ as $k\rightarrow \infty$, we have
\begin{equation}\label{a}
\left|\int_{\r^n}\phi L_{\Xi,\Lambda}Z_{1,0}-\int_{\r^n}fZ_{1,0}\right|\geq C^{-1}|c_{1,0}|.
\end{equation}

Next, we try to bound $\left|\int_{\r^n}\phi L_{\Xi,\Lambda}Z_{1,0}-\int_{\r^n}fZ_{1,0}\right|$ from above. For the second term, using the derivation of (\ref{use symmetry}), we know that
\begin{equation}\label{b}
\left|\int_{\r^n}fZ_{1,0}\right|\leq C\norm{X^{0}_{k,s+2}}{f} \int_{\r^n}\<x-P^1\>^{-(n-2)}\sum_{i}\<x-P^i\>^{-s-2}\leq C\norm{X^{0}_{k,s+2}}{f}.
\end{equation}

Then combining estimates (\ref{f_1 bound}), (\ref{a}), (\ref{b}) and Lemma \ref{smallnessZ} together, we can obtain a control of $\norm{ Y }{f_1}$,
\begin{equation}\label{c_ij}
\norm{ Y }{f_1}\leq C\max_{j,\mu}|c_{j,\mu}|=C|c_{1,0}|\leq C(\epsilon+(k/r)^{2})\norm{ X }{\phi}+C\norm{X^0_{k,s+2}}{f}=o(1).
\end{equation}

Because of Corollary \ref{Gammaboundedness}, the boundedness of $\Gamma$, and estimate (\ref{c_ij}), the forth term
$$
\left|\int_{\r^n}\frac{f_1}{|x-y|^{n-2}}\right|=o(1)\sum_{i=1}^k\<x-P^i\>^{-s}.
$$
Similarly, because of the boundedness of $\Gamma$ and  $\norm{ Y }{f}=o(1)$, we have that 
$$
\left|\int_{\r^n}\frac{f}{|x-y|^{n-2}}\right|=o(1)\sum_{i=1}^k\<x-P^i\>^{-s}.
$$
Since the estimate of $\norm{ Y }{(L_{g }-\Delta)\phi}$ in Proposition \ref{Welldefined}, 
$$
\norm{ Y }{(L_{g }-\Delta)\phi}\leq C\epsilon=o(1).
$$
Then because of the boundedness of $\Gamma$, we have
$$
\left|\int_{\r^n}\frac{(L_{g }-\Delta)\phi}{|x-y|^{n-2}}\right|=o(1)\sum_{i=1}^k\<x-P^i\>^{-s}.
$$

The second term is the tricky part: since $\frac{4}{n-2}\leq1$,
$$
\left|\int_{\r^n}\frac{n(n+2)u_{\Xi,\Lambda}^{\frac{4}{n-2}}\phi}{|x-y|^{n-2}}\right|\leq C\int_{\r^n}\frac{1}{|x-y|^{n-2}}\sum_{i=1}^k\<y-P^i\>^{-s}\sum_{j=1}^k\<y-P^j\>^{-4}
$$
Then we use Lemma \ref{lem3} with
\begin{equation*}
\tau=\left\{
\begin{array}{lr}
2-\delta_0, & \text{if }\frac{n-2}{2}<s<n-2-\delta_0\\
2-\frac{s+(n-2)}{2}, & \text{if }n-2-\delta_0\leq s<n-2
\end{array}
\right.
\end{equation*} 
for some $\delta_0$ small, to obtain 
$$
\left|\int_{\r^n}\frac{n(n+2)u_{\Xi,\Lambda}^{\frac{4}{n-2}}\phi}{|x-y|^{n-2}}\right|\leq C \int_{\r^n}\frac{1}{|x-y|^{n-2}}\sum_{i=1}^k\<y-P^i\>^{-s-4+\tau}.
$$
Then because of Lemma \ref{lem4}, we have that 
$$
\left|\int_{\r^n}\frac{n(n+2)u_{\Xi,\Lambda}^{\frac{4}{n-2}}\phi}{|x-y|^{n-2}}\right|\leq C\sum_{i=1}^k\<x-P^i\>^{-s-2+\tau}.
$$ 
Adding the four terms together, we have the following pointwise estimate
\begin{equation}\label{pointwise estimate}
\left(\sum_{i=1}^k\<x-P^i\>^{-s}\right)^{-1}|\phi(x)|\leq o(1)+Cd(x)^{\tau-2}, \quad x\in \r^n.
\end{equation}
The key point here is that $\tau<2$.

In the other direction, because of the Schauder estimate (\ref{schauder of L}) and the assumption $\norm{ X }{\phi}=1$, we have
$$
1=\norm{ X }{\phi}\leq C(\norm{ Y }{f-f_1}+\norm{X^0_{k,s}}{\phi})=o(1)+C\norm{X^0_{k,s}}{\phi}.
$$
So, we obtain $\norm{X^0_{k,s}}{\phi}\geq \frac{1}{C_1}$, i.e. $\sup_{x\in\r^n}\left(\sum_{j=1}^k\<x-P^j\>^{-s}\right)^{-1}|\phi(x)|\geq \frac{1}{C_1}$.

However, because of the pointwise estimate (\ref{pointwise estimate}), we can choose $\bar{R}$ large enough but NOT depending on the sequential index $l$ such that
$$
\sup_{\r^n\setminus\cup_{i=1}^kB(P^i,\bar{R})}\left(\sum_{j=1}^k\<x-P^j\>^{-s}\right)^{-1}|\phi(x)|\leq o(1)+\bar{R}^{\tau-2}< \frac{1}{C_1}.
$$
In other words, the outside part can not reach $ \frac{1}{C_1}$. Then the sup can only be taken in the inside part, i.e., there exist $1\leq \bar{i}_l\leq k_l$,  such that 
$$
\sup_{B(P^{\bar{i}},\bar{R})}\left(\sum_{j=1}^k\<x-P^j\>^{-s}\right)^{-1}|\phi(x)|\geq  \frac{1}{C_1}.
$$
In other words,
$$
\sup_{B(P^{\bar{i}},\bar{R})}|\phi(x)|\geq  \frac{1}{C_1}\inf_{B(P^{\bar{i}},\bar{R})}\sum_{j=1}^k\<x-P^j\>^{-s}\geq \frac{(1+\bar{R})^{-s}}{C_2}>0.
$$

However, because $r/k \rightarrow\infty$, after passing to subsequence, $\bar{\phi}_l(x):=\phi_l(x+P^{\bar{i}_l})$ converges uniformly in any compact set to a solution $\Phi$ of 
$$
-\Delta \Phi-n(n+2)\sigma_{\bar{\xi},\bar{\lambda}}^{\frac{4}{n-2}}\Phi=0, \text{ in }\r^n.
$$
for $\bar{\lambda}=\lim_l\lambda_{\bar{i}_l}$ and $\bar{\xi}=\lim_l\xi^{\bar{i}_l}$ with $\bar{\lambda}\in \left[\frac{3}{4},\frac{4}{3}\right]$ and $|\bar{\xi}|\in \left[0,\frac{1}{2}\right]$. Also since $\phi_l\in N_{k_l,\Xi_l,\Lambda_l}$,
$$
\<\bar{\phi}_l(x),Z_{\bar{i}_l,\mu}(x+P^{\bar{i}_l})\>_{\mathcal{D}^{1,2} }=0, \text{ for } \mu=0,\cdots,n.
$$
Then $\Phi$ is orthogonal to all the non-trivial kernel of the following PDE on $\r^n$: $-\Delta \Phi-n(n+2)\sigma_{\bar{\xi},\bar{\lambda}}^{\frac{4}{n-2}}\Phi=0$. So $\Phi\equiv0$. But it contradicts with 
$$
\sup_{B(0,\bar{R})}|\Phi(x)|\geq \frac{(1+\bar{R})^s}{C_2}>0.
$$
So, we have proved the estimate (\ref{a prior estimate}).

Second, the closeness of $Y_2$ follows from (\ref{a prior estimate}) and (\ref{boundedness of L}). $Y_1$ is closed because $\dim Y_1=k(n+1)$. Next, we want to prove the decomposition $Y_1\cap Y_2=\{0\}$. For $f\in Y_1\cap Y_2$, there are $\phi\in X_2$ and $f_1\in Y_1$ such that $f=L_{\Xi,\Lambda}\phi=f_1$. In other words,
$$
0=L_{\Xi,\Lambda}\phi-f_1.
$$
Using (\ref{a prior estimate}), we know that $\phi=0$ and $f=0$. In other words, we know that $Y_1\cap Y_2=\{0\}$. 

Moreover, we can count the dimension to prove that $ Y =Y_1\oplus Y_2$. Because the inclusion map $\iota:X_2\rightarrow  X $ is Fredholm of index $-k(n+1)$ and $L_{\Xi,\Lambda}$ is Fredholm of index zero by Corollary \ref{fredholm}. Then we have that $L_{\Xi,\Lambda}\circ \iota:X_2\rightarrow  X \rightarrow  Y $ is Fredholm of index $-k(n+1)$, i.e.
$$
\dim \ker(L_{\Xi,\Lambda}\circ \iota)-\codim Y_2= -k(n+1).
$$
$L_{\Xi,\Lambda}\circ \iota$ has trivial kernel because we can use (\ref{a prior estimate}) with $f=0$. So we know that $\codim Y_2=k(n+1)$. Combining with $\dim Y_1=k(n+1)$ and $Y_1\cap Y_2=\{0\}$, we can obtain $ Y =Y_1\oplus Y_2$.

Third, using estimates (\ref{c_ij}) and (\ref{a prior estimate}) for $f=f_1+f_2$, where $f_j\in Y_j$, for $j=1,2$, we obtain
$$
\norm{ Y }{f_1} \leq C\norm{ X }{\phi}+C\norm{X^0_{k,s+2}}{f}\leq C\norm{ Y }{f}.
$$
Then we can obtain estimate (\ref{bound of proj}) because of triangle inequality $\norm{ Y }{f_2}\leq \norm{ Y }{f_1}+\norm{ Y }{f}$.

In the end, using (\ref{a prior estimate}) for $f\in Y_2$, we can prove estimate (\ref{invertibility}). 

\end{proof}

\section{Nonlinear Elliptic Operator $F$} 

\subsection{Implicit Function Theorem}

Lyapunov-Schmidt reduction method is to divide the process of solving the equation $F(u_{\Xi,\Lambda}+\phi)=0$ into solving the following two equations:
\begin{equation*}
\begin{aligned}
\Proj_{Y_2}\circ F(u_{\Xi,\Lambda}+\phi)&=0;\\
(\Id-\Proj_{Y_2})\circ F(u_{\Xi,\Lambda}+\phi)&=0.
\end{aligned}
\end{equation*}

We will use the following easily-verified version of implicit function theorem, see \cite{li1998dirichlet}, to find a $\phi_{\Xi,\Lambda}\in X_2$ for each $(\Xi,\Lambda)\in \mathcal{C}_k$ to solve the first equation. 

\begin{lem}\label{nirenbergLemma}
Let $\mathcal{B}_1$ and $\mathcal{B}_2$ be two Banach spaces, let $z_0\in \mathcal{B}_1$, $a>0$ and $B(z_0,a)=\{z\in \mathcal{B}_1\mid \norm{\mathcal{B}_1}{z-z_0}\leq a\}$. Suppose $G$ is a $C^1$ map from $B(z_0,a)$ to $\mathcal{B}_2$ with $G'(z_0)$ invertible, and satisfying, for some $0<\theta<1$,
$$
\norm{\mathcal{B}_1}{G'(z_0)^{-1}G(z_0)}\leq (1-\theta)a,
$$
and
$$
\norm{\mathcal{B}_1}{G'(z_0)^{-1}(G(z)-G(z'))-(z-z')}\leq \theta\norm{\mathcal{B}_1}{z-z'}, \text{ for any } z,z'\in B(z_0,a).
$$
Then we have an unique $z\in B(z_0,a)$ solving $G(z)=0$.
\end{lem}

\begin{prop}\label{solvability}
Let $n\geq 6$, $0<\alpha<\frac{4}{n-2}$, $\frac{n-2}{2}<s<n-2$, $c_0>0$ and $\tau_0\in\r$. There exist $\epsilon_0, k_0, C$ depending on $ n,s,\alpha, c_0, \tau_0$ such that for any $0<\epsilon\leq \epsilon_0$, any $k\geq k_0$, $(\Xi,\Lambda)\in \mathcal{C}_k$, there exists a solution $\phi_{\Xi,\Lambda}\in X_2$ solving
\begin{equation}\label{prob2}
\Proj_{Y_2}\circ F(u_{\Xi,\Lambda}+\phi_{\Xi,\Lambda})=0.
\end{equation}
Moreover, we have
\begin{equation}\label{phi & F(u) same order}
\norm{ X }{\phi_{\Xi,\Lambda}}\leq C\norm{ Y }{F(u_{\Xi,\Lambda})}.
\end{equation}
\end{prop}
\textbf{Remark.} We need the condition $\frac{n-2}{2}<s<n-2$ because we will use Lemma \ref{smallnessfu} in the proof of this proposition.
\begin{proof}
Let $\mathcal{B}_1=X_2$ with the norm $\norm{\mathcal{B}_1}{\cdot}=\norm{ X }{\cdot}$ and $\mathcal{B}_2=Y_2$. Consider a $C^1$ map $G: \mathcal{B}_1\rightarrow \mathcal{B}_2$
$$
G(\phi)=\Proj_{Y_2}\circ F(u_{\Xi,\Lambda}+\phi).
$$ 
Then $G(0)=\Proj_{Y_2}\circ F(u_{\Xi,\Lambda})$. Picking $z_0=0$, we know that $G'(0)=\Proj_{Y_2}\circ L_{\Xi,\Lambda}\circ \iota: \mathcal{B}_1\rightarrow \mathcal{B}_2$ is invertible because of Proposition \ref{Surjectivity}.

Using estimates (\ref{bound of proj}) and (\ref{invertibility}), we have 
\begin{equation*}
\norm{\mathcal{B}_1}{G'(0)^{-1}G(0)}\leq C\norm{Y}{G(0)}\leq C_1\norm{ Y }{F(u_{\Xi,\Lambda})}.
\end{equation*}
Pick $a=2C_1\norm{ Y }{F(u_{\Xi,\Lambda})}$. Then it is trivial that $\norm{\mathcal{B}_1}{G'(0)^{-1}G(0)}\leq \frac{1}{2}a$. 

Let $\phi_1,\phi_2\in B(0,a)$. Using estimates (\ref{bound of proj}) and (\ref{invertibility}) again, combining (\ref{Lip of F}), we have that 
\begin{align*}
&\norm{\mathcal{B}_1}{G'(0)^{-1}(G(\phi_1)-G(\phi_2))-(\phi_1-\phi_2)}\\
\leq &C\norm{ Y }{F(u_{\Xi,\Lambda}+\phi_1)-F(u_{\Xi,\Lambda}+\phi_2)-L_{\Xi,\Lambda}(\phi_1-\phi_2)}\\
\leq &C_2a^{\frac{4}{n-2}-\alpha}\norm{\mathcal{B}_1}{\phi_1-\phi_2}.
\end{align*}

Because of estimate (\ref{smallness1}), we know that $a\leq C_3(\epsilon+(k/r)^{n-s})$. So we can pick $\epsilon_0$ small, $k_0$ large such $C_2a^{\frac{4}{n-2}-\alpha}<\frac{1}{2}$. Then, for any $(\Xi,\Lambda)\in \mathcal{C}_k$, we can use Lemma \ref{nirenbergLemma} to find $\phi_{\Xi,\Lambda}\in X_2$ solving $G(\phi_{\Xi,\Lambda})=0$ with estimate
$$
\norm{ X }{\phi_{\Xi,\Lambda}}< a=2C_1\norm{ Y }{F(u_{\Xi,\Lambda})}.
$$

\end{proof}

Then solving the second equation can be considered as finding a zero point of a finite dimensional function. Since $F$ is variational, we will reduce the solvability of second equation to the existence of a critical point of a finite dimensional function $I:\mathcal{C}_k\rightarrow \r$:
$$
I(\Xi,\Lambda)=I_g(u_{\Xi,\Lambda}+\phi_{\Xi,\Lambda}),
$$
where $I_g$ is defined in (\ref{def of energy I}) and $\phi_{\Xi,\Lambda}$ is from the Proposition \ref{solvability}.

\begin{coro}\label{deduce to critical pt}
Under the same condition, $I(\Xi,\Lambda)$ is a $C^1$ function. Also, if $(\bar{\Xi},\bar{\Lambda})\in \mathcal{C}_k$ is a critical point of $I$, then 
$$
(\Id-\Proj_{Y_2})\circ F(u_{\bar{\Xi},\bar{\Lambda}}+\phi_{\bar{\Xi},\bar{\Lambda}})=0.
$$
\end{coro}
\begin{proof}
Fixing $(\bar{\Xi},\bar{\Lambda})$ as a critical point of $I$, we will simply denote $\phi=\phi_{\bar{\Xi},\bar{\Lambda}}$, $u=u_{\bar{\Xi},\bar{\Lambda}}$. Because of equation (\ref{prob2}), we can write
$$
F(u+\phi)=\sum_{j,\mu}c_{j,\mu}\Delta Z_{j,\mu},
$$
where $Z_{j,\mu}$ are from (\ref{def of tangent vectors}). Without loss of generality, we assume that $|c_{1,0}|=\max_{j,\mu}|c_{j,\mu}|$. Then, we have
\begin{align*}
0=\partial_{\lambda_1}I(\bar{\Xi},\bar{\Lambda})
&=-\int_{\r^n}F(u+\phi)(Z_{1,0}+\partial_{\lambda_1}\phi)\\
&=\sum_{j,\mu}c_{j,\mu}\left(\<Z_{j,\mu},Z_{1,0}\>_{\mathcal{D}^{1,2} }+\<Z_{j,\mu},\partial_{\lambda_1}\phi\>_{\mathcal{D}^{1,2} }\right)\\
&=\sum_{j,\mu}c_{j,\mu}\left(\<Z_{j,\mu},Z_{1,0}\>_{\mathcal{D}^{1,2} }-\<\partial_{\lambda_1}Z_{j,\mu},\phi\>_{\mathcal{D}^{1,2} }\right),
\end{align*}	
where we have used $\<Z_{j,\mu},\phi\>_{\mathcal{D}^{1,2} }=0$ in the last inequality. Using the orthogonality estimates (\ref{orthognality of Z}), we obtain
\begin{equation}\label{c_10 estimate}
a_1\lambda_1^{-2}|c_{1,0}|\leq Ck(k/r)^{n-2}|c_{1,0}|+\sum_{j,\mu}\left|\<\partial_{\lambda_1}Z_{j,\mu},\phi\>_{\mathcal{D}^{1,2} }\right||c_{1,0}|.
\end{equation}

Next, we will estimate the second term on the right hand side of inequality above. Since $|\partial_{\lambda_1}Z_{j,\mu}|\leq C\delta_{1j}\<x-P^1\>^{-(n-2)}$ and $|\Delta\phi|\leq \norm{X^{2}_{k,s}}{\phi}\sum_{i}\<x-P^i\>^{-s-2}$, we have
$$
\sum_{j,\mu}\left|\<\partial_{\lambda_1}Z_{j,\mu},\phi\>_{\mathcal{D}^{1,2} }\right|\leq C\norm{X^{2}_{k,s}}{\phi}\int_{\r^n}\<x-P^1\>^{-(n-2)}\sum_{i}\<x-P^i\>^{-s-2}.
$$
Then we can use the derivation of (\ref{use symmetry}), combining with estimates (\ref{smallness1}) and (\ref{phi & F(u) same order}), to obtain  
$$
\sum_{j,\mu}\left|\<\partial_{\lambda_1}Z_{j,\mu},\phi\>_{\mathcal{D}^{1,2} }\right|\leq  C\norm{X^{2}_{k,s}}{\phi}\leq C(\epsilon+(k/r)^{n-s}).
$$

Plugging the above estimate in (\ref{c_10 estimate}), we have
$$
|c_{1,0}|\leq C(k(k/r)^{n-2}+\epsilon+(k/r)^{n-s})|c_{1,0}|.
$$
Picking $\epsilon_0$ small and $k_0$ large such that $|c_{1,0}|<\frac{1}{2}|c_{1,0}|$, then we know that $\max_{j,\mu}|c_{j,\mu}|=|c_{1,0}|=0$.
\end{proof}

\subsection{Estimates of Leading Terms}

In the rest of this paper, we will fix $\epsilon=\epsilon_0$ where the $\epsilon_0=\epsilon_0(n,c_0)$ is from Corollary \ref{deduce to critical pt}.

Recall the definition of $I_g$ from (\ref{def of energy I}). For any $v,\phi,\psi \in \mathcal{D}^{1,2} $, write that 
\begin{align*}
I'_g(v)\phi&=\int_{\r^n}g^{\mu\nu}D_{\mu}vD_{\nu}\phi+c(n)R_{g }v\phi-n(n-2)|v|^{\frac{4}{n-2}}v\phi;\\
I''_g(v)(\phi,\psi)&=\int_{\r^n}g^{\mu\nu}D_{\mu}\phi D_{\nu}\psi+c(n)R_{g }\phi\psi-n(n+2)|v|^{\frac{4}{n-2}}\phi\psi.
\end{align*}
Specially, denote
$$
I_{\delta}(v)=\int_{\r^n}\frac{1}{2}|Dv|^2-\frac{(n-2)^2}{2}|v|^{\frac{2n}{n-2}}.
$$
According to the previous section, we need to find out the critical point of       
$$
I(\Xi,\Lambda)=I_{g }(u_{\Xi,\Lambda}+\phi_{\Xi,\Lambda})
$$ 
where $\phi_{\Xi,\Lambda}$ is found in Proposition \ref{solvability}.

To give the leading term of $I(\Xi,\Lambda)$, we need few more definitions. Recall the notations: $\sigma_{\xi,\lambda}(x)=\lambda^{\frac{n-2}{2}}(1+\lambda^2|x-\xi|^2)^{-\frac{n-2}{2}}$. Consider
\begin{equation}\label{def of G}
\hat{G}(\xi,\lambda):=\int_{\r^n}\frac{1}{4}\sum_{\mu\nu}(\hat{H}^2)_{\mu\nu}D_{\mu}\sigma_{\xi,\lambda} D_{\nu}\sigma_{\xi,\lambda}-\frac{c(n)}{8}\sum_{\mu\nu\alpha}(D_{\alpha}\hat{H}_{\mu\nu})^2\sigma_{\xi,\lambda}^2
\end{equation}
where $\hat{H}$ is given by (\ref{def of H}). Recall that $u_{\Xi,\Lambda}=\sum_{i=1}^k\sigma_i$ and $\sigma_i(x)=\sigma_{\xi^i,\lambda_i}(x-P^i)$. Similarly, consider
\begin{equation}\label{def of f}
f_{\Xi,\Lambda}:=\sum_{i=1}^k\eta_i \text{ and } \eta_i(x):=\sum_{\mu\nu}\hat{H}_{\mu\nu}(x-P^i)D_{\mu\nu}\sigma_{\xi^i,\lambda_i}(x-P^i).
\end{equation}

By direct computations, we know that $\norm{X^{0,\alpha}_{k,n-8}}{f_{\Xi,\Lambda}}\leq C$. When $n\geq 25$ and $0<c_0<\frac{n-2}{2}-8$, we can pick $\frac{n-2}{2}<s<n-10-c_0$ such that $s+2<n-8$. So, according to (\ref{step function equation}), we know that
$$
\norm{ Y }{f_{\Xi,\Lambda}}\leq\norm{X^{0,\alpha}_{k,n-8}}{f_{\Xi,\Lambda}}\leq C.
$$
According to Proposition \ref{Surjectivity}, we can find an unique $w_{\Xi,\Lambda}\in X_2$ such that 
\begin{equation}\label{def of w}
L_{\Xi,\Lambda}w_{\Xi,\Lambda}=\Proj_{Y_2}f_{\Xi,\Lambda}.
\end{equation}
Now we can give the definition of an auxiliary function:
\begin{equation*}
G_k(\Xi,\Lambda):=\sum_{i=1}^k\hat{G}(\xi^i,\lambda_i)+\frac{1}{2}\int_{\r^n}f_{\Xi,\Lambda}w_{\Xi,\Lambda}.
\end{equation*}

\begin{prop}\label{I and G_k approxi}
Let $n\geq 25$, $0<c_0<\frac{n-2}{2}-8$, and $\tau_0\in \r$. There exist positive constants $k_0, C$ depending on $n, c_0, \tau_0 $ such that for any $k\geq k_0$, $(\Xi,\Lambda)\in \mathcal{C}_k$,
$$
\left|I(\Xi,\Lambda)-kI_0-\epsilon^2t^{16+2c_0}G_k(\Xi,\Lambda)\right|
\leq C\epsilon^2t^{16+2c_0}\left(t^{\frac{16}{n-2}}+t^{c_0}+t^{\frac{1}{2}(\frac{n-2}{2}-8-c_0)}\right),
$$
where $I_0:=I_{\delta}(\sigma_{\xi,\lambda})$ does not depend on $\xi,\lambda$, so we simply use $I_0$ to denote it. 
\end{prop}

Consider an expansion of $I$:
\begin{equation*}
I(\Xi,\Lambda)=A_1(\Xi,\Lambda)+A_2(\Xi,\Lambda)+A_3(\Xi,\Lambda),
\end{equation*}
where $A_1(\Xi,\Lambda):=I_{g }(u_{\Xi,\Lambda})$, $A_2(\Xi,\Lambda):=I'_g(u_{\Xi,\Lambda})\phi_{\Xi,\Lambda}+\frac{1}{2}I''_{g }(u_{\Xi,\Lambda})(\phi_{\Xi,\Lambda},\phi_{\Xi,\Lambda})$ and
\begin{align*}
A_3(\Xi,\Lambda):=&\frac{(n-2)^2}{2}\int_{\r^n} |u_{\Xi,\Lambda}+\phi_{\Xi,\Lambda}|^{\frac{2n}{n-2}}-u_{\Xi,\Lambda}^{\frac{2n}{n-2}}\\
&-\frac{(n-2)^2}{2}\int_{\r^n}\frac{2n}{n-2}u_{\Xi,\Lambda}^{\frac{n+2}{n-2}}\phi_{\Xi,\Lambda}+\frac{n(n+2)}{(n-2)^2}u_{\Xi,\Lambda}^{\frac{4}{n-2}}\phi_{\Xi,\Lambda}^2.
\end{align*}

Using Lemma \ref{lem1} and Lemma \ref{lem3} with $\tau=\frac{5}{4}$, we have
\begin{equation*}
\begin{aligned}
|A_3(\Xi,\Lambda)|&\leq C\int_{\r^n} |\phi_{\Xi,\Lambda}|^{\frac{2n}{n-2}}\leq Ck\norm{X^0_{k,s}}{\phi_{\Xi,\Lambda}}^{\frac{2n}{n-2}}\int_{\Omega_1}\<x-P^1\>^{\frac{2n(\tau-s)}{n-2}}\\
&\leq Ck\norm{ X }{\phi_{\Xi,\Lambda}}^{\frac{2n}{n-2}}.
\end{aligned}
\end{equation*}

Next, we will divide this subsection into two parts: the estimates of $A_1(\Xi,\Lambda)$ and $A_2(\Xi,\Lambda)$.
 
We will use the following control on $h$ in the balls $B_i:=\{\<x-P^i\><r/k\}$:
\begin{equation}\label{h}
\begin{aligned}
|h(x)|&\leq Ct^{8+c_0}\<x-P^i\>^8\\
|Dh(x)|&\leq Ct^{8+c_0}\<x-P^i\>^7\\
|D^2h(x)|&\leq Ct^{8+c_0}\<x-P^i\>^6.
\end{aligned}
\end{equation}

\subsubsection{The Estimate of $A_1$}
\begin{prop}\label{I_g(u)}
Let $n\geq 20$, $c_0>0$. There exists a positive constant $C=C(n)$ such that for any $k\geq 3$, $(\Xi,\Lambda)\in \mathcal{C}_k$,
$$
\left|A_1(\Xi,\Lambda)-kI_0-\epsilon^2t^{16+2c_0}\sum_{i=1}^k\hat{G}(\xi^i,\lambda_i)\right|\leq Ck(\epsilon^3t^{\frac{2n}{n-2}(8+c_0)}+(k/r)^{n-3}).
$$

\end{prop}
\begin{proof}
For simplicity, we drop the subscript $\Xi,\Lambda$ of $u_{\Xi,\Lambda}$ in this proof. Because of Lemma \ref{lem5}, $\trace h=0$ and $D_{\mu}h_{\mu\nu}=0$, we can divide $A_1$ into four terms: 
$$
A_1(\Xi,\Lambda)=I_{\delta}(u)+\epsilon G_1(u)+\epsilon^2 G_2(u)+O(\epsilon^3 R(u)),$$
where 
\begin{align*}
G_1(u)&:=-\frac{1}{2}\int_{\r^n}h_{\mu\nu}D_{\mu}u D_{\nu}u, \\
G_2(u)&:=\frac{1}{4}\int_{\r^n}(h^2)_{\mu\nu}D_{\mu}u D_{\nu}u-\frac{c(n)}{2}(D_{\alpha}h_{\mu\nu})^2u^2,\\
R(u)&:=\int_{\r^n}|h|^3|Du|^2+(|h|^2|D^2h|+|h||Dh|^2)u^2.
\end{align*}

1. Let us deal with $R(u)$ first. Consider the inside part $B_i=\{\<x-P^i\><r/k\}$ and the outside part $\r^n\setminus\left(\cup_i B_i\right)$.

For the inside part, $\<x-P^i\><r/k\leq Ct^{-1}$. So, by the estimates (\ref{h}) and (\ref{step function gamma}), we have
\begin{equation}\label{inside part}    
\begin{aligned}
&\int_{\cup_i B_i}|h|^3|Du|^2+(|h|^2|D^2h|+|h||Dh|^2)u^2\\
\leq &Ckt^{3(8+c_0)}\int_{B_i}\<x-P^i\>^{-2(n-1)+24}\\
\leq &Ckt^{\frac{2n}{n-2}(8+c_0)}\int_{B_i}\<x-P^i\>^{-2n+2+\frac{16n}{n-2}}\leq Ckt^{\frac{2n}{n-2}(8+c_0)}.
\end{aligned}
\end{equation}

For the outside part $\r^n\setminus\left(\cup_i B_i\right)$, by (\ref{step function gamma}), we have
$$
|D^lu|(x)\leq C\gamma(d(x))d(x)^{-(n-2)-l}, \text{ for } l=0,1,2.
$$
Then by (\ref{scaling coeifficients}) and the above estimate, there is
\begin{align*}
&\int_{\r^n\setminus\left(\cup_i B_i\right)}|h|^3|Du|^2+(|h|^2|D^2h|+|h||Dh|^2)u^2\\
\leq &C\int_{\r^n\setminus\left(\cup_i B_i\right)}\gamma(d(x))^2d(x)^{-2n+2}+t^2\gamma(d(x))^2d(x)^{-2n+4}\\
\leq &Ck\int_{B(0,t^{-1})\cap \Omega_1\setminus B_1} \gamma(\<x-P^1\>)^2\<x-P^1\>^{-2n+2},
\end{align*}
where we have used $\supp h=B(0,t^{-1})$, (\ref{B inside outside argument}) and the symmetry of $\gamma(d(x)), d(x)$ in the second inequality. Translate by $P^1$ and enlarge the domain of integration,
$$
\int_{B(0,t^{-1})\cap \Omega_1\setminus B_1} \gamma(\<x-P^1\>)^2\<x-P^1\>^{-2n+2}dx\leq \int_{r/k\leq \<x\>\leq Ct^{-1}}\gamma(\<x\>)^2\<x\>^{-2n+2}dx.
$$
Change to polar coordinate and change variable $\rho=\<x\>$,
$$
\int_{r/k\leq \<x\>\leq Ct^{-1}}\gamma(\<x\>)^2\<x\>^{-2n+2}dx\leq C \int_{r/k}^{Ct^{-1}}\gamma(\rho)^2\rho^{-n+1}d\rho.
$$
Using the definition of $\gamma(\rho)$, we know that
\begin{align*}
&\int_{r/k}^{Ct^{-1}}\gamma(\rho)^2\rho^{-n+1}d\rho\\
\leq &C(r/k)^{2-n}\left(\sum_{j=1}^{k-1}(j+1)^2(j^{2-n}-(j+1)^{2-n})\right)+C(r/k)^{2-n}k^{4-n}\\
\leq & C(r/k)^{2-n}\left(1+k^{4-n}+\sum_{j=1}^{k-1}(2j+1)j^{2-n}\right)\leq C(k/r)^{n-2}.
\end{align*}
Combining the above four estimates together, we obtain the estimate of the outside part,
\begin{equation}\label{outside part}
\int_{\r^n\setminus\left(\cup_i B_i\right)}|h|^3|Du|^2+(|h|^2|D^2h|+|h||Dh|^2)u^2\leq Ck(k/r)^{n-2}. 
\end{equation}

Combining the inside part (\ref{inside part}) and outside part (\ref{outside part}), we have 
\begin{equation}\label{R(u)term}
|R(u)|\leq Ck(t^{\frac{2n}{n-2}(8+c_0)}+(k/r)^{n-2}).
\end{equation}

2. Let us deal with $G_2(u)$. The proof of the outside part $\r^n\setminus\left(\cup_i B_i\right)$ is the same with that of $R(u)$. For the inside part, we only need to control the interactions between different bubbles. When $x\in B_1$, by (\ref{tau eqn1}) with $\tau=n-3>1$, there are
\begin{align*}
\sum_{j=2}^k\sigma_j\leq &C\sum_{j=2}^k\<x-P^j\>^{-(n-2)}\leq C(k/r)^{n-3}\<x-P^1\>^{-1},\\
\left|D\sum_{j=2}^k\sigma_j\right|\leq &C\sum_{j=2}^k \<x-P^j\>^{-(n-1)}\leq C(k/r)^{n-3}\<x-P^1\>^{-2}.
\end{align*}
So, combining the above two estimates with $\<x-P^1\><r/k\leq Ct^{-1}$, we have
\begin{align*}
&\left|\int_{B_1}(h^2)_{\mu\nu}D_{\mu}u D_{\nu}u-\frac{c(n)}{2}(D_{\alpha}h_{\mu\nu})^2u^2-\int_{B_1}(h^2)_{\mu\nu}D_{\mu}\sigma_1D_{\nu}\sigma_1-\frac{c(n)}{2}(D_{\alpha}h_{\mu\nu})^2\sigma_1^2\right|\\
\leq&C(k/r)^{n-3}\int_{B_1}\<x-P^1\>^{-(n-1)-2}+t^2\<x-P^1\>^{-(n-2)-1}\\
\leq&C(k/r)^{n-3}.
\end{align*}
Combining the inside and outside parts together, we have 
\begin{equation}\label{G_2(u)term}
\left|G_2(u)-t^{16+2c_0}\sum_{i=1}^k\hat{G}(\xi^i,\lambda_i)\right|\leq Ck(k/r)^{n-3}.
\end{equation}

3. Let us deal with $G_1(u)$. Similarly, the outside part $\r^n\setminus\left(\cup_i B_i\right)$ is the same with that of $R(u)$. For the inside part, we can control the interactions between different bubbles. In other words, we have
$$
\left|-\int_{B_1}h_{\mu\nu}D_{\mu}u D_{\nu}u+\int_{B_1}h_{\mu\nu}D_{\mu}\sigma_1D_{\nu}\sigma_1\right|\leq C(k/r)^{n-3}.
$$
Next, we just need to estimate $\left|\int_{B_1}h_{\mu\nu}D_{\mu}\sigma_1D_{\nu}\sigma_1\right|$. Consider an equality 
$$
D_{\mu}\sigma_1D_{\nu}\sigma_1-c(n)D_{\mu\nu}(\sigma_1^2)=\frac{1}{n}(|D\sigma_1|^2-c(n)\Delta(\sigma_1^2))\delta_{\mu\nu}.
$$
Multiply the equality above by $h_{\mu\nu}$ and integrate over $\r^n$, do integration by parts and use the two conditions of $h$: $\trace h=0$ and $D_{\mu}h_{\mu\nu}=0$. Then,
$$
\int_{\r^n}h_{\mu\nu}D_{\mu}\sigma_1D_{\nu}\sigma_1=0.
$$
So, we can compute 
$$
\left|\int_{B_1}h_{\mu\nu}D_{\mu}\sigma_1D_{\nu}\sigma_1\right|=\left|\int_{\r^n\setminus B_1}h_{\mu\nu}D_{\mu}\sigma_1D_{\nu}\sigma_1\right|\leq C\int_{\r^n\setminus B_1}\<x\>^{-2n+2}\leq C(k/r)^{n-2}.
$$
Combining the inside and outside parts together, we have 
\begin{equation}\label{G_1(u)term}
\left|G_1(u)\right|\leq Ck(k/r)^{n-3}.
\end{equation}

4. Let us deal with $I_{\delta}(u)$. Because $-\Delta \sigma_j=n(n-2)\sigma_j^{\frac{n+2}{n-2}}$ and $u=\sum_i\sigma_i$, we have 
\begin{align*}
I_{\delta}(u)&=\frac{1}{2}\int_{\r^n}|Du|^2-\frac{(n-2)^2}{2}\int_{\r^n}u^{\frac{2n}{n-2}}\\
&=\frac{n(n-2)}{2}\sum_{i,j}\int_{\r^n}\sigma_i^{\frac{n+2}{n-2}}\sigma_j-\frac{(n-2)^2}{2}\int_{\r^n}(\sum_j\sigma_j)^{\frac{2n}{n-2}}.
\end{align*}
For the first term, use (\ref{tau eqn2}) with $\tau=n-2>1$, 
\begin{align*}
&\left|\sum_{i,j}\int_{\r^n}\sigma_i^{\frac{n+2}{n-2}}\sigma_j-k\int_{\r^n}\sigma^{\frac{2n}{n-2}}\right|\\
\leq &\left|\sum_{i\neq j}\int_{\r^n}\sigma_i^{\frac{n+2}{n-2}}\sigma_j\right|\leq C(k/r)^{n-2}\int_{\r^n}\sum_i\<x-P^i\>^{-(n+2)}\leq Ck(k/r)^{n-2}. 
\end{align*}
For the second term, by symmetry, we estimate the integration on $\Omega_1$. Using Lemma \ref{lem1}, and then using (\ref{tau eqn1}) with $\tau=n-2, \frac{n-2}{2}, \frac{(n-2)^2}{2n}$ in order, we have
\begin{align*}
&\left|\int_{\Omega_1}(\sum_j\sigma_j)^{\frac{2n}{n-2}}-\int_{\r^n}\sigma^{\frac{2n}{n-2}}\right|\\
&\leq\int_{\r^n\setminus\Omega_1}\sigma_1^{\frac{2n}{n-2}}+C\int_{\Omega_1}\sigma_1^{\frac{n+2}{n-2}}\sum_{j=2}^k\sigma_j+\sigma_1^{\frac{4}{n-2}}(\sum_{j=2}^k\sigma_j)^{2}+(\sum_{j=2}^k\sigma_j)^{\frac{2n}{n-2}}\\
&\leq \int_{\r^n\setminus B_1}\<x-P^1\>^{-2n}+C(k/r)^{n-2}\int_{\Omega_1}\<x-P^1\>^{-(n+2)}\leq C(k/r)^{n-2}.
\end{align*}
Combining the two terms together, we obtain 
\begin{equation}\label{I_0(u)term}
|I_{\delta}(u)-kI_{\delta}(\sigma)|\leq Ck(k/r)^{n-2}.
\end{equation}

Combining (\ref{R(u)term}), (\ref{G_2(u)term}), (\ref{G_1(u)term}), (\ref{I_0(u)term}) together, we can obtain the estimate in Proposition \ref{I_g(u)}.
\end{proof}

\subsubsection{The Estimate of $A_2$}

\begin{prop}\label{second part}
Let $n\geq 25$, $0<c_0<\frac{n-2}{2}-8$, and $\tau_0\in \r$. There exist positive constants $k_0, C$ depending on $n,c_0,\tau_0$ such that for any $k\geq k_0$, $(\Xi,\Lambda)\in \mathcal{C}_k$,
$$
\left|A_2(\Xi,\Lambda)+\frac{\epsilon^2t^{16+2c_0}}{2}\int_{\r^n}f_{\Xi,\Lambda}w_{\Xi,\Lambda}\right|
\leq C\epsilon^2 t^{16+2c_0}\left(t^{\frac{16}{n-2}}+t^{c_0}+t^{\frac{1}{2}(\frac{n-2}{2}-8-c_0)}\right),
$$
where $f_{\Xi,\Lambda}$ and $w_{\Xi,\Lambda}$ are defined as equations (\ref{def of f}) and (\ref{def of w}).
\end{prop}
Recall that $A_2(\Xi,\Lambda)=I'_g(u_{\Xi,\Lambda})\phi_{\Xi,\Lambda}+\frac{1}{2}I''_{g }(u_{\Xi,\Lambda})(\phi_{\Xi,\Lambda},\phi_{\Xi,\Lambda})$. Integrating by parts, we have the two terms, 
$$
I'_{g}(u_{\Xi,\Lambda})\phi_{\Xi,\Lambda}=\int_{\r^n}F(u_{\Xi,\Lambda})\phi_{\Xi,\Lambda} \text{ and }I''_{g }(u_{\Xi,\Lambda})(\phi_{\Xi,\Lambda},\phi_{\Xi,\Lambda})=\int_{\r^n}\phi_{\Xi,\Lambda}L_{\Xi,\Lambda}\phi_{\Xi,\Lambda}.
$$
So, we need to take out the leading terms from $\phi_{\Xi,\Lambda}$ and $F(u_{\Xi,\Lambda})$. 

The way to take out the leading term of $F(u_{\Xi,\Lambda})$ is purely computational, see the following lemma. 
\begin{lem}\label{basicproperty2}
Let $n\geq 25$, $0<c_0<\frac{n-2}{2}-8$, $\frac{n-2}{2}<s<n-10-c_0$, $0<\alpha<\frac{4}{n-2}$, and $\tau_0\in \r$. There exist positive constants $k_0, C$ depending on $n,s,\alpha,c_0,\tau_0$ such that for any $k\geq k_0$, $(\Xi,\Lambda)\in \mathcal{C}_k$, we have
\begin{align*}
\norm{ Y }{F(u_{\Xi,\Lambda})}&\leq C(\epsilon t^{8+c_0}+(k/r)^{n-2-s}),\\
\norm{ Y }{F(u_{\Xi,\Lambda})+\epsilon t^{8+c_0}f_{\Xi,\Lambda}}&\leq C(\epsilon t^{8+c_0}(k(k/r)^{c_0}+t^{c_0})+(k/r)^{n-2-s}).
\end{align*}	
\end{lem}
\begin{proof}
Again, for simplicity, we drop the subscript $\Xi,\Lambda$ of $u_{\Xi,\Lambda}$, $f_{\Xi,\Lambda}$ in this proof.

After using $-\Delta \sigma_j=n(n-2)\sigma_j^{\frac{n+2}{n-2}}$, both $F(u)$ and $F(u)+\epsilon t^{8+c_0}f$ have the same nonlinear term, and we can use (\ref{smallnesssup2}) and (\ref{smallnessHolder2}) to control the nonlinear term 
$$
\norm{ Y }{\sum_{i=1}^k\sigma_i^{\frac{n+2}{n-2}}-(\sum_{i=1}^k\sigma_i)^{\frac{n+2}{n-2}}}\leq C(k/r)^{n-s}.
$$

The linear term of $F(u)$ is $(L_{g }-\Delta)u$. The linear term of $F(u)+\epsilon t^{8+c_0}f$ is $(L_{g }-\Delta)u+t^{8+c_0} f$. The estimates of H{\"o}lder semi-norms of the linear terms are similar to that of the sup-norms of the linear terms. The only difference is that we need the product rule and interpolation inequalities from Lemma \ref{holder semi norm lemma}. Readers can see more details in the proof of Proposition \ref{Welldefined}. In this proof, we just focus on the sup-norms of the linear terms.

For the linear term of $F(u)$, we know that 
$$
|(L_{g }-\Delta)u|\leq C\epsilon (|h||D^2u|+|Dh||Du|+|D^2h|u).
$$
We consider two parts $B_i=\{\<x-P^i\><r/k\}$ and $\r^n\setminus\left(\cup_i B_i\right)$. Because $\supp h\subseteq B(0,t^{-1})$, we just consider the outside part as $B(0,t^{-1})\setminus\left(\cup_i B_i\right)$. 

When $x\in B_1$, using (\ref{h}), and $s+2<n-8$, we have
\begin{align*}
|(L_{g }-\Delta)u|(x)\leq& C\epsilon t^{8+c_0}\<x-P^1\>^{8-n}\leq C\epsilon t^{8+c_0}\<x-P^1\>^{-s-2}\\
\leq &C\epsilon t^{8+c_0}\sum_i\<x-P^i\>^{-s-2}.
\end{align*}

When $x\in B(0,t^{-1})\setminus\left(\cup_i B_i\right)$, using (\ref{coefficient}), (\ref{step function gamma}), (\ref{B inside outside argument}), $d(x)\geq r/k$, and (\ref{step function gamma}) again, we have
\begin{align*}
|(L_{g }-\Delta)u|(x)&\leq C\epsilon \gamma(d(x))\left(d(x)^{-n}+td(x)^{1-n}+t^2d(x)^{2-n}\right)\\
&\leq C\epsilon \gamma(d(x))d(x)^{-n}\leq C\epsilon (k/r)^{n-2-s}\sum_i\<x-P^i\>^{-s-2}.
\end{align*}

In conclusion, we have
$$
\norm{ Y }{F(u)}\leq C(\epsilon t^{8+c_0}+(k/r)^{n-2-s}).
$$

For the linear term of $F(u)+\epsilon t^{8+c_0}f$, because $\trace h=0$ and $D_{\mu}h_{\mu\nu}=0$, we know that
\begin{align*}
\left|(L_{g }-\Delta)u+\epsilon t^{8+c_0} f\right|\leq &C\epsilon|h_{\mu\nu}D_{\mu\nu}u-t^{8+c_0}f|\\
+&C\epsilon^2(|h|^2|D^2u|+|h||Dh||Du|+|h||D^2h|u+|Dh|^2u).
\end{align*}

When $x\in B_1$, using (\ref{h}), (\ref{B inside outside argument}) and $s+2<n-8$, we can control 
\begin{align*}
&(|h|^2|D^2u|+|h||Dh||Du|+|h||D^2h|u+|Dh|^2u)(x)\\
\leq &Ct^{16+2c_0}\<x-P^1\>^{16-n}\leq Ct^{8+2c_0}\<x-P^1\>^{8-n}\\
\leq &Ct^{8+2c_0}\sum_i\<x-P^i\>^{-s-2}.
\end{align*}
When $x\in B_1$, using (\ref{tau eqn1}) with $\tau=c_0>0$ and $s+2<n-8-c_0$, we have
\begin{align*}
|h_{\mu\nu}D_{\mu\nu}u-t^{8+c_0}f|(x)\leq &Ct^{8+c_0}\sum_{i=2}^{k}\left|\hat{H}_{\mu\nu}(x-P^1)-\hat{H}_{\mu\nu}(x-P^i)\right| |D_{\mu\nu}\sigma_i|\\
\leq &Ct^{8+c_0}\left(\sum_{i=2}^k\<x-P^i\>^{8-n}+\<x-P^1\>^8\sum_{i=2}^k\<x-P^i\>^{-n}\right)\\
\leq &Ct^{8+c_0}k(k/r)^{c_0}\<x-P^1\>^{8-n+c_0}\\
\leq &Ct^{8+c_0}k(k/r)^{c_0}\sum_i\<x-P^i\>^{-s-2}.
\end{align*}
Combining the two terms together, when $x\in B_1$, we have
$$
|(L_{g }-\Delta)u+\epsilon t^{8+c_0}f|(x)\leq C\epsilon t^{8+c_0}(k(k/r)^{c_0}+t^{c_0})\sum_i\<x-P^i\>^{-s-2}.
$$

When $x\in B(0,t^{-1})\setminus(\cup_iB_i)$, using (\ref{coefficient}), (\ref{step function gamma}), (\ref{B inside outside argument}), $d(x)\geq r/k$, and (\ref{step function gamma}), we have
\begin{align*}
&|(L_{g }-\Delta)u+\epsilon t^{8+c_0}f|(x)\\
\leq &C\epsilon \gamma(d(x))\left(d(x)^{-n}+td(x)^{1-n}+t^2d(x)^{2-n}+t^8d(x)^{8-n}\right),\\
\leq &C\epsilon \gamma(d(x))d(x)^{-n}\leq C\epsilon (k/r)^{n-2-s}\sum_i\<x-P^i\>^{-s-2}.
\end{align*}

In conclusion, we have
$$
\norm{ Y }{F(u)+\epsilon t^{8+c_0}f}\leq C(\epsilon t^{8+c_0}(k(k/r)^{c_0}+t^{c_0})+(k/r)^{n-2-s}).
$$
\end{proof}

The way to take out the leading term of $\phi_{\Xi,\Lambda}$ is by so called Newton's method.
\begin{lem}\label{u-v-w}
Let $n\geq 25$, $0<c_0<\frac{n-2}{2}-8$, $\frac{n-2}{2}<s<n-10-c_0$, $0<\alpha<\frac{4}{n-2}$, and $\tau_0\in \r$. There exist positive constants $k_0, C$ depending on $n,s,\alpha,c_0,\tau_0$ such that for any $k\geq k_0$, $(\Xi,\Lambda)\in \mathcal{C}_k$, we have
\begin{align*}
\norm{ X }{\phi_{\Xi,\Lambda}}&\leq C(\epsilon t^{8+c_0}+(k/r)^{n-2-s}),\\
\norm{ X }{\phi_{\Xi,\Lambda}-\epsilon t^{8+c_0}w_{\Xi,\Lambda}}&\leq C(\epsilon t^{8+c_0}((\epsilon t^8)^{\frac{4}{n-2}-\alpha}+k(k/r)^{c_0}+t^{c_0})+(k/r)^{n-2-s}).
\end{align*}
\end{lem}
\begin{proof}
For simplicity, we drop the subscript $\Xi,\Lambda$ of $w_{\Xi,\Lambda}$, $\phi_{\Xi,\Lambda}$ in this proof. From Proposition \ref{solvability}, we know that $\phi\in X_2$ solves 
$$
\Proj_{Y_2}F(u_{\Xi,\Lambda}+\phi)=0.
$$
Because of estimate (\ref{phi & F(u) same order}) and Lemma \ref{basicproperty2}, we have
$$
\norm{ X }{\phi}\leq C(\epsilon t^{8+c_0}+(k/r)^{n-2-s}).
$$

Next, we will prove the second estimate. Recall that $w\in X_2$ satisfies the following equation
$$
L_{\Xi,\Lambda}w=\Proj_{Y_2}f_{\Xi,\Lambda}.
$$
Then $\phi-\epsilon t^{8+c_0}w\in X_2$ and we can compute $L_{\Xi,\Lambda}(\phi-\epsilon t^{8+c_0}w)$:
\begin{equation}\label{identity of phi-w}
\begin{aligned}
&L_{\Xi,\Lambda}(\phi-\epsilon t^{8+c_0}w)\\
=&-\Proj_{Y_2}(F(u_{\Xi,\Lambda})+\epsilon t^{8+c_0}f_{\Xi,\Lambda})-\Proj_{Y_2}(F(u_{\Xi,\Lambda}+\phi)-F(u_{\Xi,\Lambda})-L_{\Xi,\Lambda}\phi).
\end{aligned}
\end{equation}

Firstly, from Lemma \ref{basicproperty2}, 
$$
\norm{ Y }{F(u_{\Xi,\Lambda})+\epsilon t^{8+c_0}f_{\Xi,\Lambda}}\leq C(\epsilon t^{8+c_0}(k(k/r)^{c_0}+t^{c_0})+(k/r)^{n-2-s}).
$$
Secondly, from Lemma \ref{basicproperty2} and Proposition \ref{Welldefined},
$$
\norm{ Y }{F(u_{\Xi,\Lambda}+\phi)-F(u_{\Xi,\Lambda})-L_{\Xi,\Lambda}\phi}\leq C\norm{ X }{\phi}^{\frac{n+2}{n-2}-\alpha}\leq C(\epsilon t^{8+c_0}+(k/r)^{n-2-s})^{\frac{n+2}{n-2}-\alpha}.
$$
Then, combining the above two estimates with (\ref{identity of phi-w}) and Proposition \ref{Surjectivity}, we have
$$
\norm{ X }{\phi-\epsilon t^{8+c_0}w}\leq C(\epsilon t^{8+c_0}((\epsilon t^8)^{\frac{4}{n-2}-\alpha}+k(k/r)^{c_0}+t^{c_0})+(k/r)^{n-2-s}).
$$
\end{proof}

Now, we are ready to prove Proposition \ref{second part}.
\begin{proof}
For simplicity, we drop the subscript $\Xi,\Lambda$ of $w_{\Xi,\Lambda}$, $\phi_{\Xi,\Lambda}, u_{\Xi,\Lambda}$, $f_{\Xi,\Lambda}$ in this proof. We have
\begin{align*}
&\left|\int_{\r^n}F(u)\phi+\epsilon^2t^{16+2c_0}\int_{\r^n}f w\right|\\
\leq &\left|\int_{\r^n}F(u)\phi+\int_{\r^n}\epsilon t^{8+c_0}f\phi \right|+\left|\epsilon^2t^{16+2c_0}\int_{\r^n}fw-\int_{\r^n}\epsilon t^{8+c_0}f\phi\right|\\
\leq &\left|\int_{\r^n}\left(F(u)+\epsilon t^{8+c_0}f\right)\phi\right| +\left|\int_{\r^n}\epsilon t^{8+c_0}f(\phi-\epsilon t^{8+c_0}w)\right|.
\end{align*}
Then using (\ref{dual space}) and the relation $r=e^{k}/k$, we can obtain
\begin{align*}
&\left|\int_{\r^n}F(u)\phi+\epsilon^2t^{16+2c_0}\int_{\r^n}f w\right|\\
\leq &C\left(\norm{ X }{\phi}\norm{ Y }{F(u)+\epsilon t^{8+c_0}f}+\epsilon t^{8+c_0}\norm{ X }{\phi-\epsilon t^{8+c_0}w}\norm{ Y }{f}\right).
\end{align*}

Similarly, we can use the same method to obtain,
\begin{align*}
&\left|\int_{\r^n}\phi L_{\Xi,\Lambda}\phi-\epsilon^2t^{16+2c_0}\int_{\r^n}f w\right|\\
\leq &C\left(\norm{ X }{L_{\Xi,\Lambda}\phi}\norm{ X }{\phi-\epsilon t^{8+c_0}w}+\epsilon t^{8+c_0}\norm{ X }{\phi-\epsilon t^{8+c_0}w}\norm{ Y }{f}\right).
\end{align*}

Combining the two estimates together, using Lemmas \ref{basicproperty2} and \ref{u-v-w}, we can obtain 
\begin{align*}
&\left|A_2(\Xi,\Lambda)+\frac{\epsilon^2t^{16+2c_0}}{2}\int_{\r^n}f_{\Xi,\Lambda}w_{\Xi,\Lambda}\right|\\
\leq &C(\epsilon^2 t^{16+2c_0}((\epsilon t^8)^{\frac{4}{n-2}-\alpha}+k(k/r)^{c_0}+t^{c_0})+\epsilon t^{8+c_0}(k/r)^{n-2-s}+(k/r)^{2(n-2-s)}).
\end{align*}
Recall that $\frac{n-2}{2}<s<n-10-c_0$, $0<\alpha<\frac{4}{n-2}$. Pick
$$
s=\frac{1}{2}\left(\frac{n-2}{2}+n-10-c_0\right)\text{ and } \alpha=\frac{2}{n-2}.
$$
Then picking $k_0$ large, for any $k\geq k_0$, we obtain
$$
\left|A_2(\Xi,\Lambda)+\frac{\epsilon^2t^{16+2c_0}}{2}\int_{\r^n}f_{\Xi,\Lambda}w_{\Xi,\Lambda}\right|\leq C\epsilon^2t^{16+2c_0}\left(t^{\frac{16}{n-2}}+t^{c_0}+t^{\frac{1}{2}(\frac{n-2}{2}-8-c_0)}\right).
$$
\end{proof}

\subsection{Proof of the Main Theorem}

\begin{lem}\label{G&F}
Let $n\geq 25$, $0<c_0<\frac{n-2}{2}-8$. There exist $\tau_0=\tau_0(n)\in \r$, $\epsilon=\epsilon(n, c_0)>0$, and $k_0=k_0(n, c_0)>0$ such that for any $k\geq k_0$,
$$
\min_{\partial B((\Xi_0,\Lambda_0),k^{-1})} G_k(\Xi,\Lambda)\geq G_k(\Xi_0,\Lambda_0)+\frac{1}{C}k^{-2},
$$
where $\Lambda_0=(1,\cdots,1)$, $\Xi_0=(0,\cdots,0)$.
\end{lem}
\begin{proof}
Recall the definition of the auxiliary function:
$$G_k(\Xi,\Lambda)=A_k(\Xi,\Lambda)+B_k(\Xi,\Lambda),$$
where $A_k(\Xi,\Lambda):=\sum_{i=1}^k\hat{G}(\xi^i,\lambda_i)$ and $B_k(\Xi,\Lambda):=\frac{1}{2}\int_{\r^n}f_{\Xi,\Lambda}w_{\Xi,\Lambda}$.

To deal with $A_k(\Xi,\Lambda)$, we use \cite[Corollary 23]{brendle2009blow}. The coefficient $\tau_0$ in the definition (\ref{def of H}) can be chosen such that $\hat{G}(1,0)<0$, $\nabla_{\xi,\lambda} \hat{G}(1,0)=0$. And most importantly, there exist positive constants $\delta$ and $C$ such that for any $|(\xi,\lambda)-(1,0)|<\delta$, the Hessian matrix of $\hat{G}$ is positive definite with
$$
\nabla^2_{\xi,\lambda} \hat{G}(\xi,\lambda)\geq \frac{1}{C}I_{(n+1)\times (n+1)}.
$$
Next, we compute the Hessian matrix of $A_k(\Xi,\Lambda)$. For any $i=1,\cdots,k$, the block matrices on the diagonal,
$$
\nabla^2_{\xi^i,\lambda_i}A_k(\Xi,\Lambda)=\nabla^2_{\xi,\lambda}\hat{G}(\xi^i,\lambda_i).
$$ 
All the other second order derivatives are zero. In other words,
$\partial_{\lambda_i}\partial_{\lambda_j}A_k(\Xi,\Lambda)=0$, 
$\partial_{\lambda_i}\partial_{\xi^j_{\mu}}A_k(\Xi,\Lambda)=0$, and $\partial_{\xi^i_{\nu}}\partial_{\xi^j_{\mu}}A_k(\Xi,\Lambda)=0$, for $i\neq j$. 

Consider $|(\Xi,\Lambda)|^2:=\sum_i|(\xi^i,\lambda_i)|^2$. Then we have that for any $|(\Xi,\Lambda)-(\Xi_0,\Lambda_0)|<\delta$, where $\Lambda_0=(1,\cdots,1)$, $\Xi_0=(0,\cdots,0)$, the Hessian matrix of $A_k$ is positive definite with
\begin{equation}\label{hessian estimate}
\nabla^2_{\Xi,\Lambda} A_k(\Xi,\Lambda)\geq \frac{1}{C}I_{k(n+1)\times k(n+1)}.
\end{equation}

Because of $ \nabla_{\Xi,\Lambda} A_k(\Xi_0,\Lambda_0)=0$ and the Hessian estimate (\ref{hessian estimate}), we can use Taylor expansion of $A_k$ at $(\Xi_0,\Lambda_0)$ to obtain that for any $|(\Xi,\Lambda)-(\Xi_0,\Lambda_0)|<\delta$,
\begin{equation}\label{A Lam Xi}
A_k(\Xi,\Lambda)\geq A_k(\Xi_0,\Lambda_0)+\frac{1}{C}|(\Xi,\Lambda)-(\Xi_0,\Lambda_0)|^2.
\end{equation}

Secondly, we deal with $B_k(\Xi,\Lambda)$. Recall $f_{\Xi,\Lambda}=\sum_i\eta_i$, where 
$$
\eta_i(x)=\lambda_i^{\frac{n-2}{2}}\hat{H}_{\mu\nu}(x-P^i)D_{\mu\nu}(1+\lambda_i^2|x-P^i-\xi^i|^2)^{-\frac{n-2}{2}}.
$$
Then we can use two properties of $\hat{H}$: $x_{\mu}\hat{H}_{\mu\nu}(x)=0$ and $\trace \hat{H}(x)=0$ to obtain
\begin{align*}
\eta_i(x+P^i)&=\frac{n(n-2)}{4}\lambda_i^{\frac{n-2}{2}}(1+\lambda_i^2|x-\xi^i|^2)^{-\frac{n+2}{2}}\hat{H}_{\mu\nu}(x)(x_{\mu}-\xi^i_{\mu})(x_{\nu}-\xi^{i}_{\nu})\\
&=\frac{n(n-2)}{4}\lambda_i^{\frac{n-2}{2}}(1+\lambda_i^2|x-\xi^i|^2)^{-\frac{n+2}{2}}\hat{H}_{\mu\nu}(x)\xi^i_{\mu}\xi^{i}_{\nu}.
\end{align*}
So, for any $(\Xi,\Lambda)\in\mathcal{C}_k$,
$$
\norm{ Y }{f_{\Xi,\Lambda}}\leq C\norm{X^{0,\alpha}_{k,n-6}}{f_{\Xi,\Lambda}}\leq C\max_{1\leq i\leq k}|\xi^i|^2\leq C|(\Xi,\Lambda)-(\Xi_0,\Lambda_0)|^2.
$$
Recall that $L_{\Xi,\Lambda}w_{\Xi,\Lambda}=\Proj_{Y_2}f_{\Xi,\Lambda}$. According the estimate (\ref{invertibility}), we have
$$
\norm{ X }{w_{\Xi,\Lambda}}\leq C\norm{ Y }{f_{\Xi,\Lambda}}\leq C|(\Xi,\Lambda)-(\Xi_0,\Lambda_0)|^2.
$$
Then because of the inequality (\ref{dual space}), we know that for any $(\Xi,\Lambda)\in\mathcal{C}_k$,
\begin{equation}\label{B Lam Xi}
|B_k(\Xi,\Lambda)|\leq Ck|(\Xi,\Lambda)-(\Xi_0,\Lambda_0)|^4.
\end{equation}

Combining estimates (\ref{A Lam Xi}) and (\ref{B Lam Xi}), picking $k_0$ large, we have that for any $k\geq k_0$,
$$
\min_{\partial B((\Xi_0,\Lambda_0),k^{-1})} G_k(\Xi,\Lambda)\geq A_k(\Xi_0,\Lambda_0)+\frac{1}{C}k^{-2}=G_k(\Xi_0,\Lambda_0)+\frac{1}{C}k^{-2},
$$
where we have used $B_k(\Xi_0,\Lambda_0)=0$ in the equality above. 
\end{proof}

Now, we are ready to prove Theorem \ref{blowup3}.
\begin{thm}
Let $n\geq 25$, $0<c_0<\frac{n-2}{2}-8$. There exist $\tau_0=\tau_0(n)\in \r$, $\epsilon=\epsilon(n, c_0)>0$, and $k_0=k_0(n, c_0)>0$ such that for any $k\geq k_0$, equation (\ref{afterscale}) has a positive smooth solution $v_{k}$ of the form
$$
v_k=u_{\bar{\Xi},\bar{\Lambda}}+\phi_k,
$$ 
where $(\bar{\Xi},\bar{\Lambda})\in \mathcal{C}_k$. Moreover, we have the following estimates: for any $i=1,\cdots,k$, 
\begin{align*}
I_{g}(v_k)&<kI_0,\\
\norm{L^{\frac{2n}{n-2}}(\r^n)}{\phi_k}&\leq C t^{8+c_0},\\
\left|I_{g}(v_k)-kI_0\right|&\leq Ct^{16+2c_0},\\
\sup_{\<x-P^i\><r/k}v_k(x)&\geq C^{-1}.
\end{align*}
\end{thm}
\begin{proof}
We fix $\alpha=\frac{2}{n-2}$, $s=\frac{1}{2}\left(\frac{n-2}{2}+n-10-c_0\right)$, and $\epsilon=\epsilon_0(n,c_0)$ from Corollary \ref{deduce to critical pt}. According to Proposition \ref{solvability}, we know that there exists $k_0=k_0(n,c_0)$ such that for any $k\geq k_0$, $(\Xi,\Lambda)\in \mathcal{C}_k$, there exists $\phi_{\Xi,\Lambda}\in X_2$, defining $v_{\Xi,\Lambda}:=u_{\Xi,\Lambda}+\phi_{\Xi,\Lambda}$, such that
\begin{equation}\label{last equation}
\Proj_{Y_2}\circ F(v_{\Xi,\Lambda})=0.
\end{equation}

According to Proposition \ref{I and G_k approxi}, for any $(\Xi,\Lambda)\in \mathcal{C}_k$,
$$
I(\Xi,\Lambda)=kI_0+\epsilon^2 t^{16+2c_0}(G_k(\Xi,\Lambda)+o(t^{\tau})),
$$
where $\tau=\min\{\frac{16}{n-2},\frac{1}{2}(\frac{n-2}{2}-8-c_0),c_0\}$.
Then pick $\tau_0$ from Lemma \ref{G&F} and pick $k_0$ larger. There exists $C_1>0$ such that for $k\geq k_0$,
$$
\min_{\partial B((\Xi_0,\Lambda_0),k^{-1})} I(\Xi,\Lambda)\geq I(\Xi_0,\Lambda_0)+\frac{1}{C_1}\epsilon^2t^{16+2c_0}k^{-2}.
$$
Then there exists $(\bar{\Xi},\bar{\Lambda})\in B((\Xi_0,\Lambda_0),k^{-1})$ such that it is a local minimal point of $I(\Xi,\Lambda)$. 

In conclusion, because of Corollary \ref{deduce to critical pt}, we know that 
\begin{equation*}
(\Id-\Proj_{Y_2})\circ F(v_{\bar{\Xi},\bar{\Lambda}})=0.
\end{equation*}
Then picking $v_k=v_{\bar{\Xi},\bar{\Lambda}}$ and $\phi_k=\phi_{\bar{\Xi},\bar{\Lambda}}$, combining with (\ref{last equation}), we can solve equation (\ref{afterscale}).

For simplicity, we will drop the subscript $\bar{\Xi},\bar{\Lambda}$ for $v_{\bar{\Xi},\bar{\Lambda}}$, $u_{\bar{\Xi},\bar{\Lambda}}$, and $\phi_{\bar{\Xi},\bar{\Lambda}}$.

Next, we will prove that $v$ is smooth positive. First, if $v$ has non-trivial negative part, by scaling $h(x)=\hat{h}(tx)$, $\supp \hat{h}\subseteq B(0,1)$, and H{\"o}lder inequality, we have that 
\begin{align*}
&\left|\int_{\{v<0\}}|D_{g }v|^2-|Dv|^2+c(n)R_{g }v^2\right|\\
&\leq C\epsilon\norm{C^{2,\alpha}}{\hat{h}}\int_{\{v<0\}}|Dv|^2+\left(\int_{\{v<0\}}|R_g|^{\frac{n}{2}}\right)^{\frac{2}{n}}\left(\int_{\{v<0\}}|v|^{\frac{2n}{n-2}}\right)^{\frac{n-2}{n}}\\
&\leq C\epsilon\norm{C^{2,\alpha}}{\hat{h}}\int_{\{v<0\}}|Dv|^2.  
\end{align*}
Then, using $\norm{C^{2,\alpha}}{\hat{h}}\leq 1$, picking $\epsilon_0$ from Corollary \ref{deduce to critical pt} small, we obtain
$$
\int_{\{v<0\}}|D_{g }v|^2+c(n)R_{g }v^2\geq \frac{1}{C}\norm{L^{\frac{2n}{n-2}}(\{v<0\} )}{v}^2.
$$
Since equation (\ref{afterscale}), we have
$$
\int_{\{v<0\} }|D_{g }v|^2+c(n)R_{g }v^2=n(n-2)\int_{\{v<0\} }|v|^{\frac{2n}{n-2}}=C\norm{L^{\frac{2n}{n-2}}(\{v<0\} )}{v}^{\frac{2n}{n-2}}.
$$
Combining the two estimates above, we obtain
\begin{equation}\label{lower bound in last proof}
\norm{L^{\frac{2n}{n-2}}(\{v<0\} )}{v}\geq \frac{1}{C_2}.
\end{equation}
In the other direction, since $u>0$, Lemma \ref{u-v-w}, and (\ref{D^12 embedding}), picking $k_0$ large, we have
$$
\norm{L^{\frac{2n}{n-2}}(\{v<0\})}{v}\leq \norm{L^{\frac{2n}{n-2}}(\{v<0\})}{v-u}\leq C\norm{L^{2}(\{v<0\})}{D\phi }\leq C\norm{ X }{\phi}\leq \frac{1}{2C_2}.
$$
It is contradicted with (\ref{lower bound in last proof}). So we have proved $v>0$. Then the smoothness of $v$ follows from Trudinger \cite[Theorem 3]{trudinger1968remarks}. 

Second, we will prove the integral estimates. By Lemma \ref{u-v-w} and (\ref{D^12 embedding}), picking $k_0$ large, we have
$$
\norm{L^{\frac{2n}{n-2}}(\r^n)}{\phi}\leq C\norm{\mathcal{D}^{1,2}}{\phi}\leq C\norm{X}{\phi}\leq C t^{8+c_0}.
$$
According to Proposition \ref{I and G_k approxi}, picking $k_0$ large, we have an energy estimate
$$
\left|I(\bar{\Xi},\bar{\Lambda})-kI_0\right|\leq Ct^{16+2c_0}(1+|G_k(\Xi_0,\Lambda_0)|)\leq C t^{16+2c_0}.
$$
According to Proposition \ref{I and G_k approxi} and $G_k(\Xi_0,\Lambda_0)=k\hat{G}(1,0)<0$, picking $k_0$ large, we obtain the estimate 
$$
I_{g}(v)=I(\bar{\Xi},\bar{\Lambda})<kI_0.
$$

In the end, we will prove the sup estimate. For any $i=1,\cdots,k$, consider $B_i=\{\<x-P^i\><r/k\}$. Because of (\ref{tau eqn1}), we have
$$
\sup_{B_i} u\geq C^{-1}\sup_{B_i} \<x-P^i\>^{-(n-2)}=\frac{1}{C_3}.
$$
Also, when $x\in B_i$, using Lemma \ref{u-v-w}, picking $k_0$ large, we have the pointwise difference
$$
|v-u|(x)\leq C\norm{ X }{\phi}\<x-P^i\>^{-s}\leq \frac{1}{2C_3}.
$$
Combining the above two estimates, we have $\sup_{B_i} v\geq \frac{1}{2C_3}$.
\end{proof}

\appendix

\section{H{\"o}lder Semi-norms}
\begin{lem}\label{holder semi norm lemma}
Let $0<\alpha<\beta<1$, $\delta>0$, $B\subseteq \r^n$ be an open ball with radius $R>0$. There exists a positive constant $C=C(n)$ such that for any $u,v\in C^1(B)$,
\begin{description}
    \item [Triangle Inequality] $[u+v]_{\alpha,B}\leq [u]_{\alpha,B}+[v]_{\alpha,B}$;
    \item [Product Rule] $[uv]_{\alpha,B}\leq |u|_B[v]_{\alpha,B}+|v|_B[u]_{\alpha,B}$;
    \item [Interpolation Inequality 1] $[u]_{\alpha,B}\leq C\left((\delta R)^{-\alpha}|u|_B+(\delta R)^{1-\alpha}|Du|_B\right)$;
    \item [Interpolation Inequality 2] $[u]_{\alpha,B}\leq C\left((\delta R)^{-\alpha}|u|_B+(\delta R)^{\beta-\alpha}[u]_{\beta,B}\right)$.
\end{description}
\end{lem}

\section{Energy Expansion}
Let us consider the metric $g=\exp(\epsilon h)$, where $h$ is a symmetric two tensor. In this appendix, the positive constant in the $O(\epsilon^3)$ notation only depends on $n$ and $\norm{C^3}{h}$. Then we have the following expansion:
\begin{enumerate}
\item $g_{\mu\nu}=\delta_{\mu\nu}+\epsilon h_{\mu\nu}+\frac{\epsilon^2}{2}h_{\mu\alpha}h_{\alpha\nu}+O(\epsilon^3)$, $\mu,\nu,\alpha=1,\cdots,n$.
\item $g^{\mu\nu}=\delta_{\mu\nu}-\epsilon h_{\mu\nu}+\frac{\epsilon^2}{2}h_{\mu\alpha}h_{\alpha\nu}+O(\epsilon^3)$, $\mu,\nu,\alpha=1,\cdots,n$.
\item $\det(g)=\exp(\epsilon \trace h)$.
\end{enumerate}

The following two lemmas are similar with \cite[Lemma 2.1 and Lemma 2.2]{ambrosetti1999multiplicity}. Readers can also see from \cite[Proposition 11]{brendle2007convergence}.
\begin{lem}
$$
R_{g}=\epsilon R_1+\epsilon^2 R_2+O(\epsilon^3),
$$
where 
\begin{align*}
R_1:=&D_{\mu\nu}h_{\mu\nu}-D_{\mu\mu}h_{\nu\nu};\\
R_2:=&h_{\mu\nu}(D_{\mu\nu}h_{\alpha\alpha}-D_{\mu\alpha}h_{\alpha\nu})\\
&+D_{\mu}h_{\mu\nu}D_{\nu}h_{\alpha\alpha}-\frac{1}{2}D_{\mu}h_{\mu\nu}D_{\alpha}h_{\alpha\nu}-\frac{1}{4}D_{\mu}h_{\nu\nu}D_{\mu}h_{\alpha\alpha}-\frac{1}{4}(D_{\mu}h_{\nu\alpha})^2.
\end{align*}
\end{lem}

\begin{proof}
Follow the definitions:
\begin{align*}
\Gamma^{\alpha}_{\mu\nu}&=\frac{1}{2}(D_{\mu}g_{\beta\nu}+D_{\nu}g_{\beta\mu}-D_{\beta}g_{\mu\nu})g^{\beta\alpha};\\
R^{\alpha}_{\beta\mu\nu}&=D_{\mu}\Gamma^{\alpha}_{\nu\beta}-D_{\nu}\Gamma^{\alpha}_{\mu\beta}+\Gamma^{\alpha}_{\mu \gamma}\Gamma^{\gamma}_{\nu\beta}-\Gamma^{\alpha}_{\nu \gamma}\Gamma^{\gamma}_{\mu\beta};\\
R_{g}&=g^{\mu\nu}R^{\alpha}_{\mu\alpha\nu}.
\end{align*}
The computation is not related to the main topic. Here we just write out the computation for $\Gamma^{\alpha}_{\mu\nu}$:
\begin{align*}
\Gamma^{\alpha}_{\mu\nu}=&\frac{\epsilon}{2} (D_{\mu}h_{\alpha\nu}+D_{\nu}h_{\alpha\mu}-D_{\alpha}h_{\mu\nu})\\
+&\frac{\epsilon^2}{4}(D_{\mu}(h^2)_{\alpha\nu}+D_{\nu}(h^2)_{\alpha\mu}-D_{\alpha}(h^2)_{\mu\nu}-2h_{\alpha\beta}(D_{\mu}h_{\beta\nu}+D_{\nu}h_{\beta\mu}-D_{\beta}h_{\mu\nu}))\\
+&O(\epsilon^3).
\end{align*} 
Then we can obtain the conclusion.
\end{proof}
Consider 
\begin{align*}
I_{g}(u)&=\frac{1}{2}\int_{\r^n}|D_{g}u|^2+c(n)R_{g}u^2dV_g-\frac{(n-2)^2}{2}\int_{\r^n}u^{\frac{2n}{n-2}}dV_g;\\
I_{\delta}(u)&=\frac{1}{2}\int_{\r^n}|Du|^2dx-\frac{(n-2)^2}{2}\int_{\r^n}|u|^{\frac{2n}{n-2}}dx.
\end{align*}
\begin{lem}\label{lem5}
$$I_{g}(u)=I_{\delta}(u)+\epsilon G_1(u)+\epsilon^2 G_2(u)+O(\epsilon^3),$$
where
\begin{align*}
G_1(u)&=\frac{1}{2}\int -h_{\mu\nu}D_{\mu}u D_{\nu}u+c(n)R_1u^2+\frac{\trace h}{2}\left(\frac{1}{2}|Du|^2-\frac{(n-2)^2}{2}|u|^{\frac{2n}{n-2}}\right)\\
G_2(u)&=\frac{1}{2}\int \frac{1}{2}(h^2)_{\mu\nu}D_{\mu}u D_{\nu}u+c(n)R_2u^2+\frac{(\trace h)^2}{8}\left(\frac{1}{2}|Du|^2-\frac{(n-2)^2}{2}|u|^{\frac{2n}{n-2}}\right)\\
&+\frac{\trace h}{2}(-h_{\mu\nu}D_{\mu}u D_{\nu}u+c(n)R_1u).
\end{align*}
\end{lem}
\begin{proof}
We can plug the following computations into the definition of $I_{g}(u)$ and obtain the results.
\begin{align*}
R_{g}&= \epsilon R_1+ \epsilon^2 R_2+O(\epsilon^3);\\
|D_{g}u|^2&=|Du|^2-\epsilon h_{\mu\nu}D_{\mu}uD_{\nu}u+\frac{\epsilon^2}{2}(h^2)_{\mu\nu}D_{\mu}uD_{\nu}u+O(\epsilon^3);\\
dV_g&=\left(1+\epsilon\trace h+\frac{\epsilon^2}{2}(\trace h)^2+O(\epsilon^3)\right)dx.
\end{align*}
\end{proof}

\bibliographystyle{plain}
\bibliography{Ref}

\Addresses

\end{document}